\documentclass[11pt,a4paper,english]{article}
\usepackage{mathptmx}
\usepackage[T1]{fontenc}
\usepackage[latin9]{inputenc}
\usepackage{color}
\usepackage{array}
\usepackage{amsmath}
\usepackage{amssymb}
\usepackage{graphicx}
\usepackage{esint}

\makeatletter

\pdfpageheight\paperheight
\pdfpagewidth\paperwidth

\usepackage{mathptmx,latexsym}
\usepackage[labelformat=simple]{subfig}

\usepackage{iamPreprintStyle}

\makeatother

\usepackage{babel}

\begin{document}
\preprintTitle{Time-space adaptive discontinuous Galerkin method for advection-diffusion equations with non-linear reaction mechanism} \preprintAuthor{B\"{u}lent Karas\"{o}zen%
\footnote{Department of Mathematics and Institute of Applied Mathematics, Middle East Technical University, 06800 Ankara, Turkey, \textit{email}: bulent@metu.edu.tr%
}, Murat Uzunca%
\footnote{Department of Mathematics, Middle East Technical University, 06800 Ankara, Turkey, \textit{email}: uzunca@metu.edu.tr%
}} 

\preprintAbstract{\small In this work, we apply a time-space adaptive discontinuous Galerkin method using  the elliptic reconstruction technique with a robust (in P\'{e}clet number) elliptic error estimator in space, for the convection dominated parabolic problems with non-linear reaction mechanisms.  We derive a posteriori error estimators in the $L^{\infty}(L^2)+L^2(H^1)$-type norm using backward Euler in time and discontinuous Galerkin (symmetric interior penalty Galerkin (SIPG)) in space. Numerical results for advection dominated reactive transport problems in homogeneous and heterogeneous media demonstrate the performance of the time-space adaptive algorithm.} 

\preprintKeywords{Non-linear diffusion-convection reaction, Discontinuous Galerkin, Time-space adaptivity, Elliptic reconstruction, A posteriori error estimates} 

\preprintDate{August 2014} 
\preprintNo{2014-6} 
\makePreprintCoverPage


\section{Introduction}
\label{intro}

Advection-diffusion-reaction (ADR) equations are widely used for systems in chemical reaction problems. In the linear convection or advection dominated case, stabilized continuous finite elements and discontinuous Galerkin (DG) methods are capable of handling the nonphysical oscillations. On the other hand, in the non-linear stationary case, the  non-linear reaction term produces sharp layers in addition to the spurious oscillations due to the convection. Thus, an accurate and efficient numerical resolution of such layers is a challenge as the exact location of the layers are not known a priori. In the non-stationary case, the resolution of such layers is more critical since the nature of the sharp layers may vary as time progresses. Recently, several stabilization and shock/discontinuity capturing techniques were developed for linear and non-linear steady-state problems \cite{bause12ash}. In contrast to the stabilized  continuous Galerkin finite element methods, DG methods produce stable discretizations without the need for stabilization parameters. The DG combines the best properties of the finite volume and continuous finite elements methods. Finite volume methods can only use lower degree polynomials, and continuous finite elements methods require higher regularity due to the continuity requirements. The DG method is in particular suitable for non-matching grids and hp (space and order) adaptivity, detecting sharp layers and singularities. They are easily adapted locally for nonconforming finite elements requiring less regularity.  Higher order DG approximation can be easily used by hierarchical bases \cite{burger09pkd}, and the convergence rates are limited by the consistency error which makes the DG suitable for complex fluid flows. The DG methods are robust with respect to the variation of the physical parameters like diffusion constant and permeability. The stability of the DG approximation retained by adjusting the penalty parameter to penalize the jumps at the interface of the elements. Various choices of the penalty parameter is suggested in the literature \cite{castillo12pdg,dobrev08psi,epshteyn07epp,slingerland14fls}. A unified analysis of the interior penalty DG methods for elliptic PDEs are given in \cite{arnold02uad}. Other advantages of the DG methods are conservation of mass and fluxes and parallelization. Moreover, DG methods are  better suited for adaptive strategies which are based on (residual-based) a posteriori error estimation. The disadvantages are the resulting larger and denser matrices and ill-conditioning with increasing degree of the discontinuous polynomials.

The main tool in this paper is the adaptivity applied to the DG discretized (in space) systems using (residual-based) a posteriori error estimates. The aim of the a posteriori estimates is to derive bounds between the known numerical solution and the unknown exact solution. For elliptic diffusion-convection-reaction equations, there are various studies applying adaptivity using a posteriori error estimates in the literature \cite{pietro14rra,houston02dhp,klieber06ast,schotzau09rae}, which are well-understood. In this paper, we derive a posteriori error estimates for semi-linear ADR problems using the well-know elliptic a posteriori estimates, by which in contrast to the standard energy techniques we do not need to try to adapt the estimates case by case in order to compare the exact solution with numerical solution directly. For this reason, we use the {\it elliptic reconstruction} technique in \cite{makridakis03era} which allows us to utilize available a posteriori estimates derived for elliptic equations to control the main part of the spatial error. The idea of the elliptic reconstruction technique is to construct an auxiliary solution whose difference to the numerical solution can be estimated by a known (elliptic) a posteriori estimate, and the constructed auxiliary solution satisfies a variant of the given problem with a right hand side which can be controlled in an optimal way. By this way, we are able to obtain results being optimal order in both $L^2(H^1)$ and $L^{\infty}(L^2)$-type norms, while the results obtained by the standard energy methods are only optimal order in $L^2(H^1)$-type norms, but sub-optimal order in $L^{\infty}(L^2)$-type norms. In \cite{cangiani13adg}, a posteriori error estimates in the $L^{\infty}(L^2)+L^2(H^1)$-type norm are derived for linear parabolic diffusion-convection-reaction equations using backward Euler in time and discontinuous Galerkin in space utilizing the elliptic reconstruction technique. In this paper, we extend the study in \cite{cangiani13adg} in a similar way by deriving and implementing a posteriori error estimates in the $L^{\infty}(L^2)+L^2(H^1)$-type norm using backward Euler in time and discontinuous Galerkin (symmetric interior penalty (SIPG)) in space for the convection dominated parabolic problems with non-linear reaction mechanisms. To derive the a posteriori error estimates, we use the modification of the robust (in P\'{e}clet number) a posteriori error estimator in \cite{schotzau09rae} for linear steady-state diffusion-convection-reaction equations to the steady-state diffusion-convection equations with non-linear reaction mechanisms utilizing the elliptic reconstruction technique \cite{makridakis03era}.

Application of the adaptive discontinuous Galerkin methods and a posteriori error estimates to the problems in geoscience are reviewed recently in \cite{pietro14rra}. Most of the applications of DG methods in geoscience concern reactive transport with advection \cite{bastian11udg,klieber06ast,sun05lna} and strong permeability contrasts such as layered reservoirs \cite{slingerland14fls} or vanishing and varying  diffusivity posing challenges in computations \cite{proft09dgm}. The permeability in heterogeneous porous and fractured media varies over orders of magnitude in space, which results in highly variable flow field, where the local transport is dominated by advection or diffusion \cite{tambue10eia}. Accurate and efficient numerical solution of the ADR equations to predict the macroscopic mixing, anomalous transport of the solutes and contaminants  for a wide range of parameters like permeability and P\'{e}clet numbers, different flow velocities and reaction rates and reaction rates are challenging problems \cite{tambue10eia}. In order to resolve the complex flow patterns accurately, higher order time stepping methods like exponential time stepping methods are used \cite{tambue10eia}. We show here using time-space adaptive first order backward Euler and DG in space, the same results can be obtained.

The rest of this paper is organized as follows. In the next section, we introduce the function spaces and related norms which are used in our analysis, and the model problem with the assumptions to have unique solution. In Section \ref{SIPG}, we give the symmetric discontinuous interior penalty Galerkin discretized semi-discrete system, and fully-discrete system using Backward Euler in time. The derivation of a posteriori error estimator is given in Section~\ref{aposteriori}, first for steady-state problems (with proofs), then for parabolic problems utilizing the elliptic reconstruction technique. In Section~\ref{algorithm}, we state the adaptive algorithm procedure using the a posteriori error estimators derived in Section~\ref{aposteriori}. The solution procedure of the fully-discrete system by Newton method and the structures of the arising matrix and vectors are discussed in Section~\ref{newton}. After demonstrating the performance of the algorithm by some numerical studies in Section~\ref{numeric}, the paper ends with conclusions.

\section{Model problem}
\label{secmodel}

Let $\Omega \subset\mathbb{R}^2$ be a bounded, open and convex domain with boundary $\partial\Omega$. For a Banach space $X$, define the spaces $L^p(0,T;X)$
\begin{subequations}\label{model}
\begin{align}
  \| v\|_{L^p(0,T;X)} &= \left( \int_0^T \| v(t)\|_X^p dt \right)^{1/p}<\infty \; , \qquad \text{for } 1\leq p < +\infty \nonumber \\
  \| v\|_{L^{\infty}(0,T;X)} &= \underset{0\leq t\leq T}{\text{esssup}} \|v(t)\|_X <\infty \; , \qquad \text{for } p = +\infty \nonumber
\end{align}
\end{subequations}
Also define the space
$$
H^1(0,T;X)=\{ v\in L^2(0,T;X)| \; v_t\in L^2(0,T;X)\}.
$$
We denote by $C(0,T;X)$ and $C^{0,1}(0,T;X)$, the spaces of continuous and Lipschitz-continuous functions $v:[0,T]\mapsto X$, respectively, equipped with the norms
\begin{subequations}
\begin{align}
  \| v\|_{C(0,T;X)} &= \underset{0\leq t\leq T}{\text{max}} \|v(t)\|_X <\infty \nonumber \\
  \| v\|_{C^{0,1}(0,T;X)} &= \text{max} \left\{ \| v\|_{C(0,T;X)}, \| v_t\|_{C(0,T;X)}\right\} <\infty \nonumber
\end{align}
\end{subequations}
We consider the system of semi-linear diffusion-convection-reaction equations
\begin{equation}\label{org}
  \frac{\partial u_i}{\partial t}-\nabla\cdot (\epsilon_i\nabla u_i) + \vec{\beta}_i\cdot\nabla u_i + r_i(\vec{u}) = f_i \;\; , \quad i=1,2,\ldots , J
\end{equation}
in $\Omega\times (0,T]$ for the vector of unknowns $\vec{u}=(u_1,u_2,\ldots , u_J)^T$ with appropriate boundary and initial conditions. We assume that the source functions $f_i\in C(0,T;L^2(\Omega))$, and the velocity fields $\vec{\beta }_i\in C\left(0,T;W^{1,\infty}(\Omega)\right)^2$ either given or computed. For flow in
heterogeneous media in Section~\ref{ex4}, the symmetric dispersion tensors $\epsilon_i$ are taken of the form
$$
\epsilon_i=
\begin{bmatrix}
D_i^1 & 0 \\
0 & D_i^2
\end{bmatrix}
$$
with $0<D_i^1,D_i^2\ll 1$. Moreover, we assume that the non-linear reaction terms are bounded, locally Lipschitz continuous and monotone, i.e. satisfies for any $s, s_1, s_2\ge 0$, $s,s_1, s_2 \in \mathbb{R}$ the following conditions  \cite{uzunca14adg}
\begin{subequations}
\begin{align}
|r_i(s)| &\leq C_s , \quad C_s>0 \label{nonl1}\\
\| r_i(s_1)-r_i(s_2)\|_{L^2(\Omega)} &\leq L\| s_1-s_2\|_{L^2(\Omega)} , \quad L>0\label{nonl2} \\
 r_i\in C^1(\mathbb{R}_0^+), \quad r_i(0) =0, &\quad r_i'(s)\ge 0. \label{nonl3}
\end{align}
\end{subequations}
We also assume that there are $\kappa , \kappa_*\geq 0$ satisfying for $i=1,2,\ldots ,J$
\begin{equation} \label{ass}
-\frac{1}{2}\nabla\cdot\vec{\beta }_i (x) \geq \kappa, \qquad \| -\nabla\cdot\vec{\beta }_i\|_{C(0,T;L^{\infty}(\Omega))}\leq \kappa^*\kappa,
\end{equation}
The first inequality in (\ref{ass}) ensures the well-posedness of the problem, and the letter is used in the a posteriori error analysis.

The weak formulation of the system (\ref{org}) reads as: for any $v\in H_0^1(\Omega )$, find $u_i\in L^2(0,T;H_0^1(\Omega))\cap H^1(0,T;L^2(\Omega))$, $i=1,2,\ldots , J$, such that for all $t\in (0,T]$

\begin{equation} \label{weak}
\int_{\Omega}\frac{\partial u_i}{\partial t}vdx + a(t;u_i,v)+b_i(t;\vec{u},v)=l_i(v)
\end{equation}
\begin{subequations}
\begin{align}
a(t;u_i, v)=& \int_{\Omega}(\epsilon_i\nabla u_i\cdot\nabla v+\vec{\beta }_i\cdot\nabla u_iv)dx , \label{weaka} \\
b_i(t;\vec{u}, v) =& \int_{\Omega}r_i(\vec{u})v dx ,   \\
l_i( v)=&  \int_{\Omega}f_iv dx
\end{align}
\end{subequations}
which have unique solutions $\{ u_i\}_{i=1}^J$ in the space $C(0,T;L^2(\Omega))$ under the given regularity assumptions and the conditions \eqref{nonl1}-\eqref{nonl3}.

In the sequel, for simplicity, we just consider a single equation of the system \eqref{org} (J=1) without any subscript to construct the discontinuous Galerkin discretization in the next section and thereafter, and we also continue with a homogeneous dispersion tensor leading to a simple diffusivity constant $0<\epsilon\ll 1$. We, furthermore, take into account the homogeneous Dirichlet boundary conditions to simplify the notations. It can be proceeded with heterogeneous dispersion tensor and other type of boundary conditions in a standard way.

\section{Discontinuous Galerkin discretization}
\label{SIPG}

For space discretization of  a single equation of (\ref{org}), we use the symmetric discontinuous interior penalty Galerkin (SIPG) method \cite{arnold02uad,riviere08dgm} with  the upwinding  for the convection part \cite{ayuso09dgm,houston02dhp}.

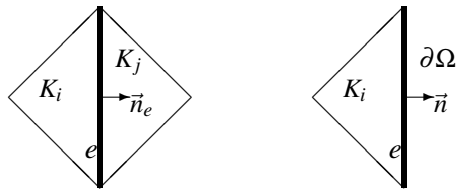
\begin{figure}[ht]
\centering
\setlength{\unitlength}{2.0mm}
\begin{picture}(35, 15)
\put(5,7){\line(1,1){6}}
\put(5,7){\line(1,-1){6}}
\put(17,7){\line(-1,1){6}}
\put(17,7){\line(-1,-1){6}}
\put(7,7){\small$K_{i}$ }
\put(12,9){\small$K_{j}$ }
\put(10,3){$e$ }
\put(11,7){\vector(1,0){2}}
\put(13,6){\small$\vec{n}_{e}$}

\put(25,7){\line(1,1){6}}
\put(25,7){\line(1,-1){6}}
\put(27,7){\small$K_{i}$ }
\put(32,9){\small$\partial\Omega$ }
\put(30,3){$e$ }
\put(31,7){\vector(1,0){2}}
\put(33,6){\small$\vec{n}$}

\linethickness{0.7mm} \put(11,1){\line(0,1){12}}
\linethickness{0.7mm} \put(31,1){\line(0,1){12}}
\end{picture}
\caption{Two adjacent elements sharing an edge (left); an element near to domain boundary (right)}
\label{jump}
\end{figure}

Let $\{\xi_h\}$ be a family of shape regular meshes with the elements (triangles) $K_i\in\xi_h$ satisfying $\overline{\Omega}=\cup \overline{K}$ and $K_i\cap K_j=\emptyset$ for $K_i$, $K_j$ $\in\xi_h$. Let us denote by $\Gamma^0$ and $\Gamma^{\partial}$ the set of interior and Dirichlet boundary edges, respectively, so that $\Gamma =\Gamma^0\cup\Gamma^{\partial}$ forms the skeleton of the mesh. For any $K\in\xi_h$, let $\mathbb{P}_k(K)$ be the set of all polynomials of degree at most $k$ on $K$. Then, set the finite dimensional space
$$
V_h(\xi_h)=\left\{ v\in L^2(\Omega ) : v|_{K}\in\mathbb{P}_k(K) ,\; \forall K\in \xi_h \right\}\not\subset H_0^1(\Omega).
$$
The functions in $V_h(\xi_h)$, in contrast to the standard (continuous) finite elements, are discontinuous along the inter-element boundaries causing that along an interior edge, there are two different traces from the adjacent elements sharing that edge. In the light of this fact, let us first introduce some notations before giving the SIPG formulation. Let $K_i$, $K_j\in\xi_h$ ($i<j$) be two adjacent elements sharing an interior edge $e=K_i\cap K_j\subset \Gamma_0$ (see Fig.\ref{jump}). We denote the trace of a scalar function $v$ from inside $K_i$ by $v_{i}$ and from inside $K_j$ by $v_{j}$, and  then we set the jump and average values of $v$ on the edge $e$
$$
[v]= v_{i}\vec{n}_e- v_{j}\vec{n}_e , \quad \{ v\}=\frac{1}{2}(v_{i}+ v_{j}).
$$
Here $\vec{n}_e$ denotes  the unit normal to the edge $e$ oriented from $K_i$ to $K_j$. Similarly, we set the jump and average values of a vector valued function $\vec{q}$ on e
$$
[\vec{q}]= \vec{q}_{i}\cdot \vec{n}_e- \vec{q}_{j}\cdot \vec{n}_e , \quad \{ \vec{q}\}=\frac{1}{2}(\vec{q}_{i}+ \vec{q}_{j}),
$$
Observe that $[v]$ is a vector for a scalar function $v$, while, $[\vec{q}]$ is scalar for a vector valued function $\vec{q}$. On the other hand, along any boundary edge $e=K_i\cap \partial\Omega$, we set
$$
[v]= v_{i}\vec{n} , \quad \{ v\}=v_{i}, \quad [\vec{q}]=\vec{q}_{i}\cdot \vec{n}, \quad \{ \vec{q}\}=\vec{q}_{i}
$$
where $\vec{n}$ is the unit outward normal to the boundary at $e$. We define the inflow and outflow boundary parts at a time $t$ by
$$
\Gamma_t^-=\{ x\in \partial\Omega | \; \vec{\beta}(\vec{x},t)\cdot \vec{n}(x)<0\}\; , \qquad \Gamma_t^+=\partial\Omega\setminus\Gamma_t^-
$$
At a time $t$, the inflow and outflow parts of an element $K$ are defined by
$$
\partial K_t^-=\{ x\in \partial K | \; \vec{\beta}(\vec{x},t)\cdot \vec{n}_K(x)<0\}\; , \qquad \partial K_t^+=\partial K\setminus\partial K_t^-
$$
where $\vec{n}_K(x)$ denotes the outward unit normal to the boundary of the element $K$ at $x$.

Under the given definitions, the semi-discrete problem reads as: for $t=0$, set $u_h(0)\in V_h(\xi_h)$ as the projection (orthogonal $L^2$-projection) of $u_0$ onto $V_h(\xi_h)$; for each $t\in (0,T]$, for all $v_h\in V_h(\xi_h)$, find $u_h\in C^{0,1}(0,T;V_h(\xi_h))$ such that

\begin{equation} \label{dg}
\int_{\Omega}\frac{\partial u_{h}}{\partial t}v_{h}dx + a_{h}(t;u_{h},v_{h}) + K_{h}(u_{h},v_{h})+b_{h}(t;u_{h}, v_{h})=l_{h}(v_{h})  ,
\end{equation}
\begin{subequations}
\begin{align}
a_{h}(t;u_{h}, v_{h})=& \sum \limits_{K \in {\xi}_{h}} \int_{K} \epsilon \nabla u_{h}\cdot\nabla v_{h} dx + \sum \limits_{K \in {\xi}_{h}} \int_{K} \vec{\beta } \cdot \nabla u_{h} v_{h} dx \label{dga} \\
&+ \sum \limits_{K \in {\xi}_{h}}\int_{\partial K_t^-\setminus\partial\Omega } \vec{\beta }\cdot \vec{n}_K (u_{h}^{out}-u_{h})  v_{h} ds - \sum \limits_{K \in {\xi}_{h}} \int_{\partial K_t^-\cap \Gamma_t^{-}} \vec{\beta }\cdot \vec{n}_K u_{h} v_{h} ds  \nonumber  \\
&+ \sum \limits_{ e \in \Gamma}\frac{\sigma \epsilon}{h_{e}} \int_{e} [u_{h}]\cdot[v_{h}] ds, \nonumber  \\
K_{h}(u_{h},v_{h}) =& - \sum \limits_{ e \in \Gamma} \int_{e} ( \{\epsilon \nabla v_{h} \}\cdot[u_{h}] + \{\epsilon \nabla u_{h} \}\cdot [v_{h}] )ds \label{kh}  \\
b_{h}(t;u_{h}, v_{h}) =& \sum \limits_{K \in {\xi}_{h}} \int_{K} r(u_{h}) v_{h} dx, \\
l_{h}( v_{h})=&  \sum \limits_{K \in {\xi}_{h}} \int_{K} f_h v_{h} dx
\end{align}
\end{subequations}
where $u_{h}^{out}$ denotes the value on an edge from outside of an element $K$. The parameter $\sigma\in\mathbb{R}_0^+$ is called the penalty parameter which should be sufficiently large; independent of the mesh size $h$ and the diffusion coefficient $\epsilon$ \cite{riviere08dgm} [Sec. 2.7.1]. We choose the penalty parameter $\sigma$ for the SIPG method depending on the polynomial degree $k$ as
$\sigma=3k(k+1)$ on interior edges and $\sigma=6k(k+1)$ on boundary edges.

Upon integration by parts on the convective term, bilinear form (\ref{dga}) will be
\begin{subequations}
\begin{align}
a_{h}(t;u_{h}, v_{h})=& \sum \limits_{K \in {\xi}_{h}} \int_{K} \epsilon \nabla u_{h}\cdot\nabla v_{h} dx - \sum \limits_{K \in {\xi}_{h}} \int_{K} ( \vec{\beta } u_h\cdot\nabla v_h + \nabla\cdot\vec{\beta } u_hv_h) dx \label{bilbyp} \\
&+ \sum \limits_{K \in {\xi}_{h}}\int_{\partial K_t^+\setminus\partial\Omega } \vec{\beta }\cdot \vec{n}_K u_{h}(v_h-v_{h}^{out}) ds + \sum \limits_{K \in {\xi}_{h}} \int_{\partial K_t^+\cap \Gamma_t^{+}} \vec{\beta }\cdot \vec{n}_K u_{h} v_{h} ds  \nonumber  \\
&+ \sum \limits_{ e \in \Gamma}\frac{\sigma \epsilon}{h_{e}} \int_{e} [u_{h}]\cdot[v_{h}] ds \nonumber
\end{align}
\end{subequations}
Note that the bilinear form $a_{h}(t;u,v)$ is well-defined for the functions $u,v\in H_0^1(\Omega )$, and equal to
$$
a_{h}(t;u, v)=\int_{\Omega}(\epsilon\nabla u\cdot\nabla v + \vec{\beta}\cdot\nabla uv)dx
$$
Thus, the weak formulation (\ref{weaka}) can be rewritten for any $t\in (0,T]$ as
\begin{equation}\label{eqweak}
\int_{\Omega}\frac{\partial u}{\partial t}vdx + a_h(t;u,v)+b(t;u,v)=l(v) \; , \qquad \forall v\in H_0^1(\Omega )
\end{equation}

For the fully discrete system, consider a subdivision of $[0,T]$ into $n$ time intervals $I_k=(t^{k-1},t^{k}]$ of length $\tau_k$, $k=1,2,\ldots , n$. Set $t^0=0$ and for $n\geq 1$, $t^k=\tau_1+\tau_2+\cdots + \tau_k$. Denote by $\xi_h^0$ an initial mesh and by $\xi_h^k$ the mesh associated to the $k^{th}$ time step for $k>0$, which is obtained from $\xi_h^{k-1}$ by locally refining/coarsening. Moreover, we assign the finite element space $V_h^k=V_h(\xi_h^k)$ to each mesh $\xi_h^k$. Then, the fully discrete problem, backward Euler in time, reads as: for $t=0$, set $u_h^0\in V_h^0$ as the projection (orthogonal $L^2$-projection) of $u_0$ onto $V_h^0$; for $k=1,2,\ldots , n$, find $u_h^k\in V_h^{k}$ such that for all $v_h^k\in V_h^{k}$
\begin{equation}\label{fullydisc}
\int_{\Omega}\frac{u_{h}^{k}-u_{h}^{k-1}}{\tau_k}v_{h}^k dx + a_{h}(t^k;u_{h}^k,v_{h}^k)+ K_{h}(u_{h}^k,v_{h}^k)+b_{h}(t^k;u_{h}^k, v_{h}^k)=\int_{\Omega}f^{k}v_{h}^k dx
\end{equation}

\section{A posteriori error analysis}
\label{aposteriori}

The problems concerned in this paper are convection/reaction dominated non-stationary models which often produce internal/boundary layers where the solutions have steep gradients. Thus, an efficient solution of such equations are needed. The most popular technique to solve the convection/reaction dominated problems is the adaptivity which requires an estimator to locate the regions where the error is too large. In this work, we construct the residual-based, robust (in P\'{e}clet number) a posteriori error estimators for non-stationary convection domianted diffusion-convection equations with non-linear reaction mechanisms. To construct the a posteriori error estimators, we extend the a posteriori error estimators for linear non-stationary problems constructed in \cite{cangiani13adg} which uses the a posteriori error estimators for linear stationary models constructed in \cite{schotzau09rae} utilizing the elliptic reconstruction technique \cite{makridakis03era} to make connection between the stationary and non-stationary error. As in \cite{cangiani13adg}, first we construct and prove the a posteriori error bounds for non-linear stationary models in Section~\ref{stationary}. Then, in Section~\ref{semidiscerror}, we give the a posteriori error bounds for the semi-discrete system of the non-stationary problems with non-linear reaction term. Finally, in Section~\ref{fullydiscerror}, we give the a posteriori error bounds for the fully-discrete system of the non-stationary problems with non-linear reaction term. The main contribution to the error analysis lies in the construction of the error bounds for non-linear stationary models in Section~\ref{stationary}. The remaining utilize the constructions in \cite{cangiani13adg}.

\subsection{A posteriori error bounds for non-linear stationary system}
\label{stationary}

We consider first  the stationary problem, i.e. for a given $t\in (0,T]$, let $u^s \in H_0^1(\Omega )$ be the unique solution to the weak problem
\begin{equation} \label{stweak}
a(t;u^s,v)+b(t;u^s,v)=l(v) \; , \qquad \forall v\in H_0^1(\Omega )
\end{equation}
and let $u_h^s\in V_h(\xi_h)$ be the unique solution to the discrete problem
\begin{equation} \label{stdg}
a_h(t;u_h^s,v_h)+K_{h}(u_{h}^s,v_{h})+b_h(t;u_h^s,v_h)=l(v_h) \; , \qquad \forall v\in V_h(\xi_h)
\end{equation}
In order to measure the spatial error, we use the energy norm

$$
||| v |||^2=\sum \limits_{K \in {\xi}_{h}}(\| \epsilon\nabla v\|_{L^2(K)}^2 +\kappa \| v\|_{L^2(K)}^2) + \sum \limits_{e \in \Gamma}\frac{\epsilon\sigma}{h_{e}}\| [v]\|_{L^2(e)}^2,
$$
and the semi-norm
\begin{equation} \label{semin}
|v|_C^2=|\vec{\beta }v|_*^2+\sum\limits_{e \in \Gamma}(\kappa h_e+ \frac{h_e}{\epsilon})\| [v]\|_{L^2(e)}^2,
\end{equation}
where
$$
|u|_*=\mathop{\text{sup}}_{w\in H_0^1(\Omega )\setminus \{ 0\}}\frac{\int_{\Omega}u\cdot \nabla w dx}{||| w|||},
$$
For each element $K \in {\xi}_{h}$, we define the local error indicators $\eta_K^2$ as

\begin{eqnarray} \label{res}
\eta_K^2= \eta_{R_K}^2  + \eta_{E_K}^2 + \eta_{J_K}^2,
\end{eqnarray}
where
\begin{eqnarray}
\eta_{R_K}^2&=&\rho_K^2\| f_h +\epsilon\Delta u_h^s-\vec{\beta }_h\cdot\nabla u_h^s-r(u_h^s)\|_{L^2(K)}^2, \nonumber \\
\eta_{E_K}^2 &=& \frac{1}{2} \sum \limits_{e \in \partial K\setminus\Gamma^{\partial} }\epsilon^{-\frac{1}{2}}\rho_e\| [\epsilon\nabla u_h^s]\|_{L^2(e)}^2, \nonumber \\
\eta_{J_K}^2 &=& \frac{1}{2}\sum \limits_{e \in \partial K\setminus\Gamma^{\partial} }\left( \frac{\epsilon\sigma}{h_e}+\kappa h_e+\frac{h_e}{\epsilon}\right) \| [u_h^s]\|_{L^2(e)}^2 +\sum \limits_{e \in \partial K\cap\Gamma^{\partial}}\left( \frac{\epsilon\sigma}{h_e}+\kappa h_e+\frac{h_e}{\epsilon}\right) \| u_h^s\|_{L^2(e)}^2, \nonumber \\
\end{eqnarray}
with the weights $\rho_K$ and $\rho_e$, on an element $K$, given by
$$
\rho_{K}=\min\{h_{K}\epsilon^{-\frac{1}{2}}, \kappa^{-\frac{1}{2}}\}, \; \rho_{e}=\min\{h_{e}\epsilon^{-\frac{1}{2}}, \kappa^{-\frac{1}{2}}\},
$$
for $\kappa \neq 0$. When $\kappa =0$, we take $\rho_{K}=h_{K}\epsilon^{-\frac{1}{2}}$ and $\rho_{e}=h_{e}\epsilon^{-\frac{1}{2}}$. Then, our a posteriori error indicator is given by
\begin{equation} \label{errind}
\eta=\left( \sum \limits_{K\in{\xi}_{h}}\eta_K^2\right)^{1/2}.
\end{equation}
We also introduce the data approximation terms,
$$
\Theta_K^2 =\rho_K^2(\| f-f_h\|_{L^2(K)}^2 + \| (\vec{\beta } -\vec{\beta }_h)\cdot\nabla u_h^s\|_{L^2(K)}^2).
$$
and  the data approximation error
\begin{equation} \label{data}
\Theta=\left( \sum \limits_{K\in{\xi}_{h}}\Theta_K^2\right)^{1/2}.
\end{equation}
Then, we have the a posteriori error bounds

\begin{eqnarray}
\| u^s-u_h^s\|_{DG} \lesssim \eta + \Theta  \qquad &(\text{reliability}), \label{rel} \\[0.2cm]
\eta \lesssim \| u^s-u_h^s\|_{DG} + \Theta \qquad &(\text{efficiency}). \label{eff}
\end{eqnarray}
with
$$
\| v\|_{DG}= ||| v||| + | v|_C
$$

\subsubsection{Proof of a posteriori error bounds}

The proof  for the bounds of spatial error estimates  are analogous to the ones in \cite{schotzau09rae} for the linear problems. Therefore, only the proofs for  non-linear reaction term are stated explicitly . In the following, we use the symbols $\lesssim$ and $\gtrsim$ to denote the bounds that are valid up to positive constants independent of the local mesh size $h$, the diffusion coefficient $\epsilon$ and the penalty parameter $\sigma$.

The spatial error $\| u^s-u_h^s\|_{DG}$ is not well-defined, since $u^s\in H_0^1(\Omega)$ and $u_h^s\in V_h(\xi_h)\nsubseteq H_0^1(\Omega)$. Therefore, we split the stationary SIPG solution $u_h^s$ as

$$
u_h^s=u_h^c+u_h^r
$$
with $u_h^c\in H_0^1(\Omega)\cap V_h(\xi_h)$ being the conforming part of the solution and $u_h^r\in V_h(\xi_h)$ is the remainder term. In this way, we have $u_h^s\in H_0^1(\Omega)+V_h(\xi_h)$, and from the triangular inequality
$$
\| u^s-u_h^s\|_{DG}\leq \| u^s-u_h^c\|_{DG}+\| u_h^r\|_{DG}
$$
Now, both the terms on the right hand side are well-defined norms, and our aim is to find bounds for them. Next, we introduce the following auxiliary forms:

\begin{subequations}
\begin{align}
D_{h}(u, v)=& \sum \limits_{K \in {\xi}_{h}} \int_{K} \left( \epsilon \nabla u\cdot\nabla v -\nabla\cdot \vec{\beta }uv\right) dx \label{d} \\
O_h(u,v)=& -\sum \limits_{K \in {\xi}_{h}} \int_{K}\vec{\beta }u\cdot\nabla v dx +\sum \limits_{K \in {\xi}_{h}} \int_{\partial K^+\cap\Gamma^+}\vec{\beta }\cdot\vec{n}_K uv ds \nonumber \\
& +\sum \limits_{K \in {\xi}_{h}} \int_{\partial K^+\setminus\partial\Omega }\vec{\beta }\cdot\vec{n}_Ku(v-v^{out}) ds \label{o} \\
J_h(u,v) =& \sum \limits_{ e \in \Gamma_{0}\cup\Gamma }\frac{\sigma \epsilon}{h_{e}} \int_{e} [u]\cdot [v] ds. \label{j}
\end{align}
\end{subequations}
We note that, for a specific $t\in (0,T]$, the SIPG bilinear form (\ref{bilbyp}) satisfies
$$
a_h(t;u,v)=D_h(u,v)+O_h(u,v)+J_h(u,v)
$$
and is well-defined on $H_0^1(\Omega )+V_h(\xi_h)$. Using the first identity in (\ref{ass}), we can easily have for any $u\in H_0^1(\Omega)$
$$
{a}_h(t;u,v) \geq ||| u|||^2
$$
Moreover, the auxiliary forms are continuous \cite{schotzau09rae}[Lemma 4.2]:
\begin{eqnarray}
|D_h(u,v)| &\lesssim& |||u|||\; |||v||| \; , \qquad u,v\in H_0^1(\Omega )+V_h(\xi_h) ,\label{cd} \\
|O_h(u,v)| &\lesssim& |\vec{\beta }u|_* \;|||v||| \; ,  \qquad u\in H_0^1(\Omega )+V_h(\xi_h), v\in H_0^1(\Omega) ,\label{co} \\
|J_h(u,v)| &\lesssim& |||u|||\; |||v||| \; ,  \qquad u,v\in H_0^1(\Omega )+V_h(\xi_h) \label{cj} ,
\end{eqnarray}
and for $u\in V_h(\xi_h)$, $v\in V_h(\xi_h)\cap H_0^1(\Omega )$ \cite{schotzau09rae}[Lemma 4.3]
\begin{eqnarray}
|K_h(u,v)| &\lesssim& \gamma^{-1/2}\left( \sum_{e\in \Gamma_0\cup\Gamma^D} \frac{\gamma\epsilon}{h_e}\| [u]\|_{L^2(e)}\right)^{1/2} |||v|||.\label{ck}
\end{eqnarray}
We also have for the non-linear form $b_h(t;u,v)$, for a specific time $t$, using the boundedness assumption (\ref{nonl1}),
\begin{eqnarray}
|b_h(t;u,v)| &\lesssim& |||v||| \; ,  \qquad u,v\in H_0^1(\Omega )+V_h(\xi_h). \label{cb}
\end{eqnarray}
Now, we give some auxiliary results and conditions which are used in the proofs.
\begin{itemize}

\item {\bf Inf-sup condition:} \cite{schotzau09rae}[Lemma 4.4] For a nonzero $u\in H_0^1(\Omega )$, for a constant $C>0$, we have
\begin{eqnarray}
|||u|||+|\vec{\beta }u|_* \lesssim \underset{v\in H_0^1(\Omega ) \setminus \{ 0\}}{\text{sup}}\frac{{a}_h(t;u,v)}{|||v|||}. \label{infsup}
\end{eqnarray}

\item {\bf Approximation operator:} Let $V_h^c=V_h(\xi_h)\cap H_0^1(\Omega )$ be the conforming subspace of $V_h(\xi_h)$. For any $u\in V_h(\xi_h)$, there exists an approximation operator $A_h:\; V_h(\xi_h) \mapsto V_h^c$ satisfying
\begin{eqnarray}
\sum_{K\in\xi}\| u-A_hu\|_{L^2(K)}^2 \lesssim \sum_{e\in \Gamma_0\cup\Gamma^D}\int_e h_e|[u]|^2ds, \label{app}  \\
\sum_{K\in\xi}\| \nabla (u-A_hu)\|_{L^2(K)}^2 \lesssim \sum_{e\in \Gamma_0\cup\Gamma^D}\int_e \frac{1}{h_e}|[u]|^2ds. \label{dapp}
\end{eqnarray}

\item  {\bf Interpolation operator:} For any $u\in H_0^1(\Omega )$, there exists an interpolation operator
$$
I_h:\; H_0^1(\Omega ) \mapsto \{ w\in C(\overline{\Omega}): \; w|_K\in\mathbb{P}_1(K), \forall K\in\xi, w=0 \; \text{on }\Gamma \}
$$
that satisfies
\begin{eqnarray}
|||I_h u||| \lesssim |||u||| , \label{int1} \\
\left( \sum_{K\in\xi}\rho_K^{-2}\| u-I_hu\|_{L^2(K)}^2\right)^{1/2} \lesssim |||u|||, \label{int2}  \\
\left( \sum_{e\in\Gamma_0\cup\Gamma^D}\epsilon^{1/2}\rho_e^{-1}\| u-I_hu\|_{L^2(K)}^2\right)^{1/2} \lesssim |||u|||. \label{int3}
\end{eqnarray}
Now, consider the splitting of the stationary solution $u_h^s=u_h^c+u_h^r$ as $u_h^c=A_h u_h^s\in H_0^1(\Omega )\cap V_h(\xi_h)$ with $A_h$ is the approximation operator and $u_h^r=u_h^s-u_h^c\in V_h$.\\

\item  {\bf Bound for remainder term:} \cite{schotzau09rae}[Lemma 4.7] For the remainder term, we have the bound
\begin{eqnarray} \label{rem}
\| u_h^r\|_{DG} \lesssim \eta
\end{eqnarray}
where $\eta$ is our a posteriori error indicator (\ref{errind}).\\
\end{itemize}

\noindent {\bf Lemma:} For a given $t\in (0,T]$ and for any $v\in H_0^1(\Omega )$, we have
\begin{eqnarray}
\int_{\Omega } f(v-I_hv)dx-{a}_h(t;u_h^s,v-I_hv)-b_h(t;u_h^s,v-I_hv) \lesssim (\eta +\Theta )|||v||| \label{lemma}
\end{eqnarray}
where $I_h$ is the interpolation operator.\\

\noindent {\it Proof:} Let
$$
T=\int_{\Omega } f(v-I_hv)dx-{a}_h(t;u_h^s,v-I_hv)-b_h(t;u_h^s,v-I_hv).
$$
Integration by parts yields
\begin{eqnarray*}
T &=& \sum \limits_{K \in {\xi}_{h}} \int_K (f+\epsilon\Delta u_h^s-\vec{\beta }\cdot\nabla u_h^s-r(u_h^s))(v-I_hv)dx\\
& & -\sum \limits_{K \in {\xi}_{h}} \int_{\partial K}\epsilon\nabla u_h^s\cdot \vec{n}_K(v-I_hv)ds\\
& & +\sum \limits_{K \in {\xi}_{h}} \int_{\partial K^-\setminus \partial\Omega }\vec{\beta }\cdot \vec{n}_k(u_h^s-u_h^{s,out})(v-I_hv)ds\\
&=&T_1+T_2+T_3.
\end{eqnarray*}
Adding and substracting the data approximation terms into the term $T_1$
\begin{eqnarray*}
T_1 &=& \sum \limits_{K \in {\xi}_{h}} \int_K (f_h +\epsilon\Delta u_h^s-\vec{\beta }_h\cdot\nabla u_h^s-r(u_h^s))(v-I_hv)dx\\
& & +\sum \limits_{K \in {\xi}_{h}} \int_K ((f-f_h)-(\vec{\beta }-\vec{\beta }_h)\cdot\nabla u_h^s)(v-I_hv)dx.
\end{eqnarray*}
Using the Cauchy-Schwarz inequality and interpolation operator identity (\ref{int2})
\begin{eqnarray*}
T_1 &\lesssim & \left( \sum \limits_{K \in {\xi}_{h}} \eta_{R_K}^2\right)^{1/2} \left( \sum \limits_{K \in {\xi}_{h}} \rho_{K}^{-2}\| v-I_hv\|_{L^2(K)}^2\right)^{1/2} \\
&& +\left( \sum \limits_{K \in {\xi}_{h}} \Theta_{K}^2\right)^{1/2} \left( \sum \limits_{K \in {\xi}_{h}} \rho_{K}^{-2}\| v-I_hv\|_{L^2(K)}^2\right)^{1/2} \\
&\lesssim & \left( \sum \limits_{K \in {\xi}_{h}} (\eta_{R_K}^2+\Theta_K^2)\right)^{1/2}|||v|||.
\end{eqnarray*}
For the terms $T_2$ and $T_3$, we have \cite{schotzau09rae}[Lemma 4.8]

$$
T_2 \lesssim   \left( \sum \limits_{K \in {\xi}_{h}} \eta_{E_K}^2 \right)^{1/2}|||v|||
$$
$$
T_3 \lesssim   \left( \sum \limits_{K \in {\xi}_{h}} \eta_{J_K}^2 \right)^{1/2}|||v|||
$$
which finishes the proof.\\

\noindent {\bf Bound to the conforming part of the error:} For a given $t\in (0,T]$, there holds
\begin{eqnarray}
\| u^s-u_h^c\|_{DG}\lesssim \eta + \Theta . \label{conf}
\end{eqnarray}
\noindent {\it Proof:} Since $u^s-u_h^c\in H_0^1(\Omega )$, we have $|u^s-u_h^c|_C=|\vec{\beta} (u^s-u_h^c)|_*$. Then, from the inf-sup condition (\ref{infsup})
$$
\| u^s-u_h^c\|_{DG} = |||u^s-u_h^c||| +|u^s-u_h^c|_C \lesssim \underset{v\in H_0^1(\Omega ) \setminus \{ 0\}}{\text{sup}}\frac{{a}_h(t;u^s-u_h^c,v)}{|||v|||} .
$$
So, we need to bound the term ${a}_h(t;u^s-u_h^c,v)$. Using that $u^s-u_h^c\in H_0^1(\Omega )$, we have
\begin{eqnarray*}
{a}_h(t;u^s-u_h^c,v) &=& {a}_h(t;u^s,v)-{a}_h(t;u_h^c,v) \\
 &=& \int_{\Omega }fv dx -b_h(t;u^s,v)-{a}_h(t;u_h^c,v)\\
&=& \int_{\Omega }fv dx -b_h(t;u^s,v)-D_h(u_h^c,v)-J_h(u_h^c,v)-O_h(u_h^c,v) \\
&=& \int_{\Omega }fv dx -b_h(t;u_h^s,v)+b_h(t;u_h^s,v)-b_h(t;u^s,v)\\
&& \quad -{a}_h(t;u_h^s,v)+D_h(u_h^r,v)+J_h(u_h^r,v)+O_h(u_h^r,v).
\end{eqnarray*}
We also have from the SIPG scheme
\begin{eqnarray*}
\int_{\Omega }fI_hv dx &=& a_h(t;u_h^s,I_hv)+K_h(u_h^s,I_hv)+b_h(t;u_h^s,I_hv)\\
\end{eqnarray*}
Hence, we obtain
$$
a(t;u^s-u_h^c,v)=T_1+T_2+T_3+T_4
$$
\begin{eqnarray*}
T_1 &=& \int_{\Omega }f(v-I_hv) dx -{a}_h(t;u_h^s,v-I_hv)-b_h(t;u_h^s,v-I_hv)\\
T_2 &=& D_h(u_h^r,v)+J_h(u_h^r,v)+O_h(u_h^r,v)\\
T_3 &=& K_h(u_h^s,I_hv) \\
T_4 &=& b_h(t;u_h^s,v)-b_h(t;u^s,v)\\
\end{eqnarray*}
From the inequality (\ref{lemma}), we have
$$
T_1 \lesssim (\eta +\Theta )|||v|||
$$
The continuity results (\ref{cd}-\ref{cj}) and the bound to remainder term (\ref{rem}) yields
$$
T_2 \lesssim (|||u_h^r|||+|\vec{\beta} u_h^r|_*)|||v||| \leq \eta |||v|||
$$
Moreover, using the identities (\ref{ck}) and (\ref{int1}), we get
$$
T_3 \lesssim \gamma^{-1/2}\left( \sum_{K\in \xi} \eta_{J_K}^2\right)^{1/2} |||I_hv||| \lesssim \gamma^{-1/2}\left( \sum_{K\in \xi} \eta_{J_K}^2\right)^{1/2} |||v|||.
$$
Finally, using Cauchy-Shwarz inequality and the boundedness property (\ref{cb}), we get
\begin{eqnarray*}
T_4 &=& b_h(t;u_h^s,v)-b_h(t;u^s,v) = \int_{\Omega }r(u_h^s)v dx-\int_{\Omega }r(u^s)v dx \\
& \leq & C_1 \| v\|_{L^2(\Omega )} - C_2 \| v\|_{L^2(\Omega )} \\
& \lesssim & |||v|||.
\end{eqnarray*}
This finishes the proof.\\

\noindent {\bf Proof to the reliability:} Combining the bounds (\ref{rem}) and (\ref{conf}) to the remainder and the conforming parts of the error, respectively, we obtain
\begin{eqnarray*}
\|u^s-u_h^s \|_{DG} &\leq& \|u^s-u_h^c \|_{DG} +\|u_h^r \|_{DG}\\
&\leq& \eta + \Theta + \eta \\
&\lesssim & \eta + \Theta
\end{eqnarray*}

\noindent {\bf Proof to the efficiency:} The proof of the efficiency is similar to Theorem 3.3 in \cite{schotzau09rae}. We only use the boundedness property (\ref{nonl1}) of the non-linear reaction term to bound the terms occurring in the procedure in \cite{schotzau09rae}.

\subsection{A posteriori error bounds for the semi-discrete system}
\label{semidiscerror}

In order to measure the error for the semi-discrete problem, we use the $L^{\infty}(L^2)+L^2(H^1)$-type norm
$$
\| v\|_*^2 = \| v\|_{L^{\infty}(0,T;L^2(\Omega ))}^2 + \int_0^T |||v|||^2dt
$$
We make use the elliptic reconstruction technique in \cite{makridakis03era}: the elliptic reconstruction $w\in H_0^1(\Omega )$ is defined as the unique solution of the problem
\begin{equation}\label{cont_elliptic}
a(t;w,v)+b_h(t;w,v)=\int_{\Omega}\left( f-\frac{\partial u_h}{\partial t}\right)vdx \; , \quad \forall v\in H_0^1(\Omega ).
\end{equation}
The SIPG discretization, on the other hand, of the above system reads as: for each $t\in (0,T]$, find $w_h\in C^{0,1}(0,T;V_h(\xi_h))$ such that for all $v_h\in V_h(\xi_h)$
$$
a_h(t;w_h,v_h)+K_h(w_h,v_h)+b_h(t;w_h,v_h)= \int_{\Omega}\left( f-\frac{\partial u_h}{\partial t}\right)v_hdx
$$
which implies that $w_h=u_h$. Hence, the error bound to the term $\| w-u_h\|_{DG}$ can be found using the a posteriori error bound (\ref{rel}) for non-linear stationary problem.

We give the a posteriori error bounds for the semi-discrete system developed as in \cite{cangiani13adg}. For the error $e(t)=u(t)-u_h(t)$ of the semi-discrete problem, we set the decomposition $e(t)=\mu (t)+\nu (t)$ with $\mu (t)= u(t)-w(t)$ and $\nu = w(t)-u_h(t)$. Then, for any $t\in (0,T]$, using (\ref{eqweak}) and (\ref{cont_elliptic}), we obtain that
$$
\int_{\Omega}\frac{\partial e}{\partial t}v dx + a_h(t;\mu , v)+ b_h(t;\mu , v)=0
$$
which yields the error bound \cite{cangiani13adg}[Theorem 5.4]
$$
\| e\|_* \lesssim \tilde{\eta}
$$
where the error estimator $\tilde{\eta}$ is defined by
\begin{eqnarray*}
\tilde{\eta}^2 = \| e(0)\|^2 + \int_0^T \tilde{\eta}_{S_1}^2dt + \text{min}\left\{ \left( \int_0^T \tilde{\eta}_{S_2}^2 dt \right)^2 , \rho_T^2\int_0^T \tilde{\eta}_{S_2}^2dt \right\} + \underset{0\leq t\leq T}{\text{max}}\tilde{\eta}_{S_3}^2
\end{eqnarray*}
with
\begin{eqnarray*}
\tilde{\eta}_{S_1}^2 &=& \sum_{K\in{\xi}_{h}} \rho_K^2\left\| f- \frac{\partial u_h}{\partial t} +\epsilon\Delta u_h - \vec{\beta}\cdot\nabla u_h - r(u_h)\right\|_{L^2(K)}^2 + \sum_{e\in\Gamma_0} \epsilon^{-1/2}\rho_e\| [ \epsilon\nabla u_h ] \|_{L^2(e)}^2 \\
&& + \sum_{e\in\Gamma_0\cup\Gamma^D}\left( \frac{\epsilon\sigma}{h_e}+\kappa h_e+\frac{h_e}{\epsilon}\right)\| [ u_h ] \|_{L^2(e)}^2\\
\tilde{\eta}_{S_2}^2 &=& \sum_{e\in\Gamma_0\cup\Gamma^D} h_e\left\| \left[ \frac{\partial u_h}{\partial t}\right]\right\|_{L^2(e)}^2 \\
\tilde{\eta}_{S_3}^2 &=& \sum_{e\in\Gamma_0\cup\Gamma^D} h_e\| [ u_h]\|_{L^2(e)}^2
\end{eqnarray*}
and the weight $\rho_T=\text{min}(\epsilon^{-\frac{1}{2}}, \kappa^{-\frac{1}{2}})$.\\

\subsection{A posteriori error bounds for the fully-discrete system}
\label{fullydiscerror}

For the fully discrete case, we consider the solutions at discrete time instances. For this reason, let $A^k\in V_h^k$ be the unique solution of the stationary system
$$
a_h(t^k; u_h^k,v_h^k) + K_h(u_h^k,v_h^k) + b_h(t^k ; u_h^k,v_h^k) = \int_{\Omega} A^k v_h^k dx \; , \quad \forall v_h^k \in V_h^k
$$
Note that for $k\geq 1$, $A^k=I_h^kf^k - (u_h^k - I_h^ku_h^{k-1})/\tau_k$ with $I_h^k$ being the $L^2$-projection operator onto the space $V_h^k$. Then, the elliptic reconstruction $w^k\in H_0^1(\Omega )$ is defined as the unique solution of the stationary system
\begin{equation}\label{disc_elliptic}
a(t^k; w^k,v) + b(t^k; w^k,v)=\int_{\Omega} A^k v dx \; , \qquad \forall v\in H_0^1(\Omega ).
\end{equation}

Next, we take, as in \cite{cangiani13adg}, the discrete solution $u_h(t)$ as a piecewise continuous function so that on each interval $(t^{k-1},t^k]$, $u_h(t)$ is the linear interpolation of the values $u_h^{k-1}$ and $u_h^k$ given by
$$
u_h(t)= l_{k-1}(t)u_h^{k-1} + l_{k}(t)u_h^{k}
$$
with the linear Lagrange interpolation basis $l_{k-1}$ and $l_{k}$ defined on $[t^{k-1},t^k]$. Then, for the error $e=u-u_h$, using (\ref{eqweak}) and (\ref{disc_elliptic}), we obtain for all $v\in H_0^1(\Omega)$
\begin{eqnarray*}
\int_{\Omega} \frac{\partial e}{\partial t} v dx + a_h(t; e,v) + b_h(t; e,v) &=& \int_{\Omega}(f-f^k)v dx + \int_{\Omega}\left( f^k - \frac{\partial u_h}{\partial t}\right) v dx \\
&& - a_h(t;u_h,v) - b_h(t; u_h,v)
\end{eqnarray*}
which yields the error bound \cite{cangiani13adg}[Theorem 6.5]
$$
\| e\|_*^2 \lesssim \eta_S^2 + \eta_T^2
$$
where $\eta_S$ is the spatial estimator given by \cite{cangiani13adg}

\begin{eqnarray*}
\eta_S^2 &=& \| e(0)\|^2 + \frac{1}{3}\sum_{k=1}^n \tau_k (\eta_{S_{1,k-1}}^2+\eta_{S_{1,k}}^2) +\sum_{k=1}^n \tau_k \eta_{S_{2,k}}^2 \\
&& + \underset{0\leq k\leq n}{\text{max}}\eta_{S_{3,k}}^2 + \text{min} \left\{ \left( \sum_{k=1}^n \tau_k \eta_{S_{4,k}}\right)^2 , \rho_T^2 \sum_{k=1}^n \tau_j \eta_{S_{4,k}}^2 \right\}
\end{eqnarray*}
with
\begin{eqnarray}
{\eta}_{S_{1,k}}^2 &=& \sum_{K\in{\xi}_{h}^{k-1}\cup {\xi}_{h}^{k}} \rho_K^2\left\| A^k +\epsilon\Delta u_h^k - \vec{\beta}^k\cdot\nabla u_h^k - r(u_h^k)\right\|_{L^2(K)}^2 + \sum_{e\in\Gamma_0} \epsilon^{-1/2}\rho_e\| [ \epsilon\nabla u_h^k ] \|_{L^2(e)}^2  \nonumber \\
&& + \sum_{e\in\Gamma}\left( \frac{\epsilon\sigma}{h_e}+\kappa h_e+\frac{h_e}{\epsilon}\right)\| [ u_h^k ] \|_{L^2(e)}^2 \label{spatialest} \\
{\eta}_{S_{2,k}}^2 &=&  \sum_{K\in{\xi}_{h}^{k-1}\cup {\xi}_{h}^{k}} \rho_K^2 \left\| f^k - I_h^kf^k + \frac{u_h^{k-1}-I_h^ku_h^{k-1}}{\tau_k}\right\|_{L^2(K)}^2 \nonumber \\
{\eta}_{S_{3,k}}^2 &=& \sum_{e\in\Gamma} h_e\| [ u_h^k]\|_{L^2(e)}^2 \nonumber \\
{\eta}_{S_{4,k}}^2 &=& \sum_{e\in\Gamma} h_e \left\| \left[ \frac{u_h^k-u_h^{k-1}}{\tau_k}\right]\right\|_{L^2(e)}^2 \nonumber
\end{eqnarray}
and $\eta_T$ is the temporal estimator given by \cite{cangiani13adg}

\begin{eqnarray} \label{timeest}
\eta_T^2 = \sum_{k=1}^n\int_{t_{k-1}}^{t_k}\eta_{T_1,k}^2 dt + \text{min} \left\{ \left( \sum_{k=1}^n \int_{t_{k-1}}^{t_k} \eta_{T_{2,k}}dt\right)^2 , \rho_T^2 \sum_{k=1}^n \int_{t_{k-1}}^{t_k} \eta_{T_{2,k}}^2dt \right\}
\end{eqnarray}
with
\begin{eqnarray*}
\eta_{T_{1,k}}^2 &=&  \sum_{K\in {\xi}_{h}^{k-1}\cup {\xi}_{h}^{k}}\epsilon^{-1} \| l_k(\vec{\beta}^k-\vec{\beta})u_h^k + l_{k-1}(\vec{\beta}^{k-1}-\vec{\beta})u_h^{k-1}\|_{L^2(K)}^2  \\
\eta_{T_{2,k}}^2 &=& \sum_{K\in {\xi}_{h}^{k-1}\cup {\xi}_{h}^{k}} \| f-f^k + l_{k-1}(A^k - A^{k-1}) + l_k(\nabla\cdot\vec{\beta}^k-\nabla\cdot \vec{\beta})u_h^k + l_k(\nabla\cdot\vec{\beta}^{j-1}-\nabla\cdot \vec{\beta})u_h^{k-1}\|_{L^2(K)}^2
\end{eqnarray*}


\section{Adaptive algorithm}
\label{algorithm}

\begin{figure}
\centering
\setlength{\unitlength}{0.5mm}
\begin{picture}(173,210)

\put(94,198){\oval(40,17)}
\put(79,198){\text{Known $u_h^{k-1}$}}
\put(81,192){\text{on $\xi_h^{k-1}$}}
\put(94,190){\vector(0,-1){20}}
\put(96,178){\text{Solve for $u_h^k$ on $\xi_h^{k-1}$}}
\put(94,162){\oval(40,17)}
\put(81,162){\text{$\tilde{\eta}_{T,k}^2 > \mathbf{ttol}$}}
\put(75,158){\vector(-2,-1){40}}
\put(50,154){\text{YES}}
\put(113,158){\vector(2,-1){40}}
\put(126,154){\text{NO}}
\put(19,138){\line(1,0){43}}
\put(19,138){\line(0,-1){25}}
\put(19,113){\line(1,0){43}}
\put(62,113){\line(0,1){25}}
\put(24,131){\text{Half current}}
\put(24,125){\text{interval until }}
\put(24,118){\text{$\tilde{\eta}_{T,k}^2 \leq \mathbf{ttol}$}}
\put(40,113){\vector(0,-1){20}}
\put(19,93){\line(1,0){43}}
\put(19,93){\line(0,-1){22}}
\put(19,71){\line(1,0){43}}
\put(62,71){\line(0,1){22}}
\put(27,86){\text{Resolve $u_h^k$}}
\put(30,77){\text{on $\xi_h^{k-1}$}}
\put(62,82){\vector(3,2){65}}
\put(129,138){\line(1,0){43}}
\put(129,138){\line(0,-1){25}}
\put(129,113){\line(1,0){43}}
\put(172,113){\line(0,1){25}}
\put(139,129){\text{Refine for}}
\put(131,119){\text{$\eta_{S_1,k}^2|_K > \mathbf{stol^+}$}}
\put(152,113){\vector(0,-1){20}}
\put(129,93){\line(1,0){43}}
\put(129,93){\line(0,-1){25}}
\put(129,68){\line(1,0){43}}
\put(172,68){\line(0,1){25}}
\put(136,84){\text{Coarsen for}}
\put(131,75){\text{$\eta_{S_1,k}^2|_K < \mathbf{stol^-}$}}
\put(152,68){\vector(-2,-1){40}}
\put(74,53){\line(1,0){38}}
\put(74,53){\line(0,-1){18}}
\put(74,35){\line(1,0){38}}
\put(112,35){\line(0,1){18}}
\put(81,44){\text{Create $\xi_h^k$}}
\put(93,35){\vector(0,-1){10}}
\put(74,25){\line(1,0){38}}
\put(74,25){\line(0,-1){20}}
\put(74,05){\line(1,0){38}}
\put(112,05){\line(0,1){20}}
\put(78,18){\text{Solve for $u_h^k$}}
\put(85,10){\text{on $\xi_h^k$}}
\put(74,15){\vector(-1,0){15}}
\put(21,25){\line(1,0){38}}
\put(21,25){\line(0,-1){20}}
\put(21,05){\line(1,0){38}}
\put(59,05){\line(0,1){20}}
\put(36,18){\text{Set}}
\put(28,11){\text{$k=k+1$}}
\put(21,15){\line(-1,0){20}}
\put(1,15){\line(0,1){185}}
\put(1,200){\vector(1,0){70}}

\end{picture}
\caption{Adaptive algorithm chart\label{chart}}
\end{figure}
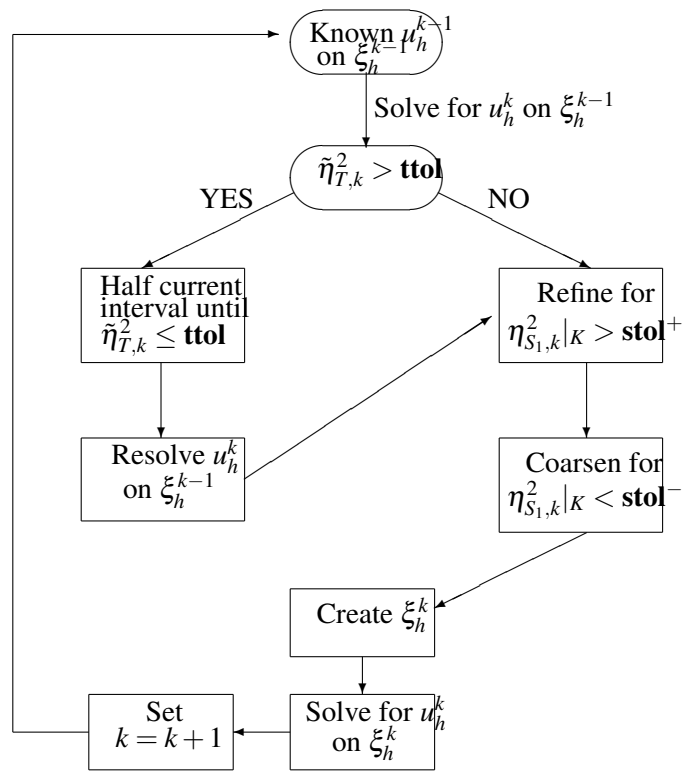

The time-space adaptive algorithm for linear convection-diffusion problems in \cite{cangiani13adg} is modified for semi-linear problems of type \eqref{model}, see Fig.~\ref{chart}. It is based on the residual-based a posteriori error estimators given in the previous sections. The algorithm starts with an initial uniform mesh in space and with a given initial time step. At each time step, the space and time-step are adaptively arranged according to the user defined tolerances $\mathbf{ttol}$ for time-step refinement, and $\mathbf{stol^+}$ and $\mathbf{stol^-}$ for spatial mesh, former corresponding to refinement and latter corresponding to the coarsening procedures in space. Note that we do not need a temporal tolerance corresponding to time-step coarsening, since we start, in our problems, with a uniform equispaced distribution of $[0,T]$ having sufficiently large time-steps. Thus, it is enough just bisecting the time intervals not satisfying the temporal tolerance $\mathbf{ttol}$. Both the refinement and coarsening processes in space are determined by the indicator $\eta_{S_{1,k}}$ (\ref{spatialest}) appearing in the spatial estimator. Since the temporal estimator $\eta_T$ (\ref{timeest}) is not easy to compute, the adaptive refinement of the time-steps are driven by the modified temporal-indicator \cite{cangiani13adg}
$$
\tilde{\eta}_{T_k}^2=\int_{t_{k-1}}^{t_k} \eta_{T_{1,k}}^2dt + \min \{\rho_T,T\}\int_{t_{k-1}}^{t_k} \eta_{T_{2,k}}^2dt
$$
sum of which gives a bound for the temporal estimator $\eta_T^2$.

Although the adaptive algorithm, Fig.~\ref{chart}, stands for a single equation of the system \eqref{org}, it is not difficult to extend the algorithm to the coupled systems. For this, say we have a system of two equations for the unknowns $u_1$ and $u_2$, the temporal-indicators $\tilde{\eta}_{T_k}^1$, $\tilde{\eta}_{T_k}^2$ and the spatial-indicators $\eta_{S_{1,k}}^1$, $\eta_{S_{1,k}}^2$ corresponding to the unknowns $u_1$ and $u_2$, respectively, are computed. To adapt the time-step size, we ask the temporal condition for the both temporal-indicators, i.e. $\tilde{\eta}_{T_k}^1\leq \mathbf{ttol}$ and $\tilde{\eta}_{T_k}^2\leq \mathbf{ttol}$. On the other hand, to select the elements to be refined, we take the set of elements which is the union of the sets of the elements satisfying $\eta_{S_{1,k}}^1> \mathbf{stol^+}$ and $\eta_{S_{1,k}}^2> \mathbf{stol^+}$, and similar procedure to select the elements to be coarsened, but not including any elements which are selected to be refined. Numerical studies demonstrate that the  adaptive algorithm is capable of resolving the layers in space as the time progresses.

\section{Solution of the fully-discrete system}
\label{newton}

In this section, we discuss the solution of the fully-discrete system \eqref{fullydisc} on an arbitrary $k^{th}$ time-step, which is solved for all $k=1,2,\ldots , n$. In order to not be confused about the notations, let us consider the system \eqref{fullydisc} on an arbitrary $k^{th}$ time-step without the superscript for the time-step of the form

\begin{equation}\label{singleeq}
\int_{\Omega}\frac{u_{h}-w_{h}}{\tau}v_{h} dx + a_{h}(u_{h},v_{h})+ K_{h}(u_{h},v_{h})+b_{h}(u_{h}, v_{h})=\int_{\Omega}fv_{h} dx \; , \quad \forall v_h\in V_h
\end{equation}
where we set $u_{h}:=u_{h}^k$, $w_h:=u_{h}^{k-1}$, $v_{h}:=v_{h}^k$, $f:=f^k$, $\tau:=\tau_k$, $a_{h}(u_{h},v_{h}):=a_{h}(t^k;u_{h}^k,v_{h}^k)$, $b_{h}(u_{h},v_{h}):=b_{h}(t^k;u_{h}^k,v_{h}^k)$ and $V_h:=V_h^k$. The approximate solution $u_h$ and the known solution (from the previous time-step) $w_h$ of \eqref{singleeq}  have the form
\begin{equation}\label{comb}
u_h=\sum_{i=1}^{Nel}\sum_{l=1}^{Nloc}u_l^i\phi_l^i \; , \quad
w_h=\sum_{i=1}^{Nel}\sum_{l=1}^{Nloc}w_l^i\phi_l^i
\end{equation}
where $\phi_l^i$'s are the basis polynomials spanning the space $V_h$, $u_l^i$'s are the unknown coefficients to be found, $w_l^i$'s are the known coefficients. $Nel$ denotes the number of triangles and $Nloc$ is the number of local dimension depending on the degree of polynomials $k$ (in 2D, $Nloc=(k+1)(k+2)/2$). In DG methods, we choose the piecewise basis polynomials $\phi_l^i$'s in such a way that each basis function has only one triangle as a support, i.e. we choose on a specific triangle $K_e$, $e\in\{ 1,2,\ldots , Nel\}$, the basis polynomials $\phi_l^e$ which are identically zero outside the triangle $K_e$, $l=1,2,\ldots , Nloc$. By this construction, the stiffness and mass matrices obtained by DG methods has a block structure, each of which related to a triangle (there is no overlapping as in continuous FEM case). The product $dof:=Nel*Nloc$ gives the degree of freedom in DG methods. Inserting the linear combinations \eqref{comb} of $u_h$ and $w_h$ in \eqref{singleeq} and choosing the test functions $v_h=\phi_l^i$, $l=1,2,\ldots , Nloc$, $i=1,2,\ldots , Nel$, the discrete residual of the system (\ref{singleeq}) in matrix vector form is given by
\begin{equation}\label{residualsystem}
Res(\vec{u})=M\vec{u} - M\vec{w} +\tau( S\vec{u} + \vec{b}(\vec{u}) - \vec{f}) = 0
\end{equation}
where $\vec{u}, \vec{w}\in\mathbb{R}^{dof}$ are the vector of unknown and known coefficients $u_l^i$'s and $w_l^i$'s, respectively, $M\in\mathbb{R}^{dof\times dof}$ is the mass matrix, $S\in\mathbb{R}^{dof\times dof}$ is the stiffness matrix corresponding to the bilinear form $\tilde{a}_h(u_h,v_h):=a_h(u_h,v_h)+K_h(u_h,v_h)$, $\vec{b}\in\mathbb{R}^{dof}$ is the vector function of $\vec{u}$ related to the non-linear form $b_h(u_h,v_h)$ and $\vec{f}\in\mathbb{R}^{dof}$ is the vector to the linear form $\int_{\Omega}fv_{h}dx$. The explicit definitions are given by
$$
S=
\begin{bmatrix}
S_{11} & S_{12} & \cdots & S_{1,Nel} \\
S_{21} & S_{22} &  & \vdots \\
\vdots &  & \ddots &  \\
S_{Nel,1} & \cdots &  & S_{Nel,Nel}
\end{bmatrix} \; , \quad M=
\begin{bmatrix}
M_{11} & M_{12} & \cdots & M_{1,Nel} \\
M_{21} & M_{22} &  & \vdots \\
\vdots &  & \ddots &  \\
M_{Nel,1} & \cdots &  & M_{Nel,Nel}
\end{bmatrix}
$$

$$
 \vec{u}=
\begin{bmatrix}
 \vec{u}_1 \\
\vec{u}_2 \\
\vdots \\
\vec{u}_{Nel}
\end{bmatrix} \; , \quad \vec{w}=
\begin{bmatrix}
 \vec{w}_1 \\
\vec{w}_2 \\
\vdots \\
\vec{w}_{Nel}
\end{bmatrix} \; , \quad \vec{b}(\vec{u})=
\begin{bmatrix}
 \vec{b}_1(\vec{u}) \\
\vec{b}_2(\vec{u}) \\
\vdots \\
\vec{b}_{Nel}(\vec{u})
\end{bmatrix} \; , \quad \vec{f}=
\begin{bmatrix}
 \vec{f}_1 \\
\vec{f}_2 \\
\vdots \\
\vec{f}_{Nel}
\end{bmatrix}
$$
where all the blocks have dimension $Nloc$:
$$
S_{ji}=
\begin{bmatrix}
\tilde{a}_h(\phi_1^i,\phi_1^j) & \tilde{a}_h(\phi_2^i,\phi_1^j) & \cdots & \tilde{a}_h(\phi_{Nloc}^i,\phi_1^j) \\
\tilde{a}_h(\phi_1^i,\phi_2^j)  & \tilde{a}_h(\phi_2^i,\phi_2^j)  &  & \vdots \\
\vdots &  & \ddots &  \\
\tilde{a}_h(\phi_1^i,\phi_{Nloc}^j)  & \cdots &  & \tilde{a}_h(\phi_{Nloc}^i,\phi_{Nloc}^j)
\end{bmatrix}
$$

$$
M_{ji}=
\begin{bmatrix}
\int_{\Omega}\phi_1^i\phi_1^jdx & \int_{\Omega}\phi_2^i\phi_1^jdx & \cdots & \int_{\Omega}\phi_{Nloc}^i\phi_1^jdx \\
\int_{\Omega}\phi_1^i\phi_2^jdx & \int_{\Omega}\phi_2^i\phi_2^jdx  &  & \vdots \\
\vdots &  & \ddots &  \\
\int_{\Omega}\phi_1^i\phi_{Nloc}^jdx  & \cdots &  & \int_{\Omega}\phi_{Nloc}^i\phi_{Nloc}^jdx
\end{bmatrix}
$$

$$
 \vec{u}_i=
\begin{bmatrix}
 u_1^i \\
u_2^i \\
\vdots \\
u_{Nloc}^i
\end{bmatrix} \; , \quad \vec{w}_i=
\begin{bmatrix}
 w_1^i \\
w_2^i \\
\vdots \\
w_{Nloc}^i
\end{bmatrix} \; , \quad \vec{b}_i(\vec{u})=
\begin{bmatrix}
 b_h(u_h,\phi_1^i) \\
b_h(u_h,\phi_2^i) \\
\vdots \\
b_h(u_h,\phi_{Nloc}^i)
\end{bmatrix} \; , \quad \vec{f}_i=
\begin{bmatrix}
 \int_{\Omega}f\phi_1^idx \\
\int_{\Omega}f\phi_2^idx \\
\vdots \\
\int_{\Omega}f\phi_{Nloc}^idx
\end{bmatrix} .
$$
The block structure of DG methods causes to the increase of the condition number of the obtained matrices by the degree $k$ of basis polynomials, which is a drawback of DG methods comparing to the classical (continuous) finite elements. However, taking into account the locality, a valuable property, of DG methods, this drawback can be handled by various preconditioners developed for DG discretized schemes in the literature. Besides, the condition number of the stiffness matrix $S$ increases linearly by the penalty parameter $\sigma$, as well. Therefor, the penalty parameter should not be chosen too large. On the other hand, it should be selected sufficiently large to ensure the coercivity of the bilinear form \cite{riviere08dgm}[Sec. 27.1], which is needed for the stability of the convergence of the DG method. It ensures that the stiffness matrix arising from the DG discretization is symmetric positive definite. In the literature, several choices of the penalty parameter are suggested. In \cite{epshteyn07epp}, computable lower bounds are derived, and in \cite{dobrev08psi}, the penalty parameter is chosen depending on the diffusion coefficient $\epsilon$. The effect of the penalty parameter on the condition number was discussed in detail for the DG discretization of the Poisson equation in \cite{castillo12pdg} and in \cite{slingerland14fls} for layered reservoirs with strong permeability contrasts. In our study, we select the penalty parameter $\sigma$ depending only on the polynomial degree $k$, as $\sigma =3k(k+1)$ on interior edges and $\sigma =6k(k+1)$ on boundary edges. The reason of the coupling of the penalty parameter $\sigma$ on boundary edges is sufficient penalization of the solution on the boundary due to the non-homogeneous Dirichlet boundary conditions which are imposed weakly in DG methods.

Next, we solve the system \eqref{residualsystem} by Newton method. In the sequel, we start with an initial guess $\vec{u}^{(0)}$ (most possibly $\vec{u}^{(0)}=\vec{w}$, i.e. the known solution from the previous time-step) and we solve the system
\begin{eqnarray*}
J^s\delta\vec{u}^{(s)} &=& -Res(\vec{u}^{(s)}) \\
\vec{u}^{(s+1)} &=&  \vec{u}^{(s)} + \delta\vec{u}^{(s)} \; , \quad s=0,1,\ldots
\end{eqnarray*}
where $J^s=M+\tau (S+J^s_{\vec{b}})$ is the Jacobian matrix of the system at the iteration $\vec{u}^{(s)}$, and  $J^s_{\vec{b}}$ denotes the Jacobian matrix to the non-linear vector $\vec{b}(\vec{u})$  at the iteration $\vec{u}^{(s)}$.

\section{Numerical studies}
\label{numeric}

In this section, we give the numerical studies demonstrating the performance of the time-space adaptive algorithm. All the computations are implemented on MATLAB-R2014a. In the problems, by the very coarse initial mesh, we mean an initial mesh which is formed, for instance on $\Omega =(0,1)^2$, by dividing the domain with $\Delta x_1=\Delta x_2=0.5$ leading to 8 triangular elements and 48 DoFs for quadratic elements.
 As the first example, Example~\ref{ex1}, we give a test example with polynomial type non-linearity having a non-moving internal layer to figure out the benchmark of the algorithm by using different tolerances and diffusion parameters $\epsilon$: rates of error, spatial and temporal estimators, and effectivity indices (proportion of the estimator to the error). We expect that the effectivity indices lie in a small band for different diffusion parameters meaning that our estimators are robust in the system P\'{e}clet number. Moreover, to demonstrate the mentioned properties, we use the weighted DoFs as in \cite{cangiani13adg}
$$
\text{Weighted DoFs} = \frac{1}{T} \sum_{k=1}^n \tau_k\lambda_k
$$
where $\lambda_k$ denotes the total number of DoFs on the union mesh $\xi_h^{k-1}\cup\xi_h^{k}$. Since the first example has a non-moving internal layer, a monotonic increase in the DoFs is expected by the time progresses. Conversely, we give problems having moving layers by the time progresses in Example~\ref{ex2}-\ref{ex3}. In this case, we expect that the refinement and coarsening procedures in space work simultaneously leading to oscillations in time vs DoFs plots. By Example~\ref{ex2}, we also test the performance of our algorithm for a coupled system. As the final example, Example~\ref{ex4}, we consider an important real life problem representing a reaction in porous media having internal layers at different locations due to the high-permeability rocks.

\subsection{Example with polynomial type non-linearity (benchmark of the algorithm)}
\label{ex1}

\begin{figure}[ht]
\centering
\includegraphics[width=0.55\textwidth]{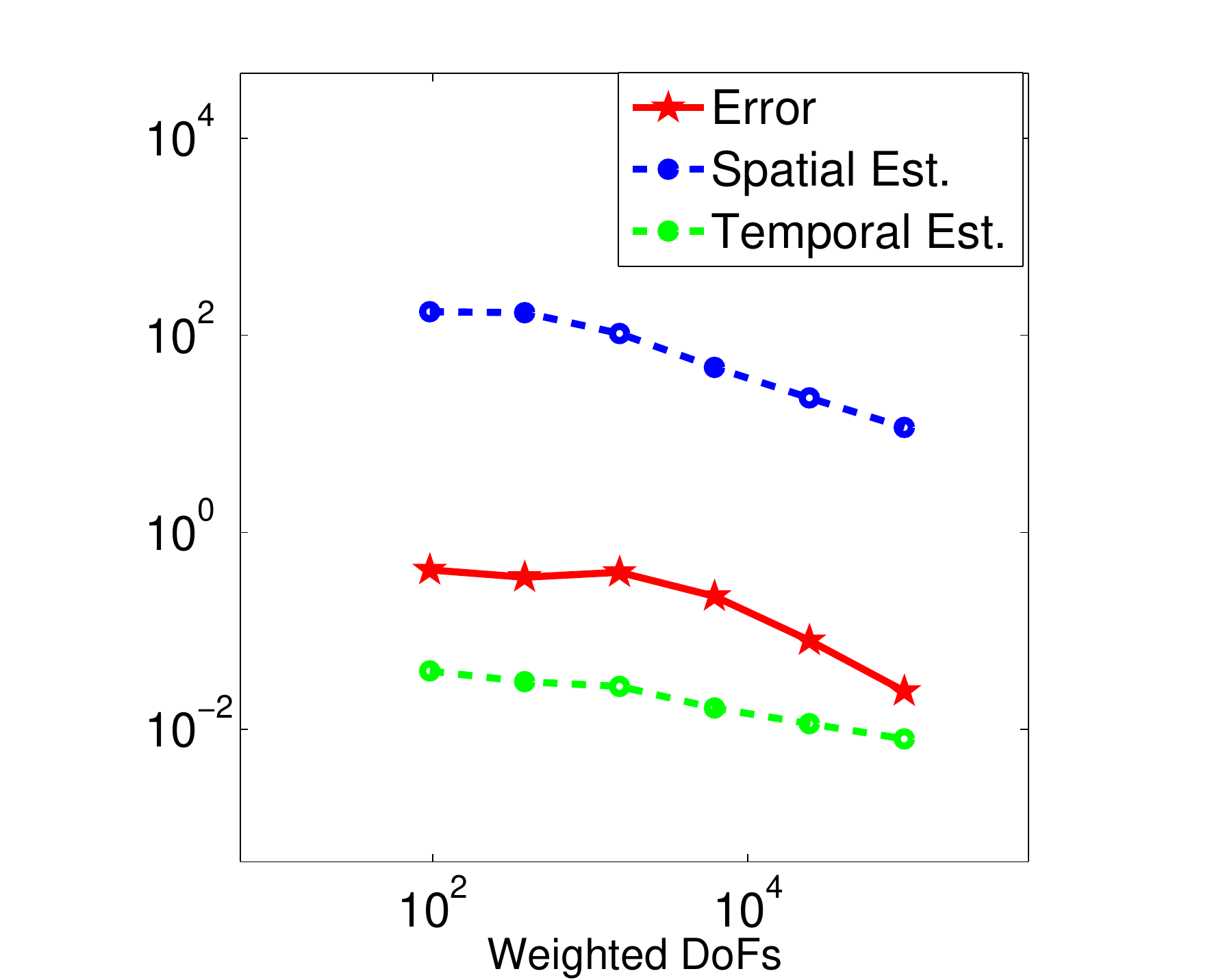}
\caption{Example \ref{ex1}: Decays of estimators and errors for uniform time-space\label{ex1:uniform} }
\end{figure}

The first example is taken from  \cite{bause12ash} with a polynomial non-linear term
$$
u_t+\vec{\beta }\cdot\nabla u-\epsilon\Delta u+r(u)=f \quad \text{in } \; \Omega =(0,1)^2,
$$
with the convection field $\vec{\beta }(x,y)=(2,3)^T$, diffusion coefficient $\epsilon =10^{-6}$, the non-linear reaction term $r(u)=u^4$. The source function $f$ and the Dirichlet boundary condition are chosen so that the exact solution is given by
\begin{eqnarray*}
u(\vec{x},t)&=& 16\sin (\pi t)x_1(1-x_1)x_2(1-x_2) \\
&& [ 0.5 + \pi^{-1}\arctan (2\epsilon^{-1/2}(0.25^2-(x_1-0.5)^2-(x_2-0.5)^2)) ].
\end{eqnarray*}
We start by demonstrating the decrease of the errors by uniform time-space refinement using linear DG elements. In Fig.~\ref{ex1:uniform}, the expected first order convergence in space and time is shown.

\begin{figure}[ht]
\centering
\subfloat{\includegraphics[width=0.5\textwidth]{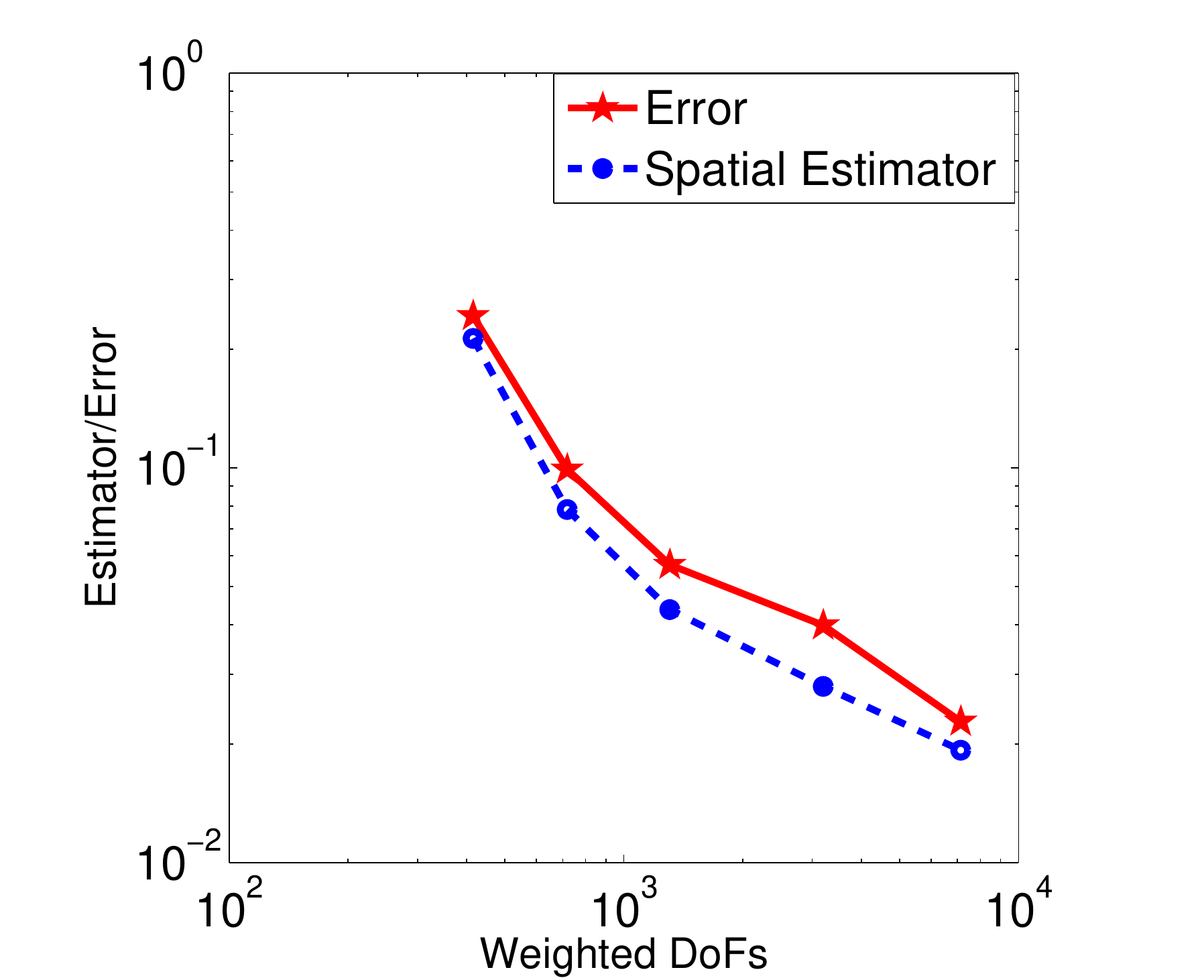}}
\subfloat{\includegraphics[width=0.5\textwidth]{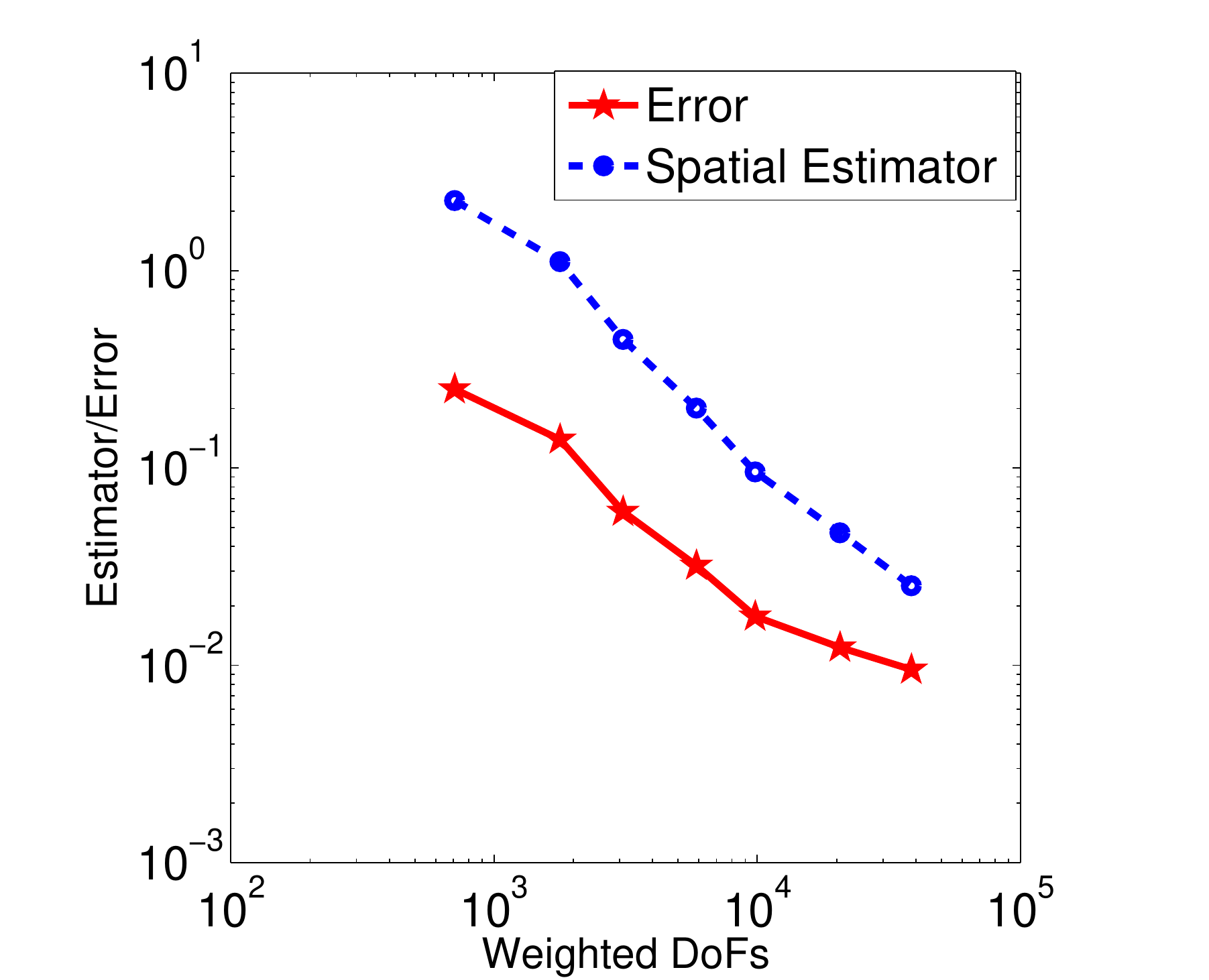}}
\caption{Example \ref{ex1}: Error/spatial estimator for $\epsilon =10^{-2}$ (left) and $\epsilon =10^{-4}$ (right)\label{ex1:errors} }
\end{figure}
\begin{figure}[ht]
\centering
\subfloat{\includegraphics[width=0.5\textwidth]{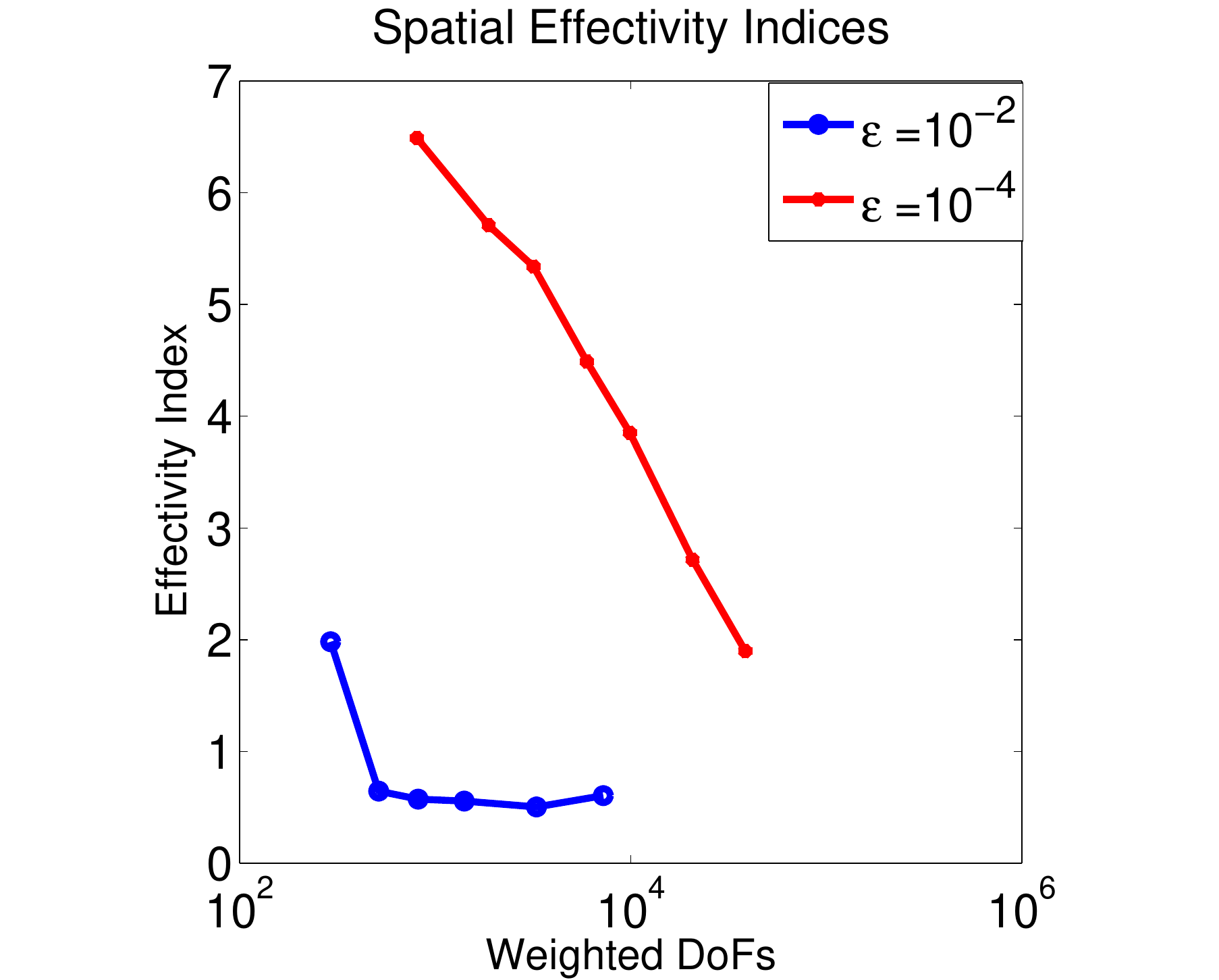}}
\subfloat{\includegraphics[width=0.5\textwidth]{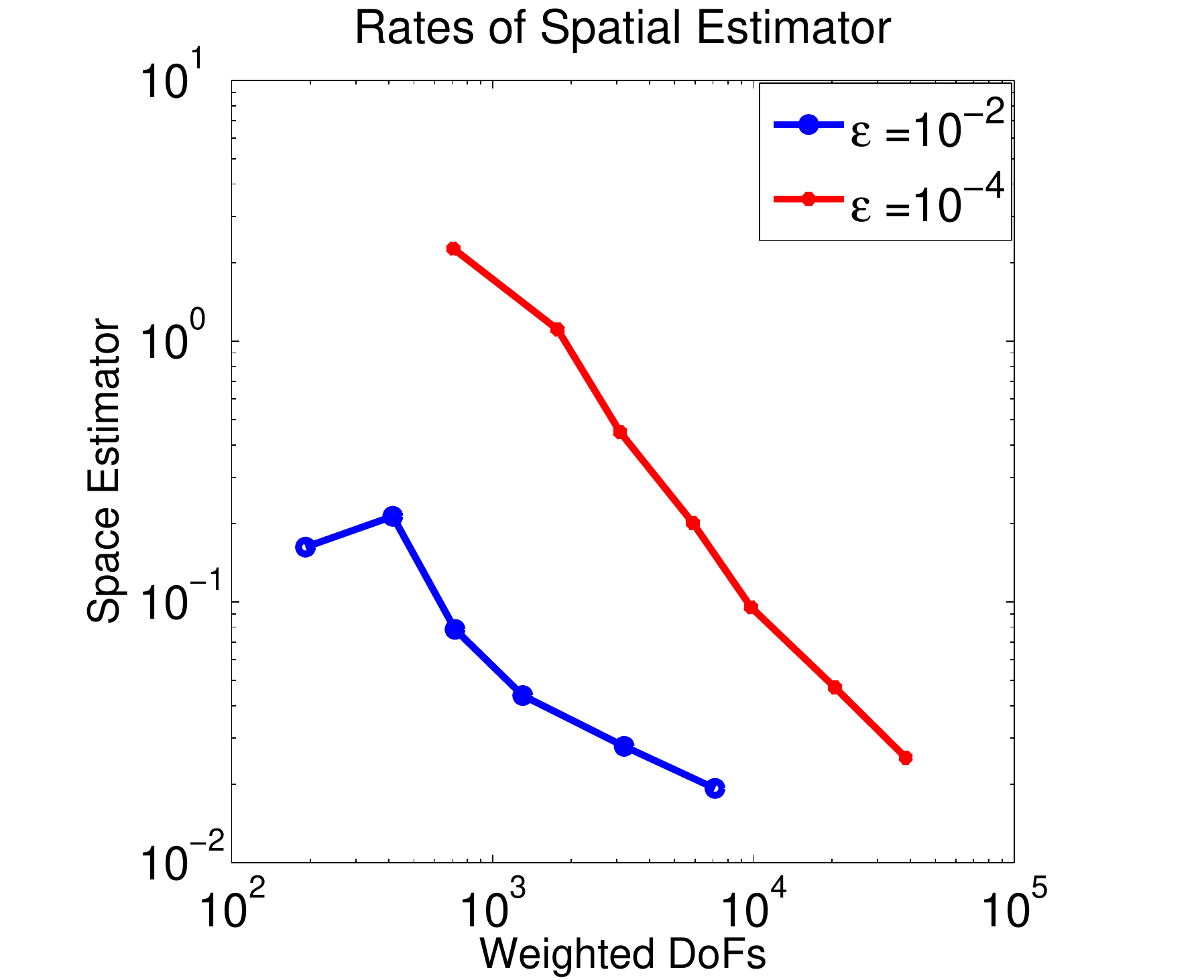}}
\caption{Example \ref{ex1}: Spatial effectivity indices (left) and estimators (right)\label{ex1:spatial}}
\end{figure}

For the time-space adaptive solution, we use quadratic DG elements. We investigate the performance of the spatial estimator by fixing the temporal time-step $\tau =0.005$ so that  the temporal error dominated by the spatial error.  We reduce the spatial estimator tolerance $\mathbf{stol^+}$ from $10^{-1}$ to $10^{-6}$. The rates of the error and the spatial estimator are similar as illustrated in Fig.~\ref{ex1:errors} for $\epsilon =10^{-2}$ and $\epsilon =10^{-4}$. Fig.~\ref{ex1:spatial} shows the spatial effectivity indices and the decrease of the spatial estimators for  the diffusion constant $\epsilon$. The effectivity indices do not exceed 7 which are acceptable as in \cite{cangiani13adg} for linear convection-diffusion problems.

\begin{figure}[ht]
\centering
\subfloat{\includegraphics[width=0.5\textwidth]{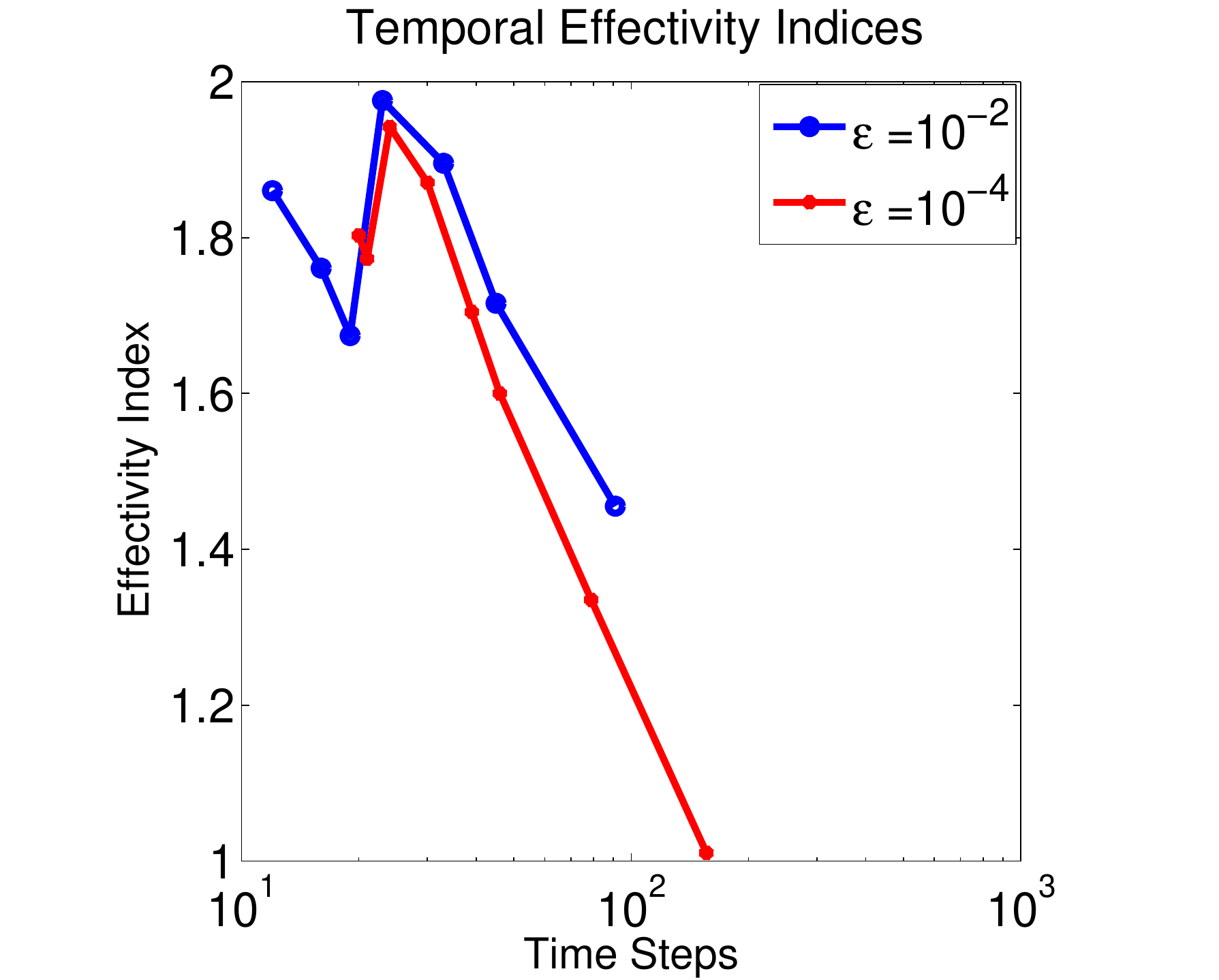}}
\subfloat{\includegraphics[width=0.5\textwidth]{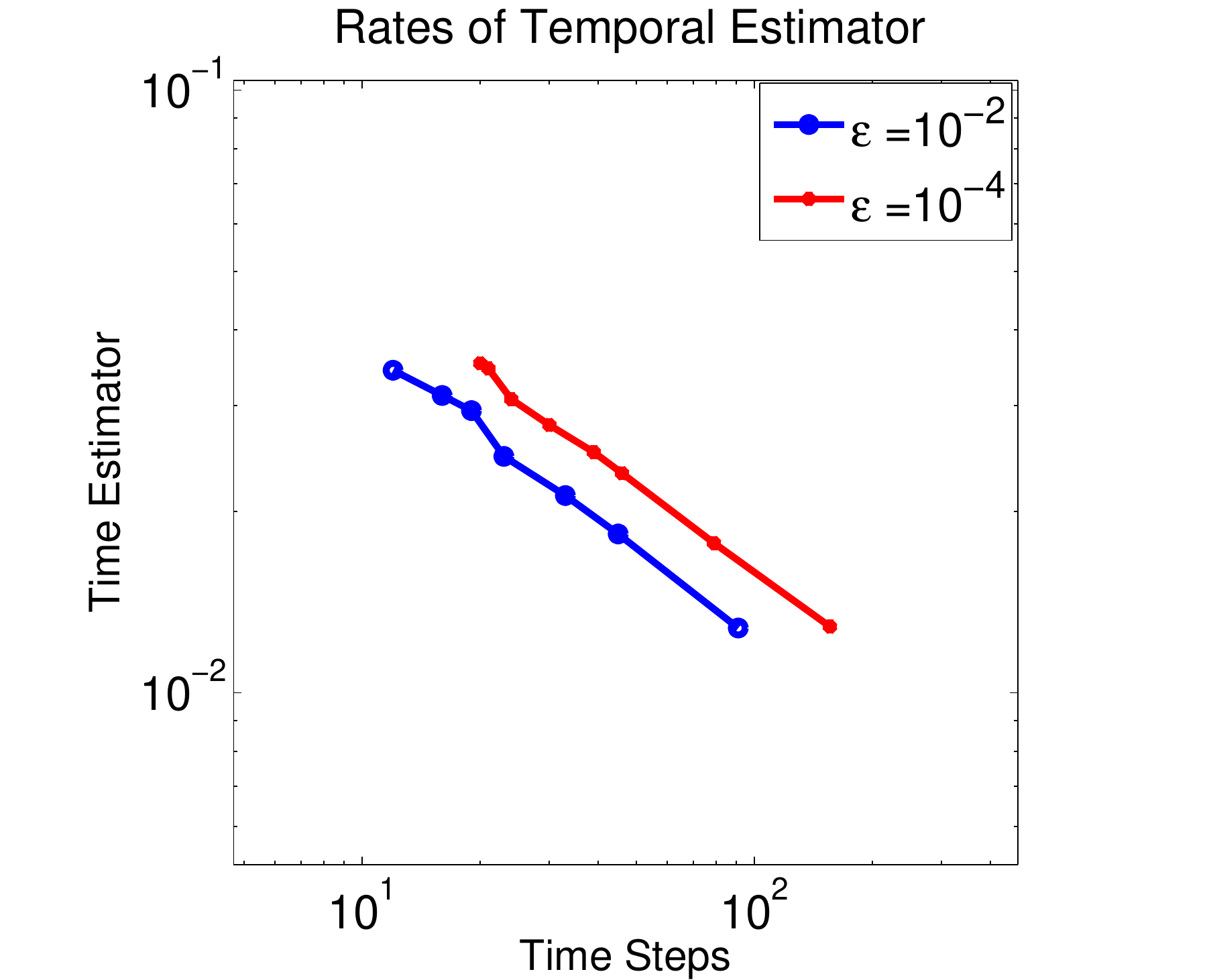}}
\caption{Example \ref{ex1}: Temporal effectivity indices (left) and estimators (right)\label{ex1:temporal}}
\end{figure}

To investigate the performance of the temporal estimator, we fix a sufficiently fine spatial mesh so that the the spatial error dominated by the temporal error, and then we reduce the temporal estimator tolerance $\mathbf{ttol}$ in the range $10^{-1}-10^{-6}$. In Fig.~\ref{ex1:temporal}, the temporal effectivity indices and the decrease of the temporal estimators are not affected by $\epsilon$, and effectivity indices are almost the same within the band 1-2.

\begin{figure}[ht]
\centering
\includegraphics[width=0.5\textwidth]{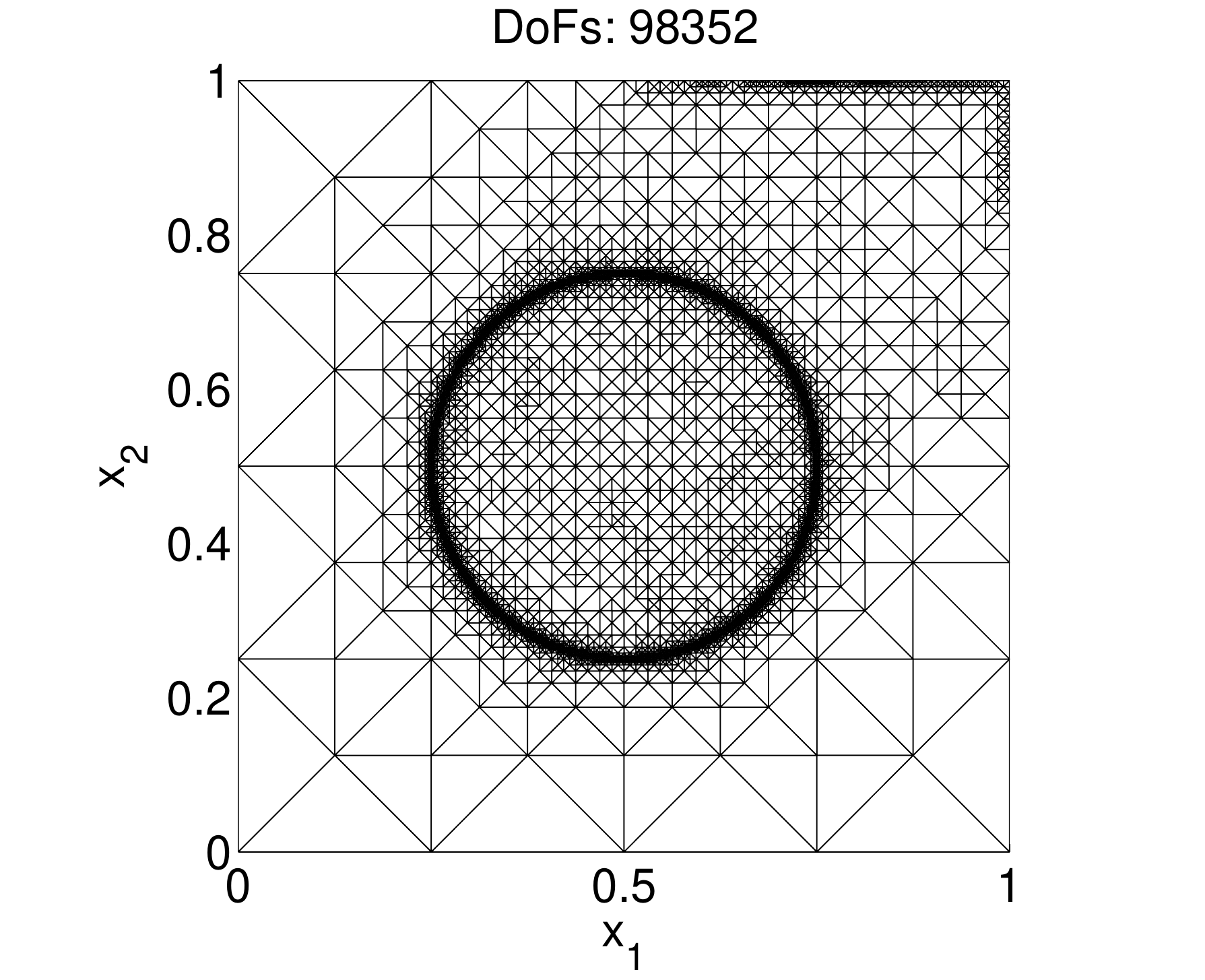}
\caption{Example \ref{ex1}: Adaptive mesh\label{ex1:mesh} }
\end{figure}

\begin{figure}[ht]
\centering
\includegraphics[width=1.2\textwidth]{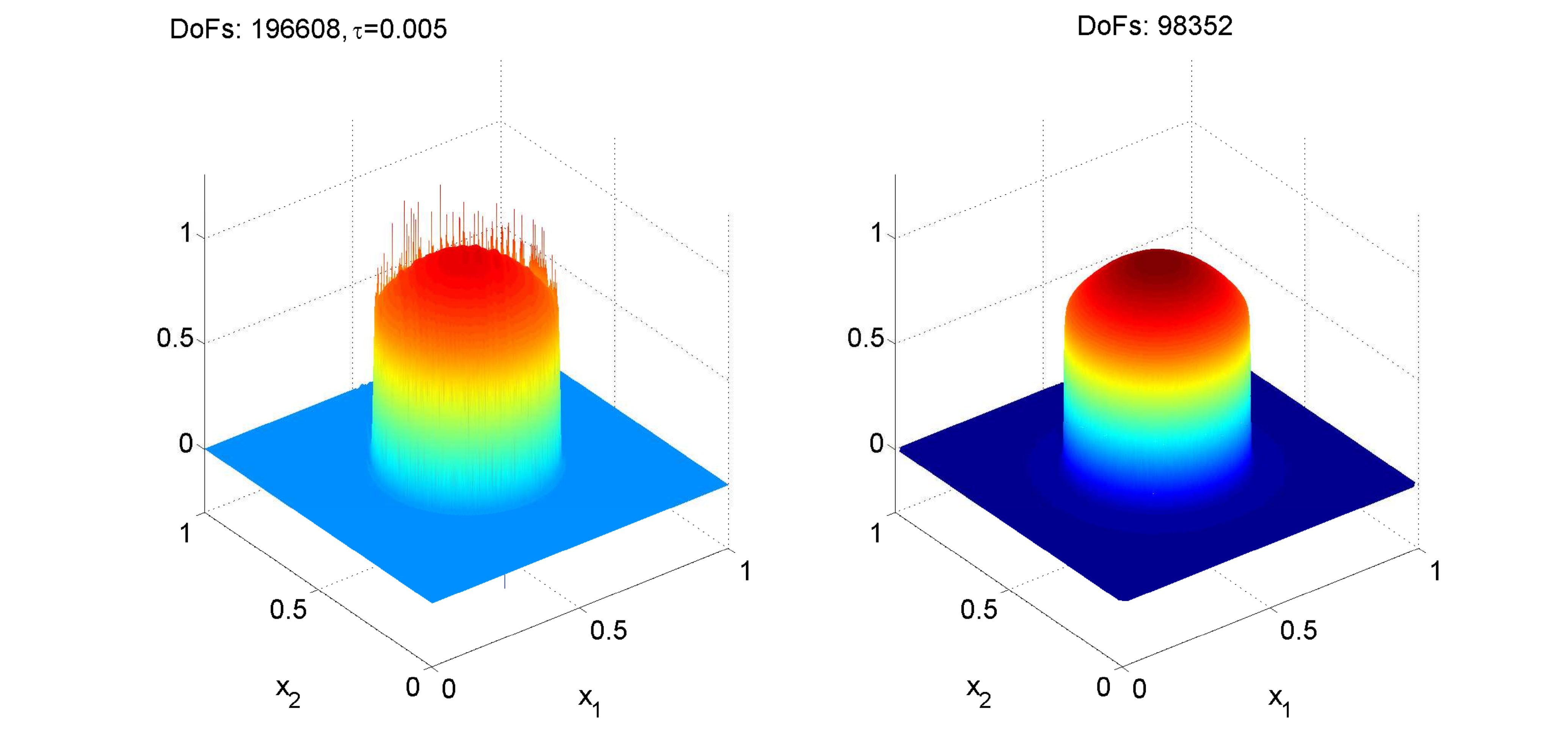}
\caption{Example \ref{ex1}: Uniform (left) and adaptive (right) solutions at T=0.5\label{ex1:plot} }
\end{figure}

\begin{figure}[ht]
\centering
\subfloat{\includegraphics[width=0.5\textwidth]{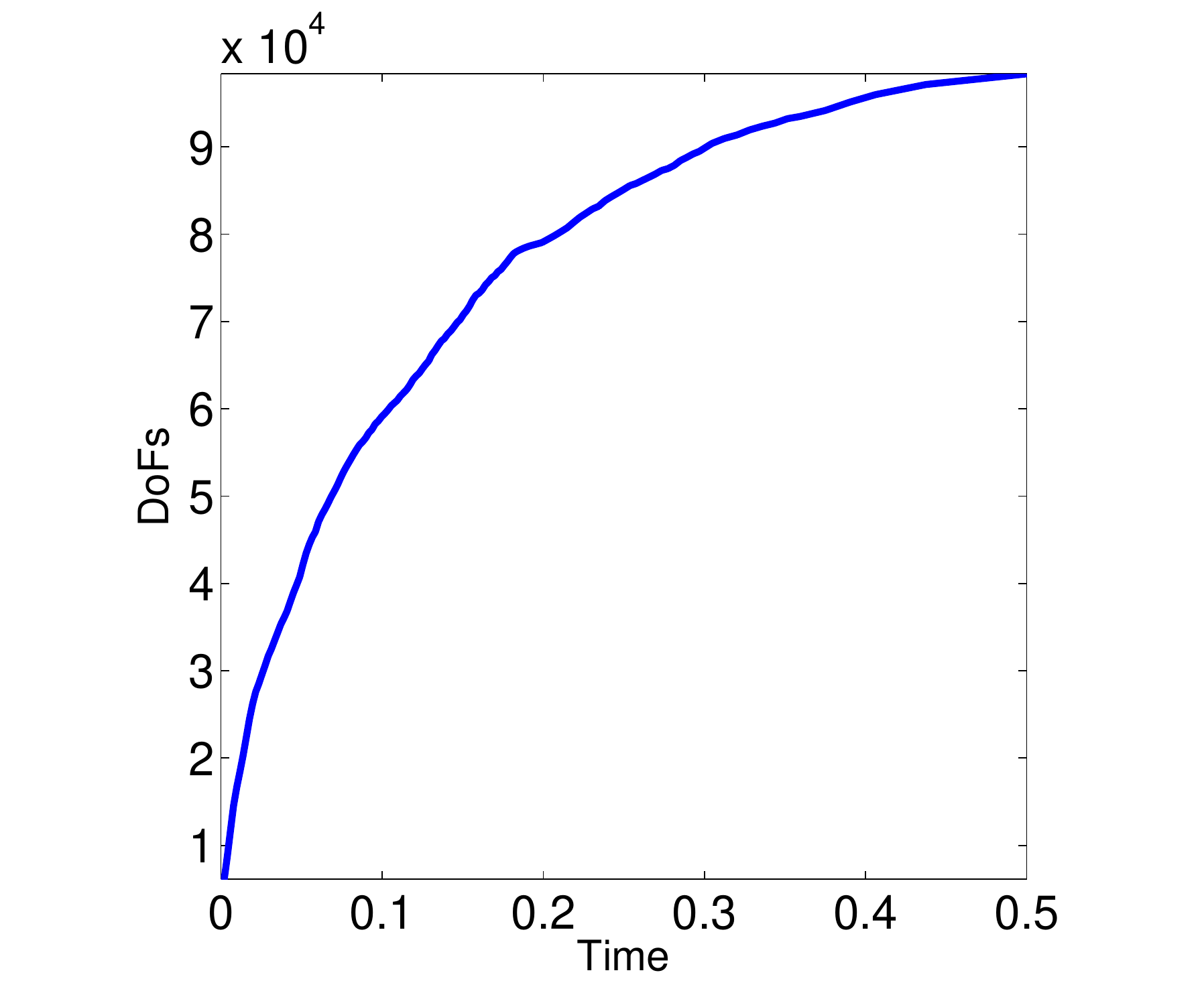}}
\subfloat{\includegraphics[width=0.5\textwidth]{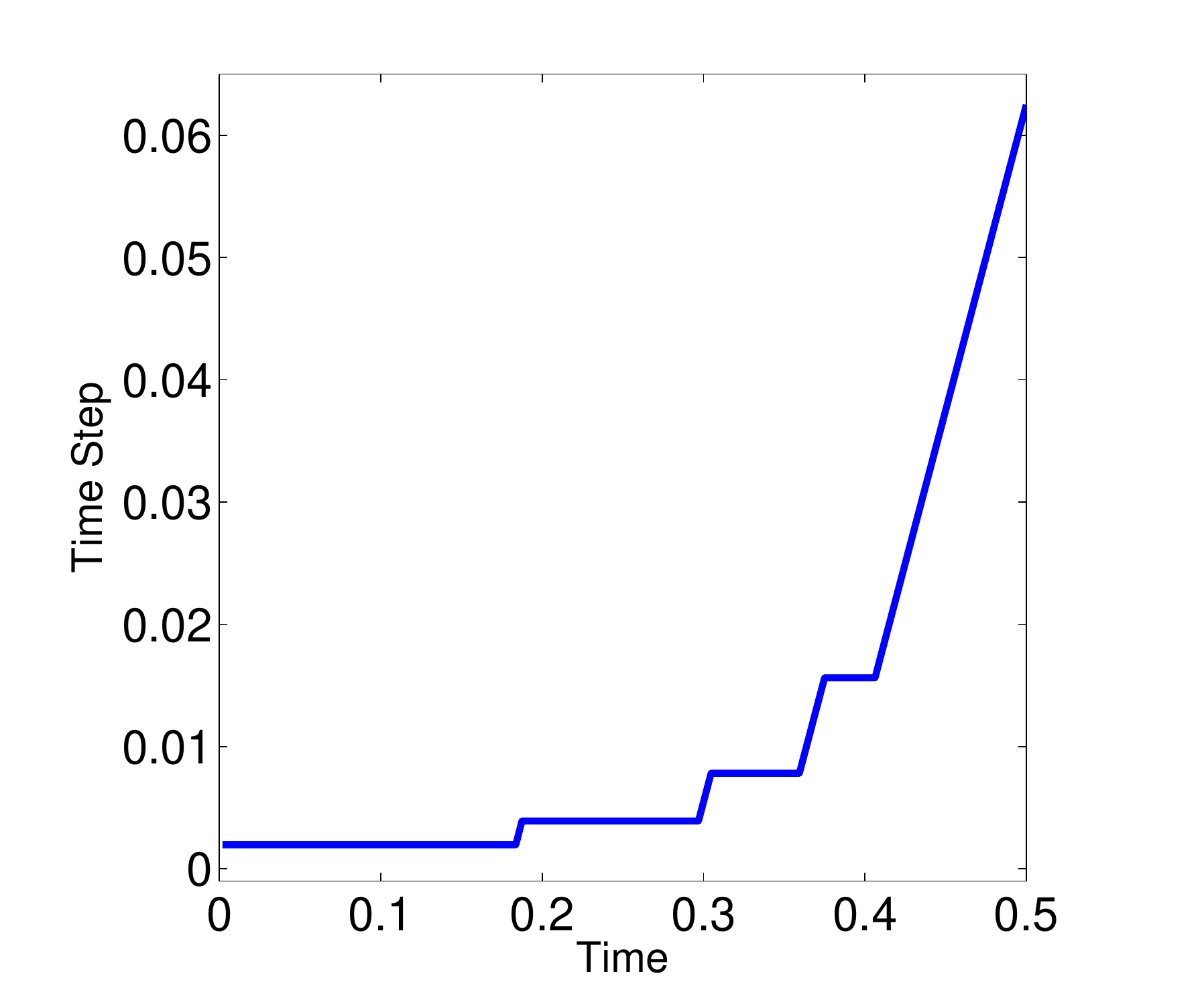}}
\caption{Example \ref{ex1}: Evolution of DoFs (left) and time-steps $\tau$ (right)\label{ex1:dofdt} }
\end{figure}

Finally, we apply the time-space adaptive algorithm  with the tolerances $\mathbf{ttol}=10^{-3}$, $\mathbf{stol^+}=3\times 10^{-4}$ and $\mathbf{stol^-}=3\times 10^{-7}$. Firstly, we prepare an initial mesh starting from a very coarse spatial mesh and a uniform partition of the time interval $[0,0.5]$ with the step-size $\tau =0.25$ until the the user defined tolerances $\mathbf{ttol}$ and $\mathbf{stol^+}$ are satisfied. The adaptive mesh at the final time $T=0.5$ is shown in Fig.~\ref{ex1:mesh}. In Fig.~\ref{ex1:dofdt} on the right, the change of the time-steps is shown, whereas the change in the DoFs is illustrated in Fig.~\ref{ex1:dofdt} on the left.  Since the layers in the problem do not move as the time progresses, the number of DoFs increases monotonically by the spatial grid refinement. In Fig.~\ref{ex1:plot}, it is shown that all the oscillations are damped out by adaptive algorithm using less DoFs compared to the uniform one.

\subsection{Coupled example with polynomial type non-linearity}
\label{ex2}

\begin{figure}[ht]
\centering
\subfloat{\includegraphics[width=0.5\textwidth]{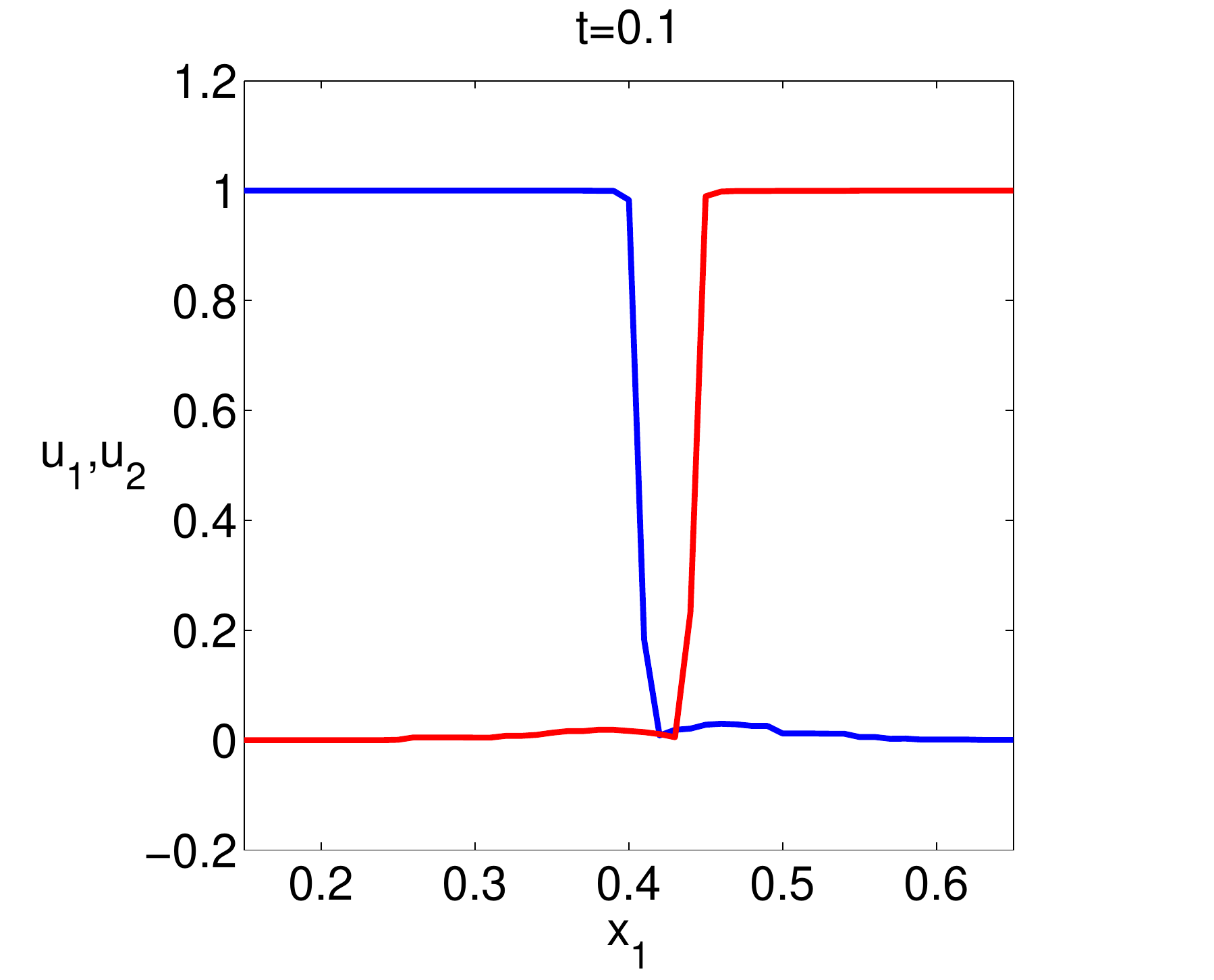}}
\subfloat{\includegraphics[width=0.5\textwidth]{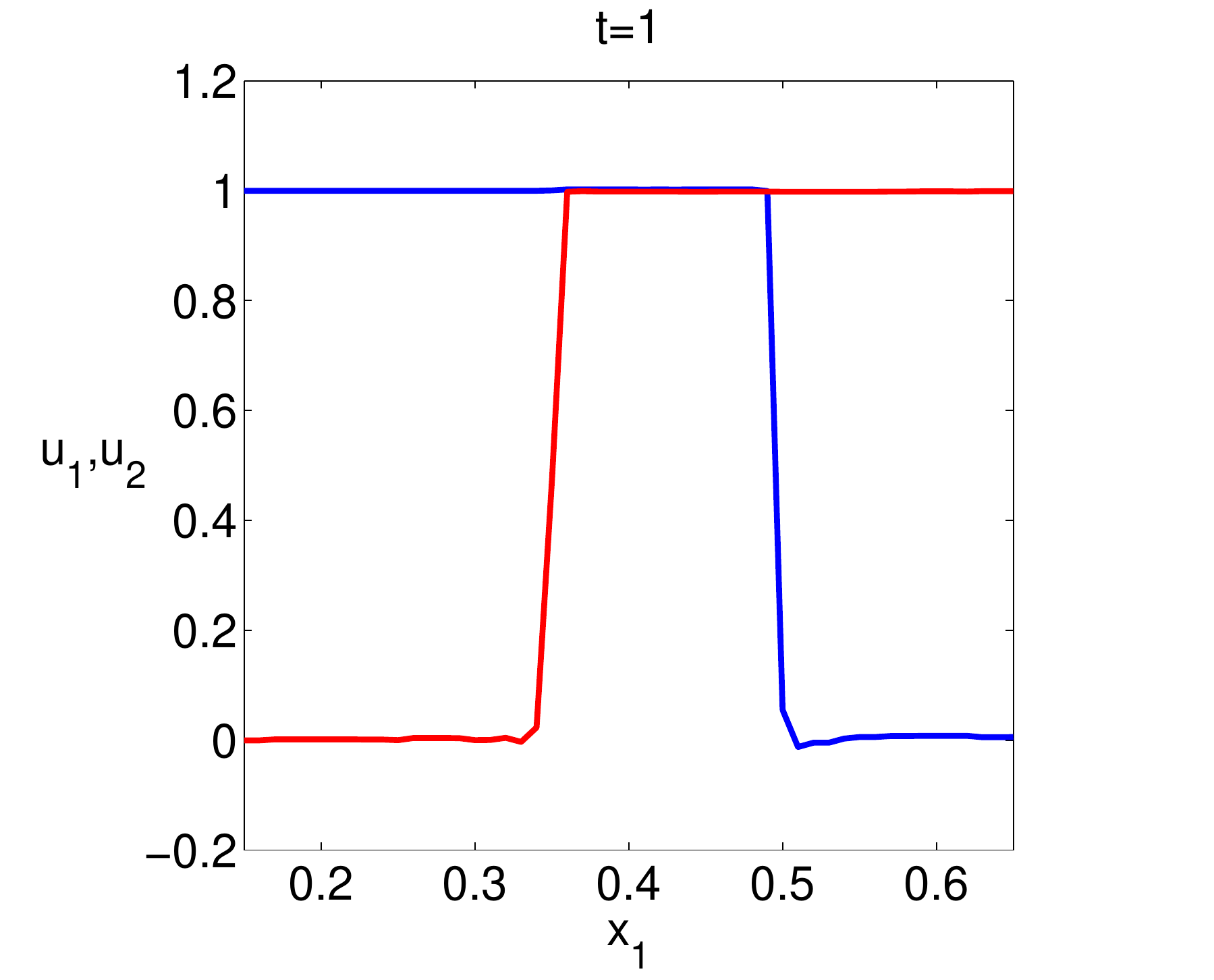}}
\caption{Example \ref{ex2}: Cross-section plots in the cross-wind direction at $t=0.1$ (left) and $t=1$ (right)\label{ex2:cross}}
\end{figure}

\begin{figure}[ht]
\centering
\subfloat{\includegraphics[width=0.42\textwidth]{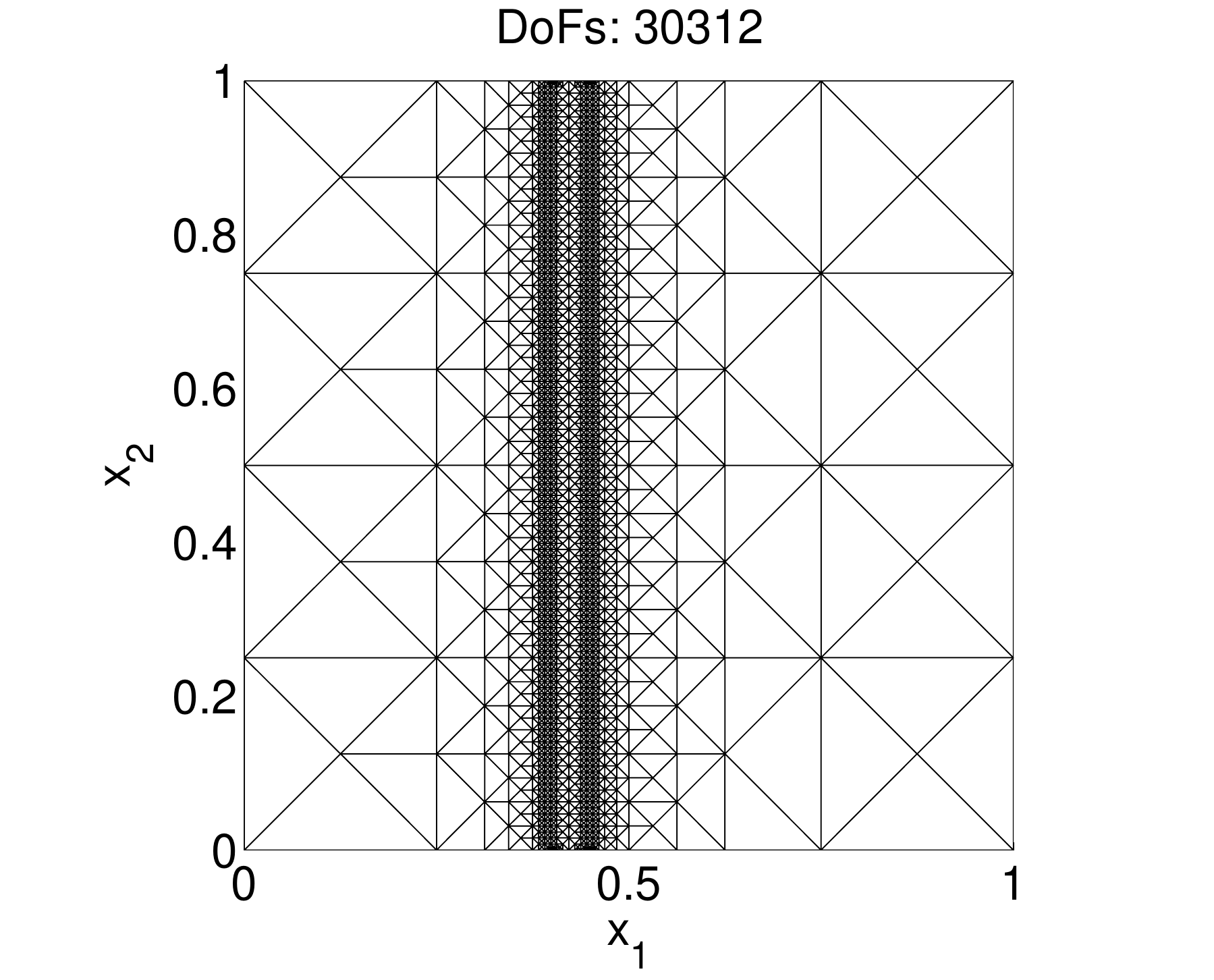}}
\subfloat{\includegraphics[width=0.42\textwidth]{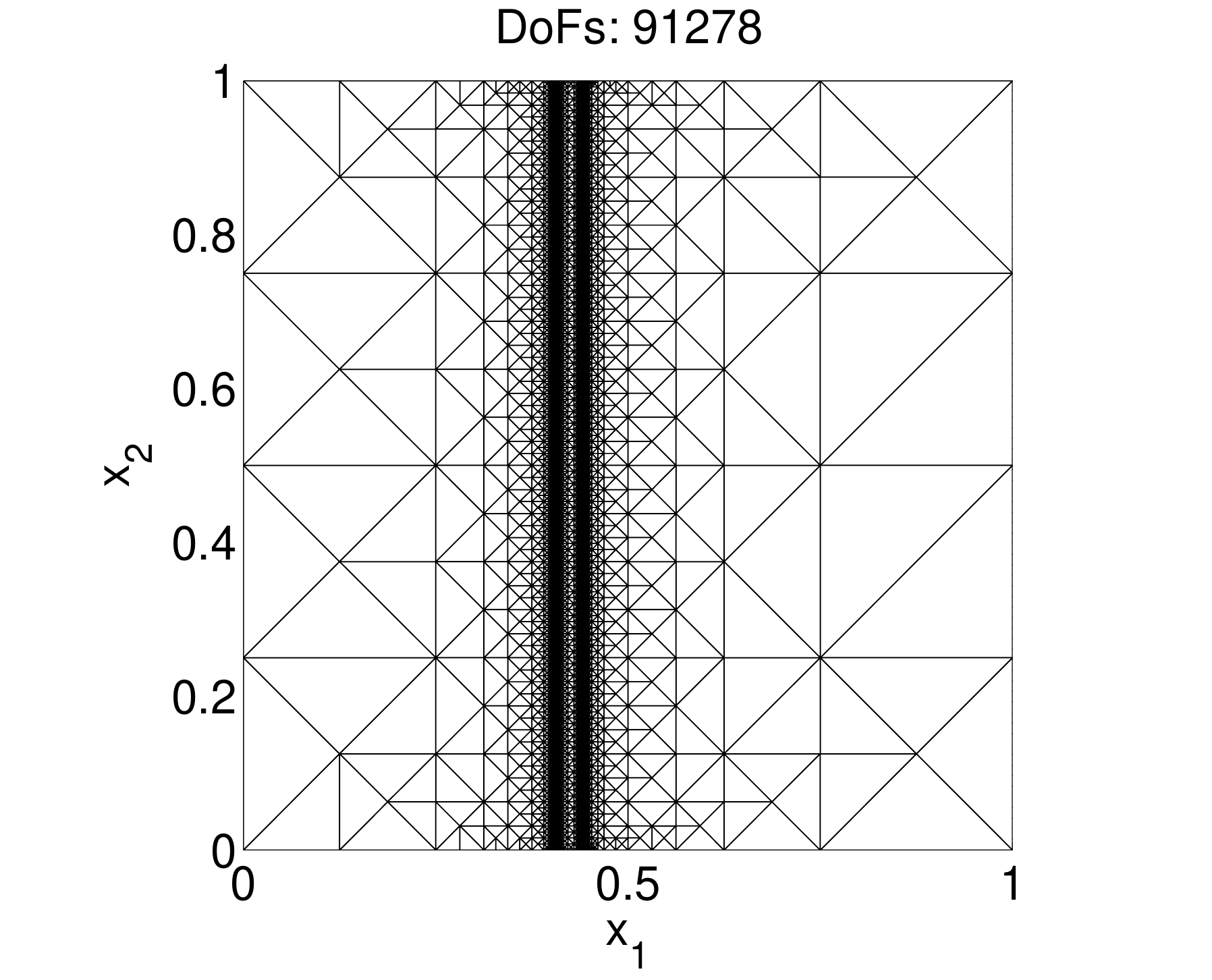}}
\subfloat{\includegraphics[width=0.42\textwidth]{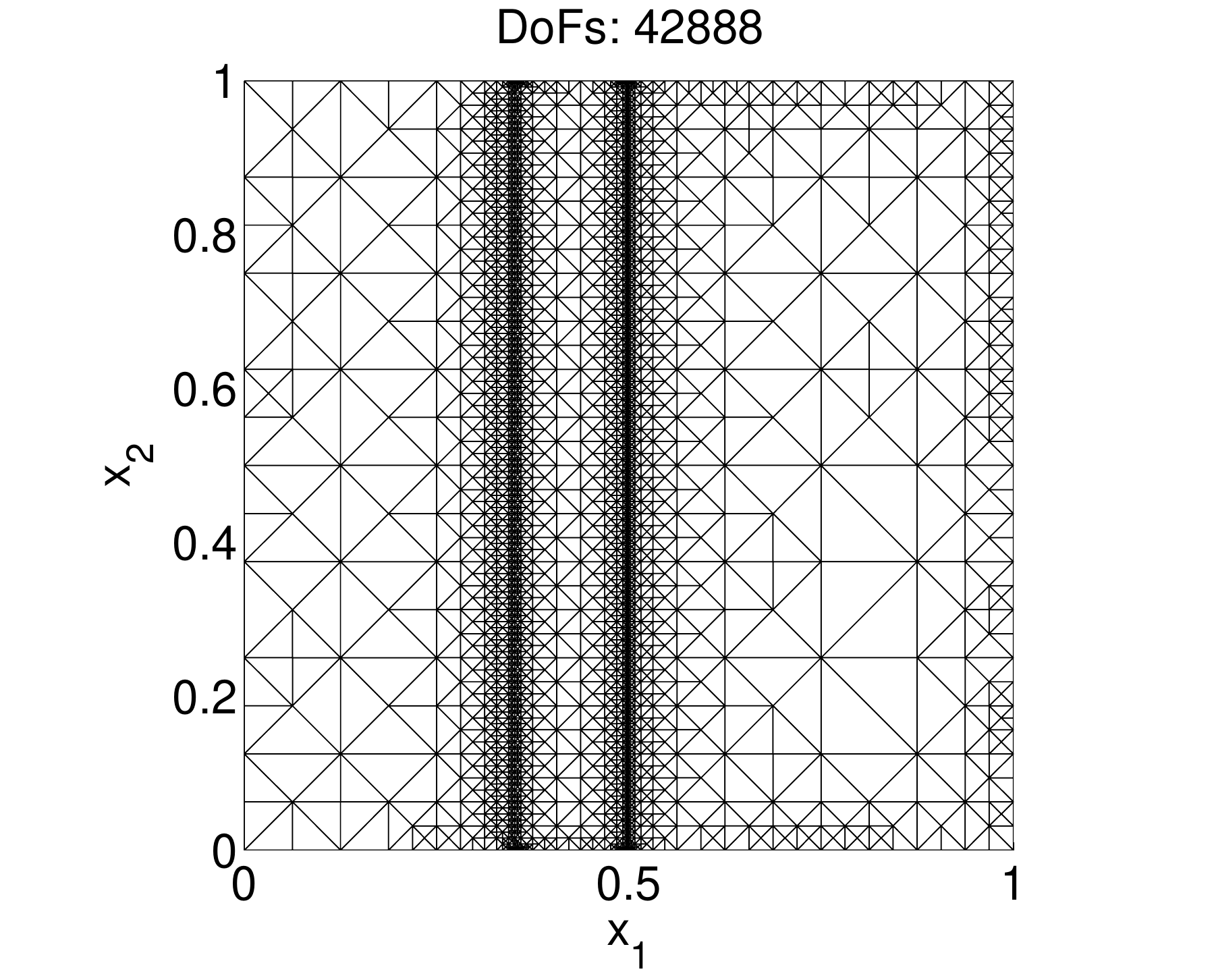}}
\caption{Example \ref{ex2}: Adaptive meshes at $t=0$, $t=0.1$ and $t=1$ (from left to right)\label{ex2:mesh}}
\end{figure}

\begin{figure}[ht]
\centering
\subfloat{\includegraphics[width=0.5\textwidth]{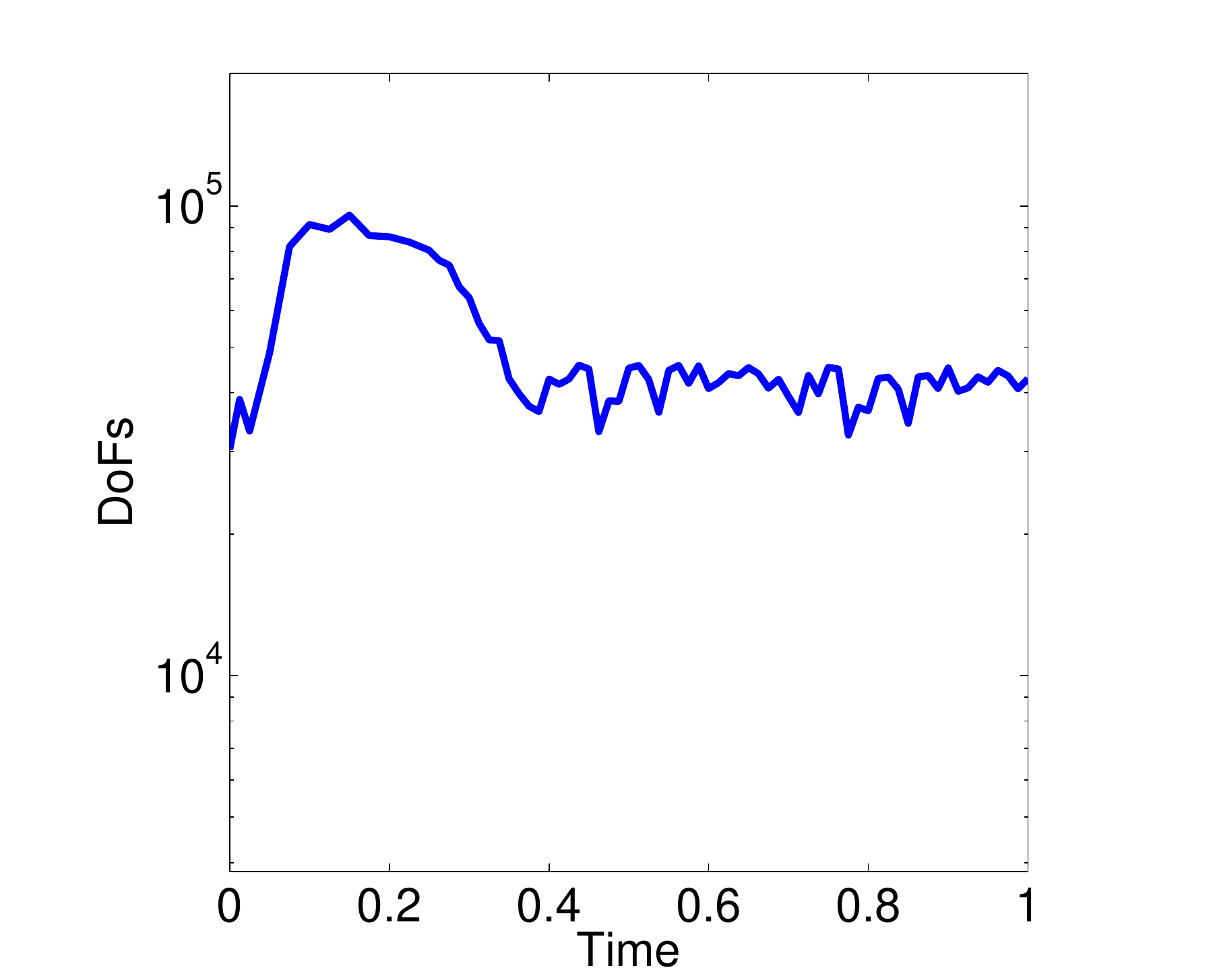}}
\subfloat{\includegraphics[width=0.5\textwidth]{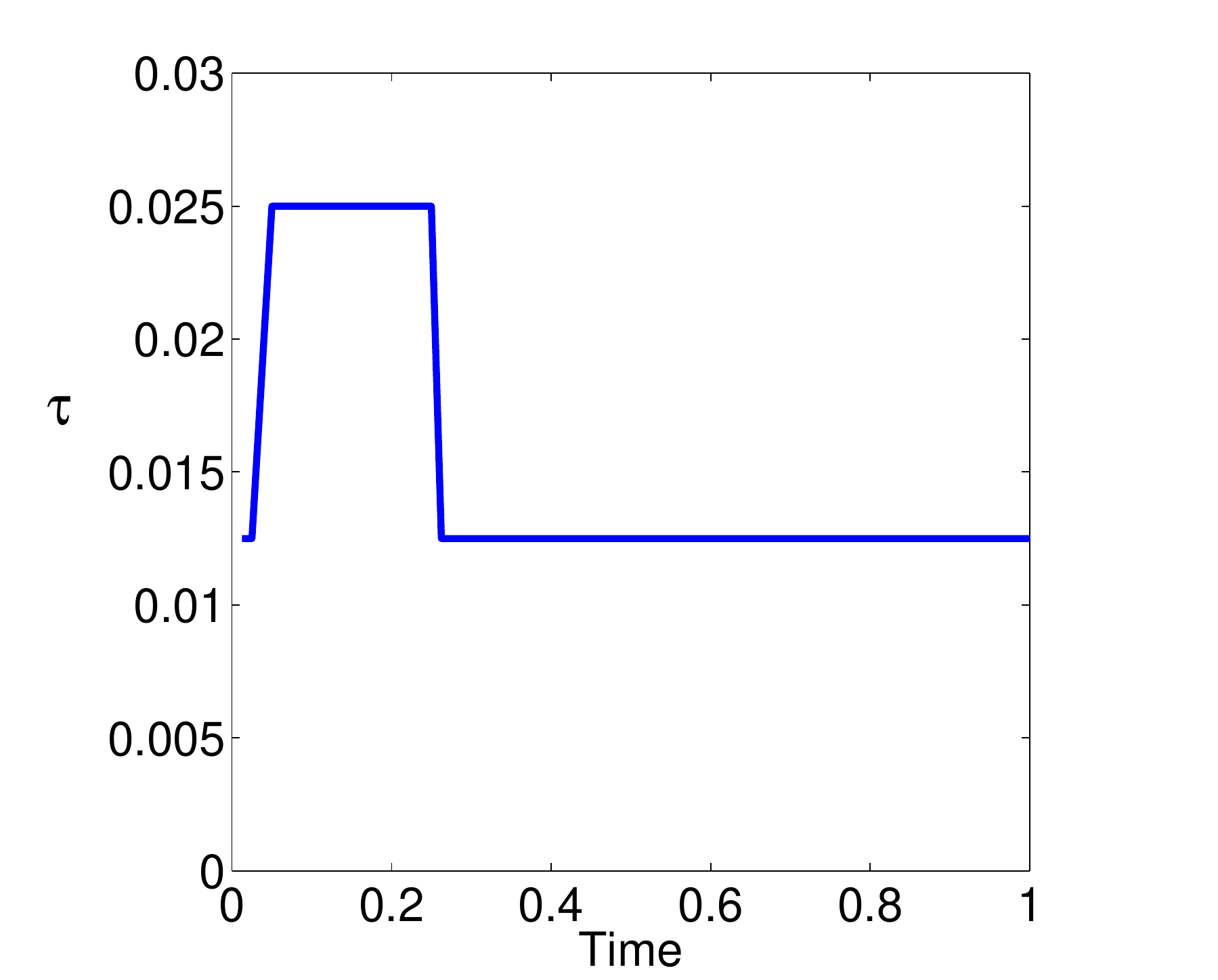}}
\caption{Example \ref{ex2}: Evolution of DoFs (left) and time-steps $\tau$ (right)\label{ex2:dtdof}}
\end{figure}

The next example is a coupled non-linear problem taken from \cite{bause13hof}.

\begin{eqnarray*}
\frac{\partial u_i}{\partial t}- \epsilon\Delta u_i+\vec{\beta }_i\cdot\nabla u_i+u_1u_2 &=& f_i , \quad i=1,2
\end{eqnarray*}
on $\Omega =(0,1)^2$ with the convection fields $\vec{\beta }_1=(1,0)^T$ and $\vec{\beta }_2=(-1,0)^T$, and the diffusion constant $\epsilon=10^{-5}$. The Dirichlet boundary conditions, initial conditions and the load functions $f_i$ are chosen so that the exact solutions are
$$
u_1(\vec{x},t)=\frac{1}{2}\left( 1-\tanh \frac{2x_1-0.2t-0.8}{\sqrt{5\epsilon }}\right)
$$
$$
u_2(\vec{x},t)=\frac{1}{2}\left( 1+\tanh \frac{2x_1+0.2t-0.9}{\sqrt{5\epsilon }}\right)
$$
We use again quadratic DG elements.
We prepare an initial mesh, Fig.~\ref{ex2:mesh} on the  left, starting with a very coarse spatial mesh and a uniform partition of the time interval $[0,1]$ with the step-size $\tau =0.1$ until the user defined tolerances $\mathbf{ttol}=10^{-3}$ and $\mathbf{stol^+}=10^{-1}$ are satisfied. Here, two sharp fronts move towards to each other and then mix directly after the time $t=0.1$, Fig.~\ref{ex2:cross}. The movement of the fronts are also visible in Fig.~\ref{ex2:mesh} claiming that refinement/coarsening of the adaptive algorithm works well. We see that the sharp fronts in the cross-wind direction $x_2=0.5x_1+0.75$ are almost damped out. Moreover, Fig.~\ref{ex2:mesh}-\ref{ex2:dtdof} show that both the spatial and temporal estimators catch the time where the two sharp fronts mix.

\subsection{Non-linear ADR equation in homogeneous porous media}
\label{ex3}

\begin{figure}[ht]
\centering
\subfloat{\includegraphics[width=0.42\textwidth]{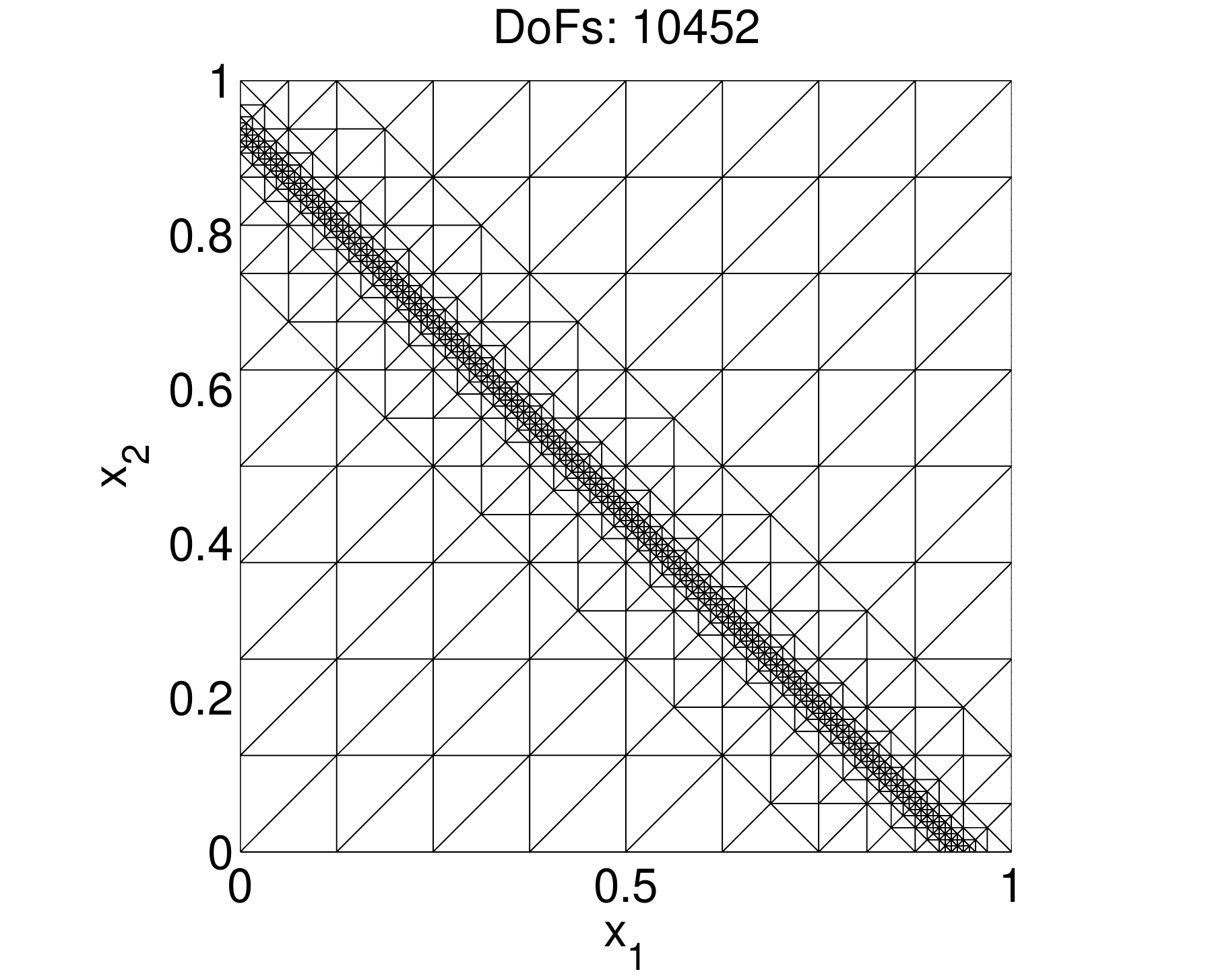}}
\subfloat{\includegraphics[width=0.42\textwidth]{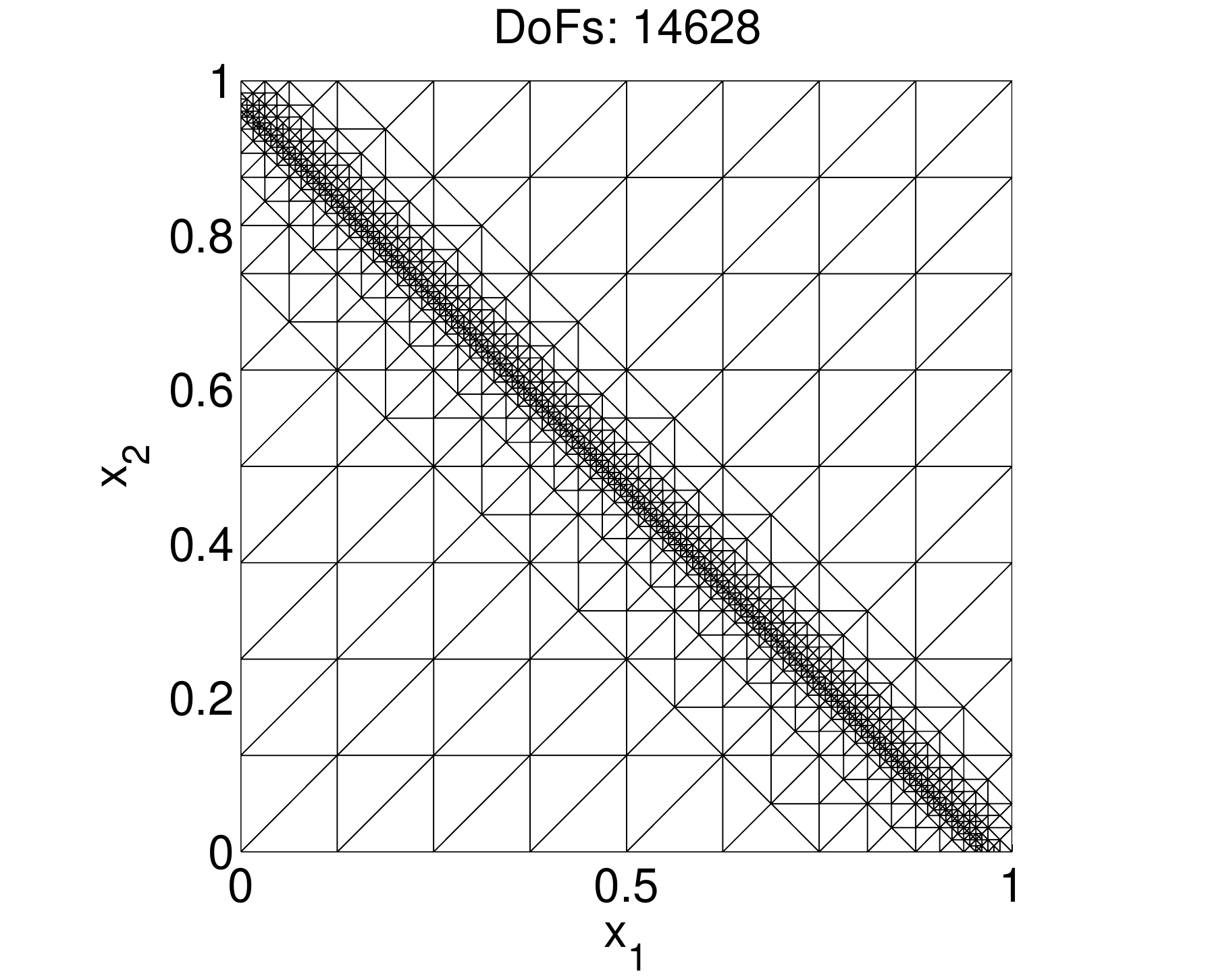}}
\subfloat{\includegraphics[width=0.42\textwidth]{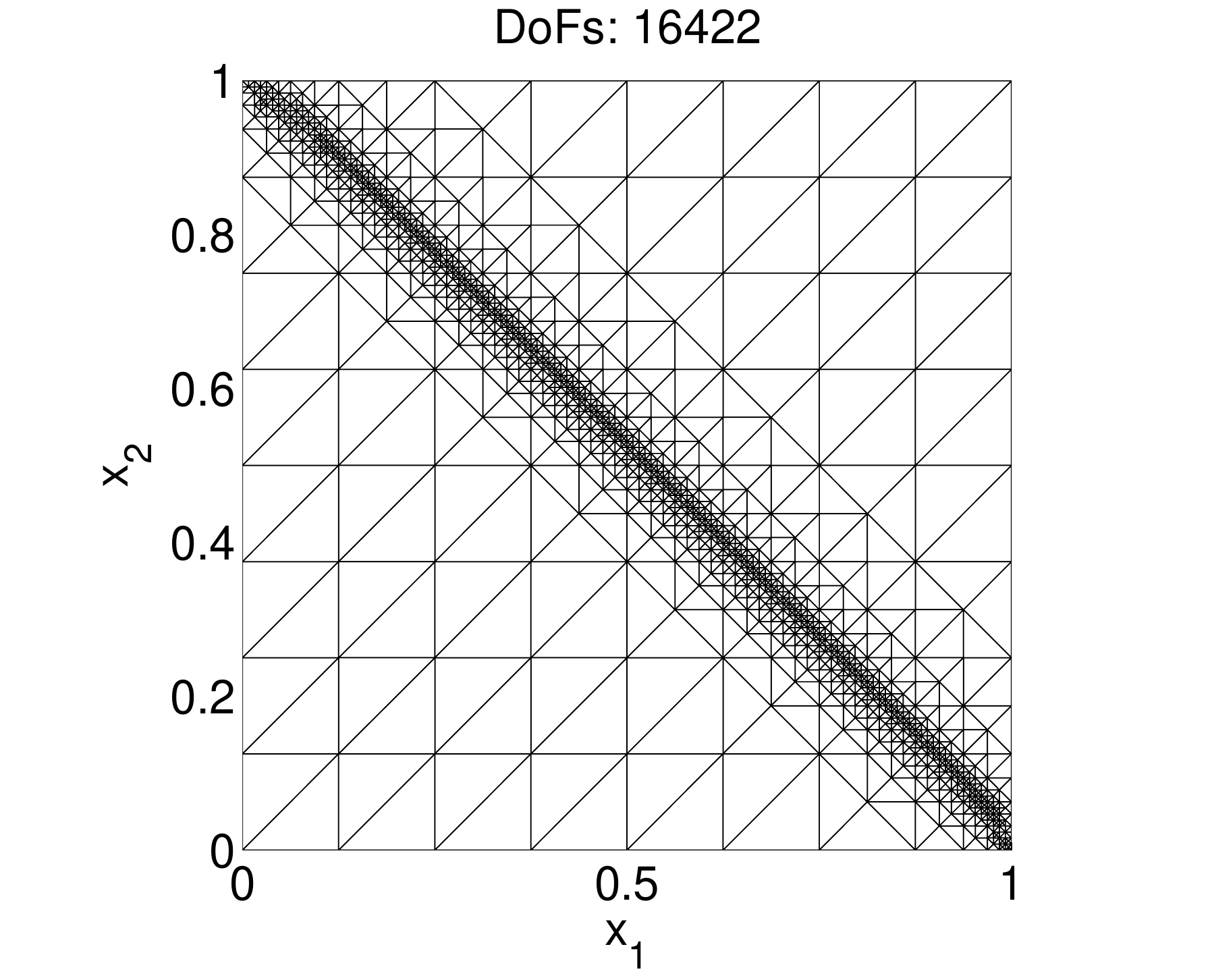}}

\subfloat{\includegraphics[width=0.42\textwidth]{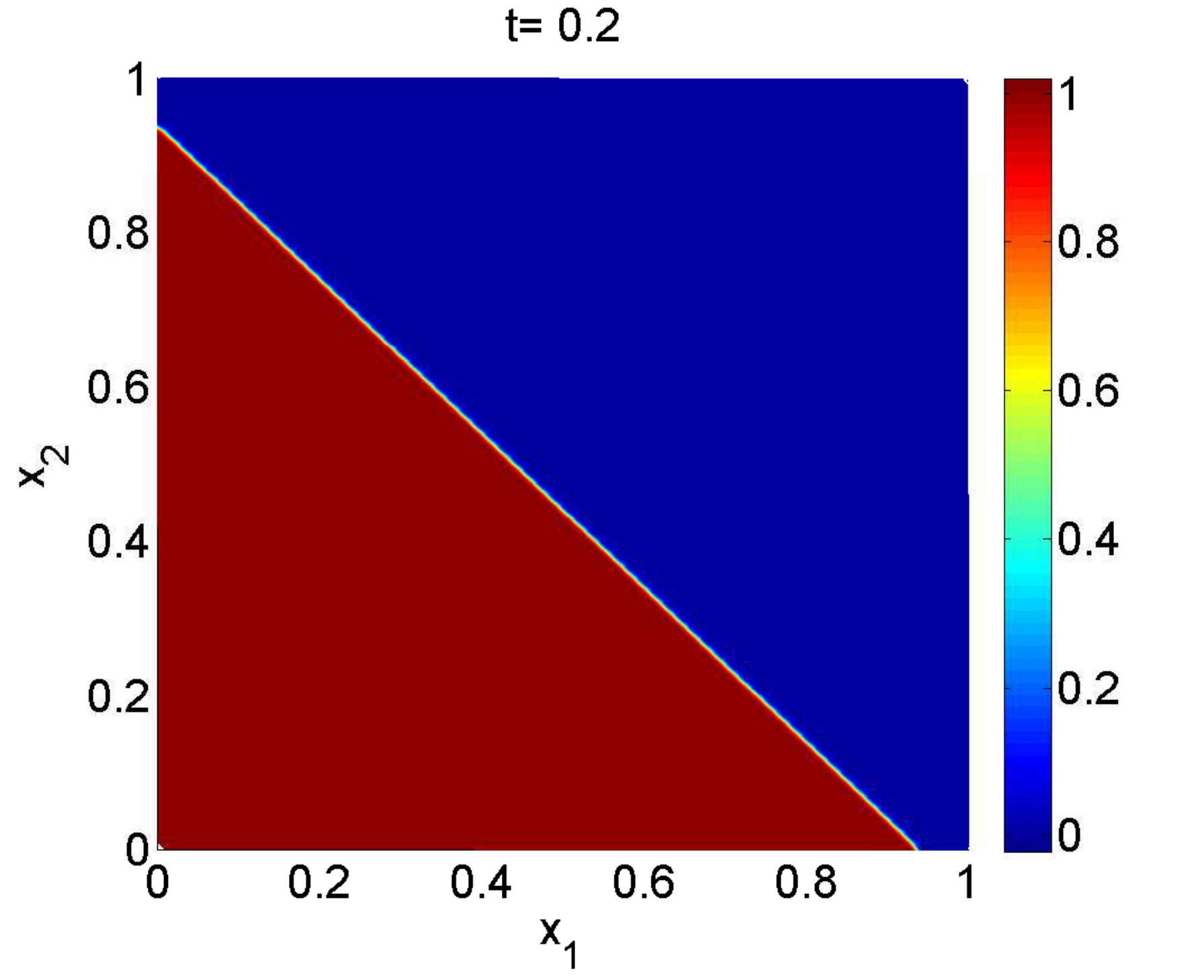}}
\subfloat{\includegraphics[width=0.42\textwidth]{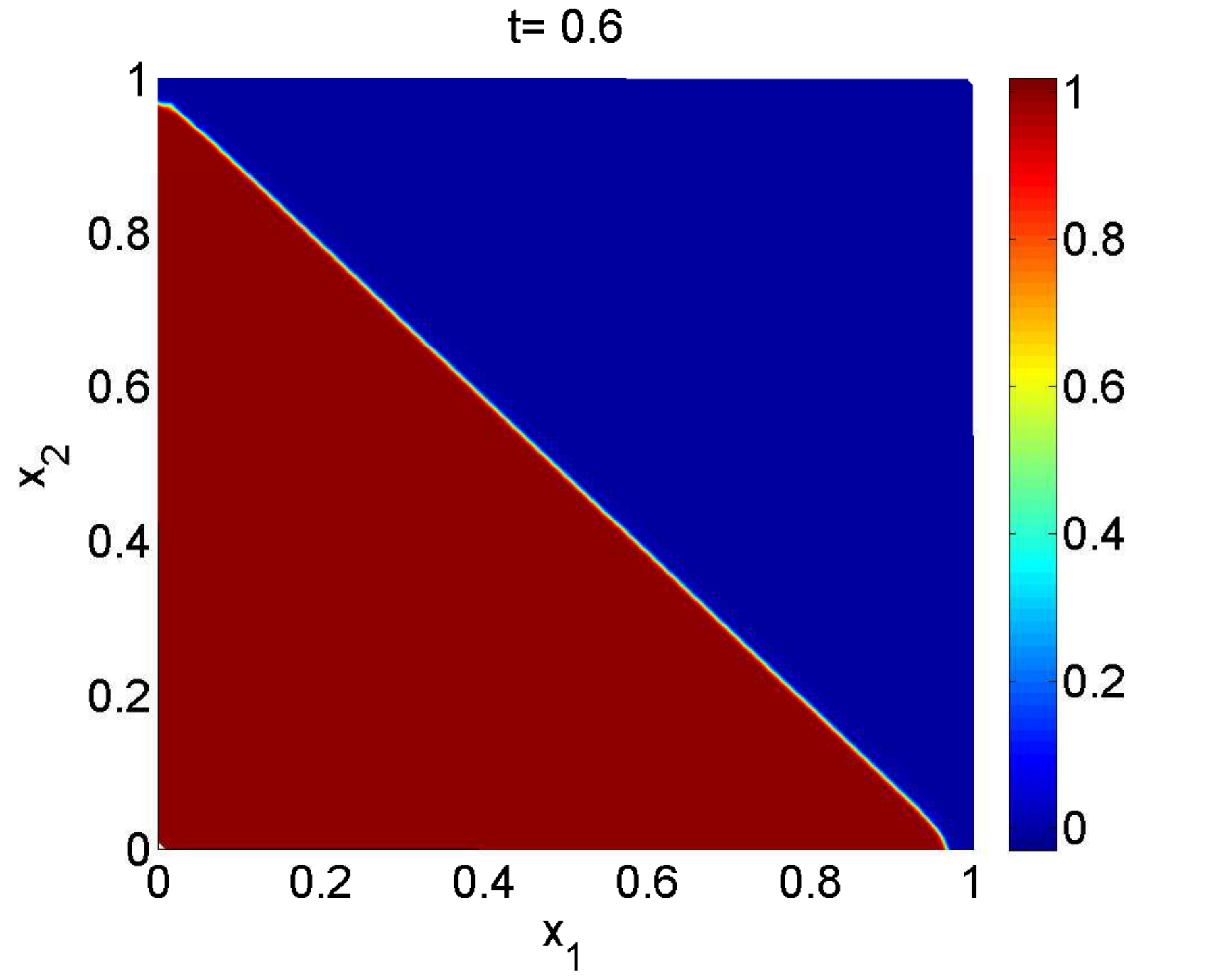}}
\subfloat{\includegraphics[width=0.42\textwidth]{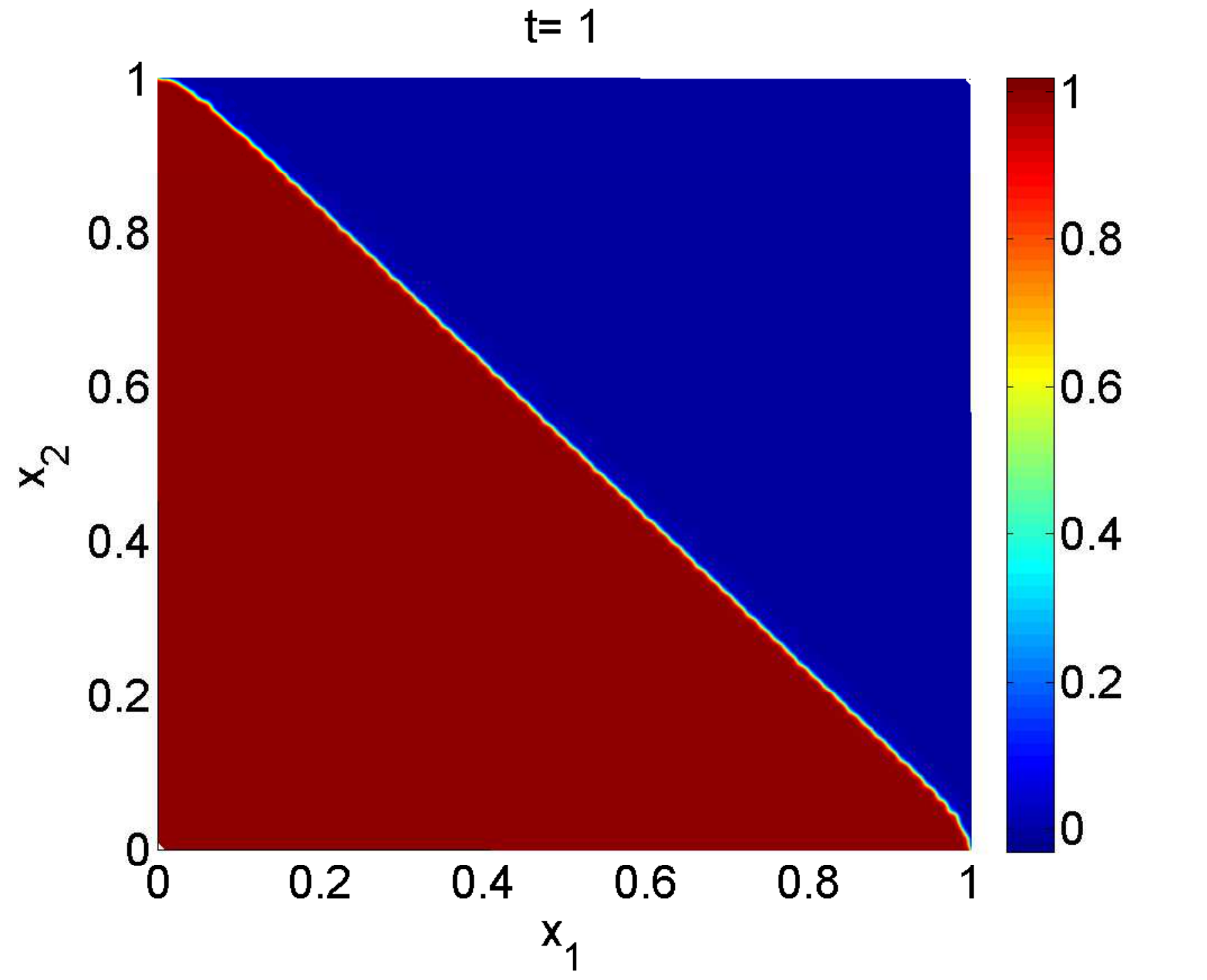}}
\caption{Example \ref{ex3}: Adaptive meshes (top) and solution profiles (bottom) at $t=0.2$, $t=0.6$ and $t=1$ (from left to right)\label{ex3:meshsol}}
\end{figure}

\begin{figure}[ht]
\centering
\subfloat{\includegraphics[width=0.5\textwidth]{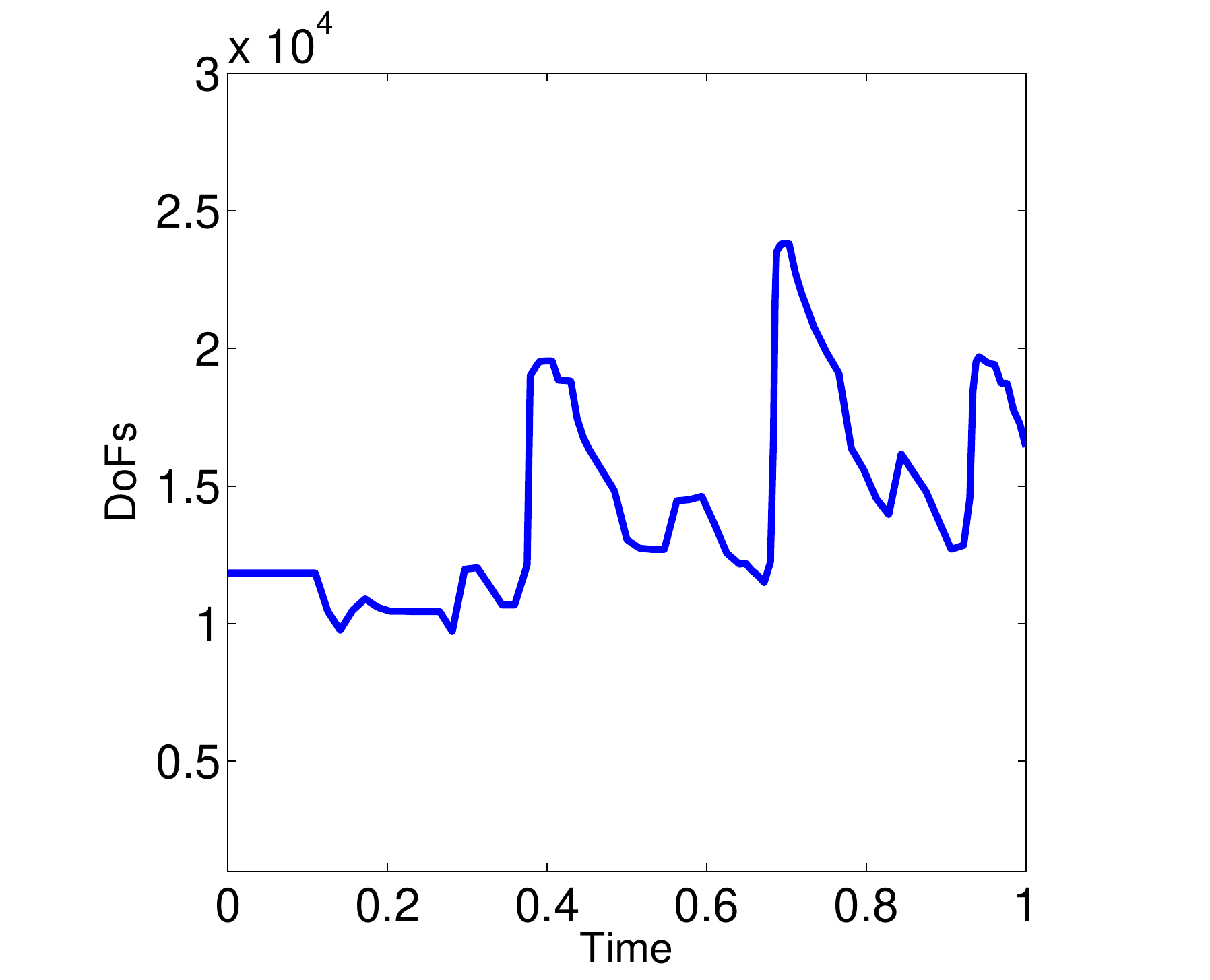}}
\subfloat{\includegraphics[width=0.5\textwidth]{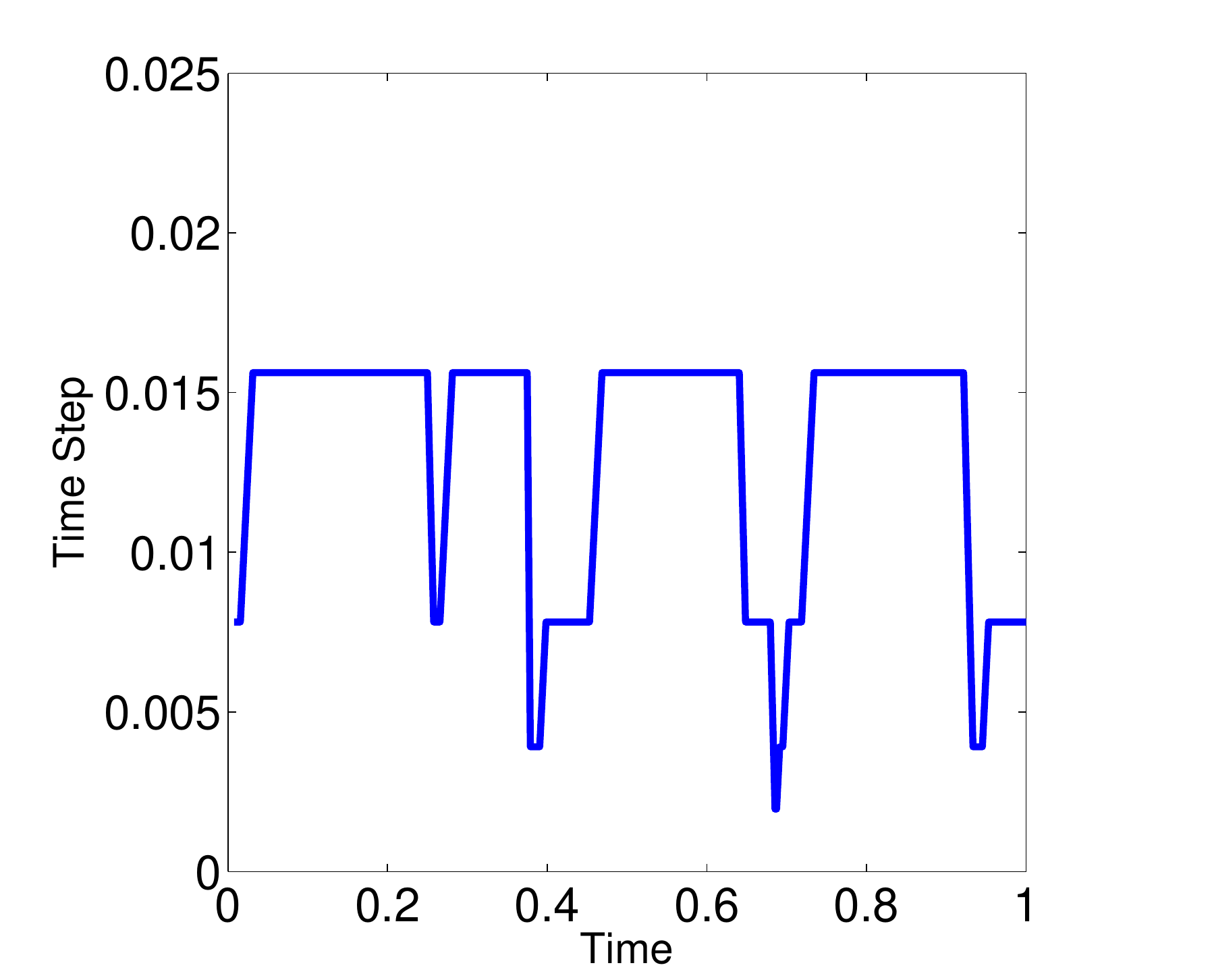}}
\caption{Example \ref{ex3}: Evolution of DoFs (left) and time-steps $\tau$ (right)\label{ex3:dtdof}}
\end{figure}

We consider the advection-diffusion-reaction (ADR) equation in \cite{tambue10eia} with polynomial type non-linear reaction
\begin{eqnarray*}
\frac{\partial u}{\partial t}- \epsilon\Delta u + \vec{\beta}\cdot\nabla u+\gamma u^2(u-1) &=& 0 \quad \text{in } \; \Omega\times (0,T]
\end{eqnarray*}
on $\Omega =(0,1)^2$. We take as in \cite{tambue10eia} the homogeneous dispersion tensor as $\epsilon =10^{-4}$, the velocity field $\vec{\beta}=(-0.01,-0.01)^T$ and $\gamma =100$. The initial and boundary conditions are chosen by the exact solution
$$
u(\vec{x},t)=[ 1+\exp ( a(x_1+x_2-bt)+a(b-1)) ]^{-1}
$$
with $a=\sqrt{\gamma /(4\epsilon)}$ and $b=-0.02+\sqrt{\gamma\epsilon}$. The problem is a transport of a front in homogeneous porous media. We simulate the given problem for the final time $T=1$, and with quadratic DG elements. We begin by preparing an initial mesh starting from a very coarse spatial mesh and a uniform partition of the time interval $[0,1]$ with the step-size $\tau =0.25$ until the user defined tolerances $\mathbf{ttol}=3\times 10^{-3}$ and $\mathbf{stol^+}=10^{-3}$ are satisfied. In  Fig.~\ref{ex3:meshsol}, the adaptive meshes and solution profiles are shown  at the times $t=\{ 0.2,0.6,1\}$ where the movement of the front can be seen. The time vs DoFs and time vs time step-size plots in Fig.~\ref{ex3:dtdof}
indicate clearly the oscillations in DoFs and time-steps due to the movement of the front.

\subsection{Non-linear ADR in deterministic heterogeneous porous media}
\label{ex4}

We consider the ADR equation in \cite{tambue10eia} with Monod or Langmuir isotherm type non-linear reaction
\begin{eqnarray*}
\frac{\partial u}{\partial t}- \nabla\cdot (\epsilon\nabla u) + \vec{\beta}(x)\cdot\nabla u+\frac{u}{1+u} &=& 0 \quad \text{in } \; \Omega\times (0,T]\\
u(x,t) &=& 1 \quad \text{on } \; \Gamma^D\times [0,T] \\
-\epsilon\nabla u(x,t)\cdot\vec{n} &=& 0 \quad \text{on } \; (\partial\Omega\setminus\Gamma^D)\times [0,T] \\
u(x,0) &=& 0 \quad \text{in } \; \Omega
\end{eqnarray*}
with $\Omega =(0,3)\times (0,2)$ and $\Gamma^D=\{ 0\}\times [0,2]$. The problem represents a reaction in porous media, for instance, transport in a highly idealized fracture pattern. Here $\epsilon$ stands for the heterogeneous dispersion tensor given by
$$
\epsilon =
\begin{bmatrix}
10^{-3} & 0\\
0 & 10^{-4}
\end{bmatrix}
$$
The velocity field $\vec{\beta}(x)$ is determined via the Darcy's law
$$
\vec{\beta} = -\frac{k(x)}{\mu}\nabla p
$$
where $p$ is the fluid pressure, $\mu$ is the fluid viscosity and $k(x)$ is the permeability of the porous medium. Using the mass conservation property $\nabla\cdot\vec{\beta}(x)=0$ under the assumption that rock and fluids are incompressible, the velocity field $\vec{\beta}(x)$ is computed by solving
\begin{eqnarray*}
\nabla\cdot\left( \frac{k(x)}{\mu}\nabla p \right) &=& 0 \quad \text{in } \; \Omega\\
p &=& 1 \quad \text{on } \; \{ 0\}\times [0,2] \\
p &=& 0 \quad \text{on } \; \{ 3\}\times [0,2] \\
-k(x)\nabla p\cdot \vec{n} &=& 0 \quad \text{on } \; (0,3)\times \{ 0,2\}
\end{eqnarray*}

\begin{figure}[ht]
\centering
\subfloat{\includegraphics[width=0.5\textwidth]{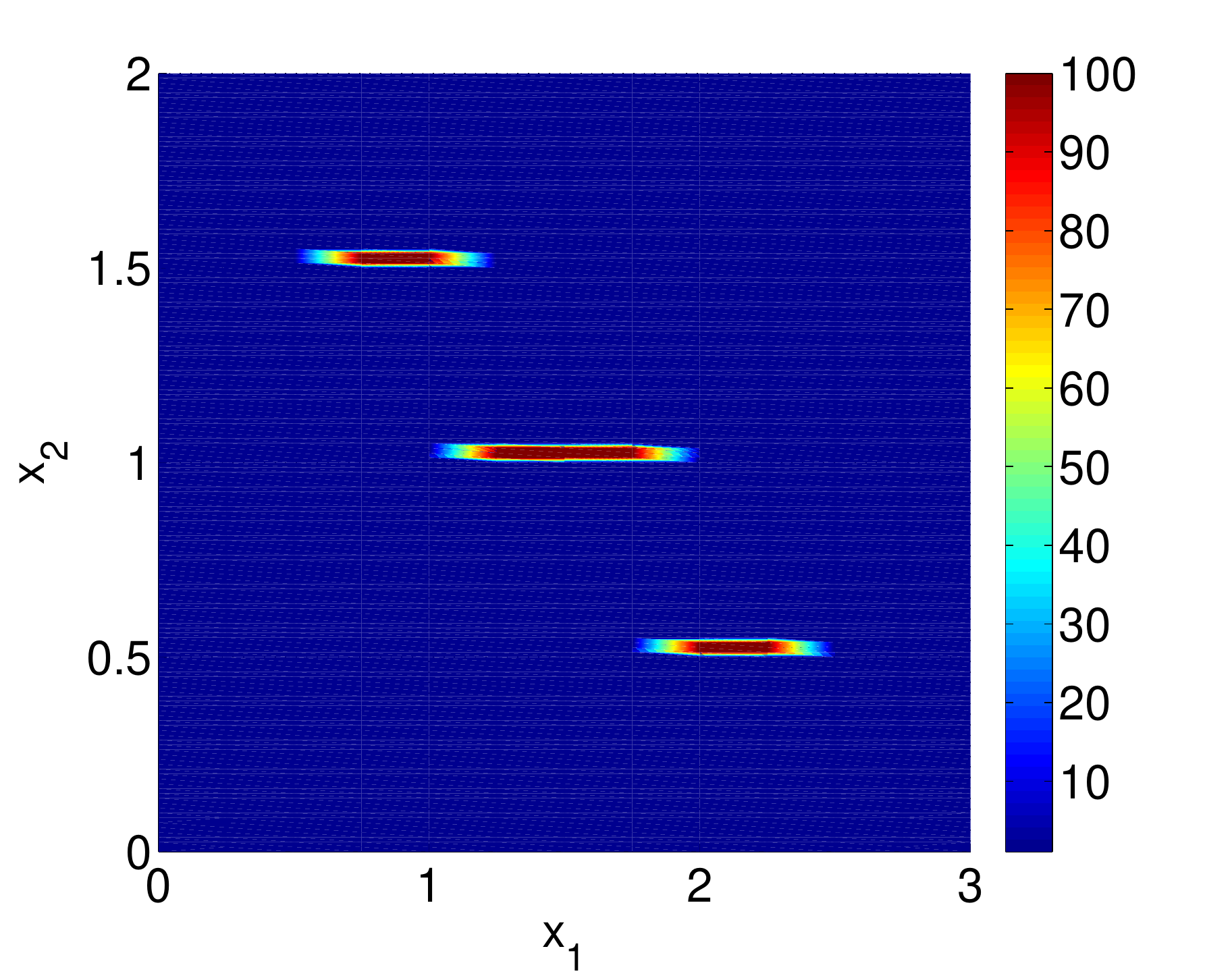}}
\subfloat{\includegraphics[width=0.5\textwidth]{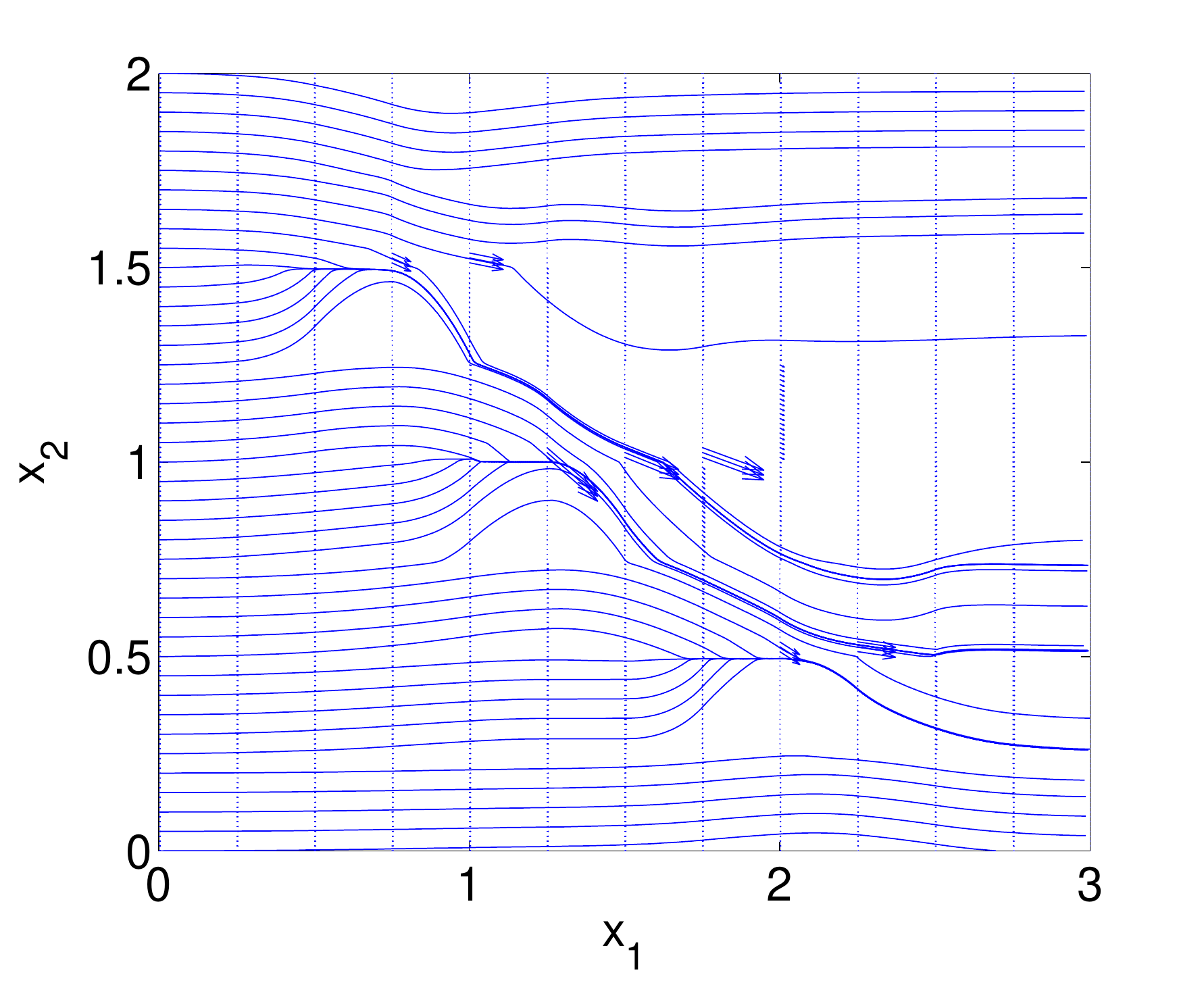}}
\caption{Example \ref{ex4}: Permeability field (left) and velocity streamlines (right)\label{ex4:perm}}
\end{figure}

\begin{figure}[ht]
\centering
\includegraphics[width=1\textwidth]{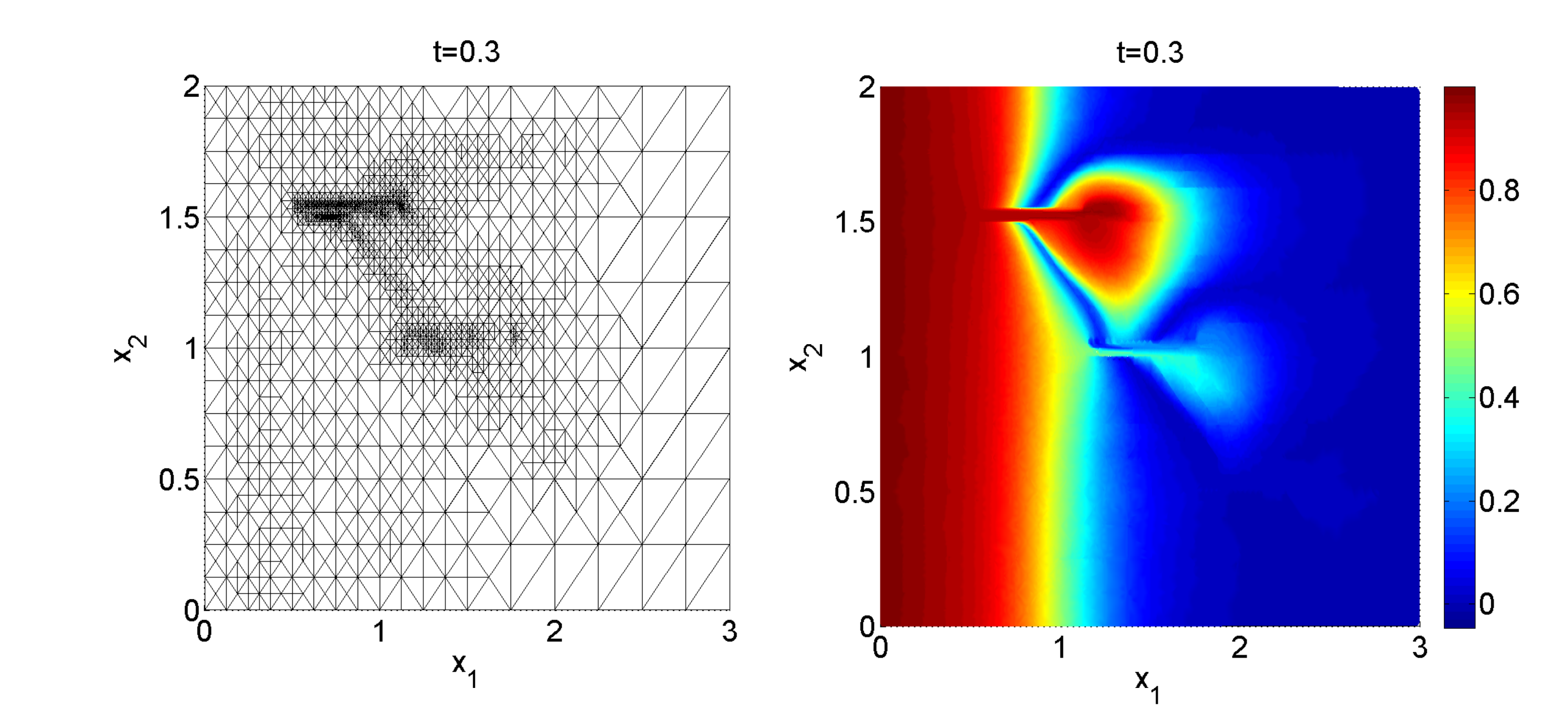}
\caption{Example \ref{ex4}: Adaptive mesh (left) and concentration (right) at $t=0.3$\label{ex4:plot03}}
\end{figure}

\begin{figure}[ht]
\centering
\includegraphics[width=1\textwidth]{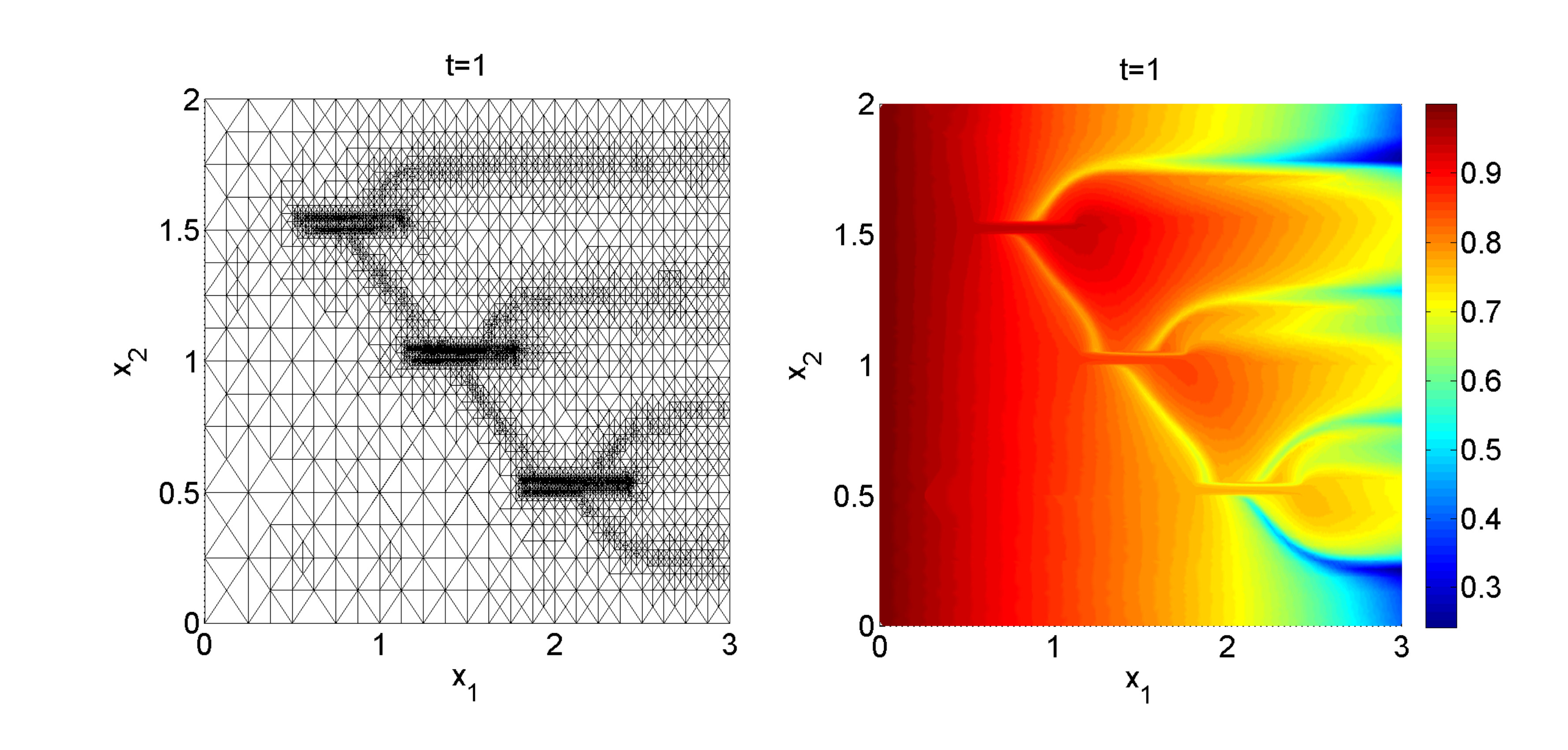}
\caption{Example \ref{ex4}: Adaptive mesh (left) and concentration (right) at $t=1$\label{ex4:plot1}}
\end{figure}

We simulate the given problem for the final time $T=1$ using linear DG elements. We take the fluid viscosity $\mu =0.1$, and the permeability field as in \cite{tambue10eia} with three parallel streaks the permeability of which are 100 times greater than the permeability of the surrounding domain, see Fig.~\ref{ex4:perm} on the left, by which the flow is canalized from the lower-permeability rocks into the high-permeability ones, Fig.~\ref{ex4:perm} on the right. For the adaptive procedure, we prepare an initial mesh starting from a very coarse spatial mesh and a uniform partition of the time interval $[0,1]$ with the step-size $\tau =0.05$ until the user defined tolerances $\mathbf{ttol}=10^{-3}$ and $\mathbf{stol^+}=3\times 10^{-4}$ are satisfied. Fig.~\ref{ex4:plot03}-\ref{ex4:plot1} show the adaptive meshes and concentrations at $t=0.3$ and $t=1$, where we can clearly see the flow-focusing due to the high-permeability.

\begin{figure}[ht]
\centering
\subfloat{\includegraphics[width=0.5\textwidth]{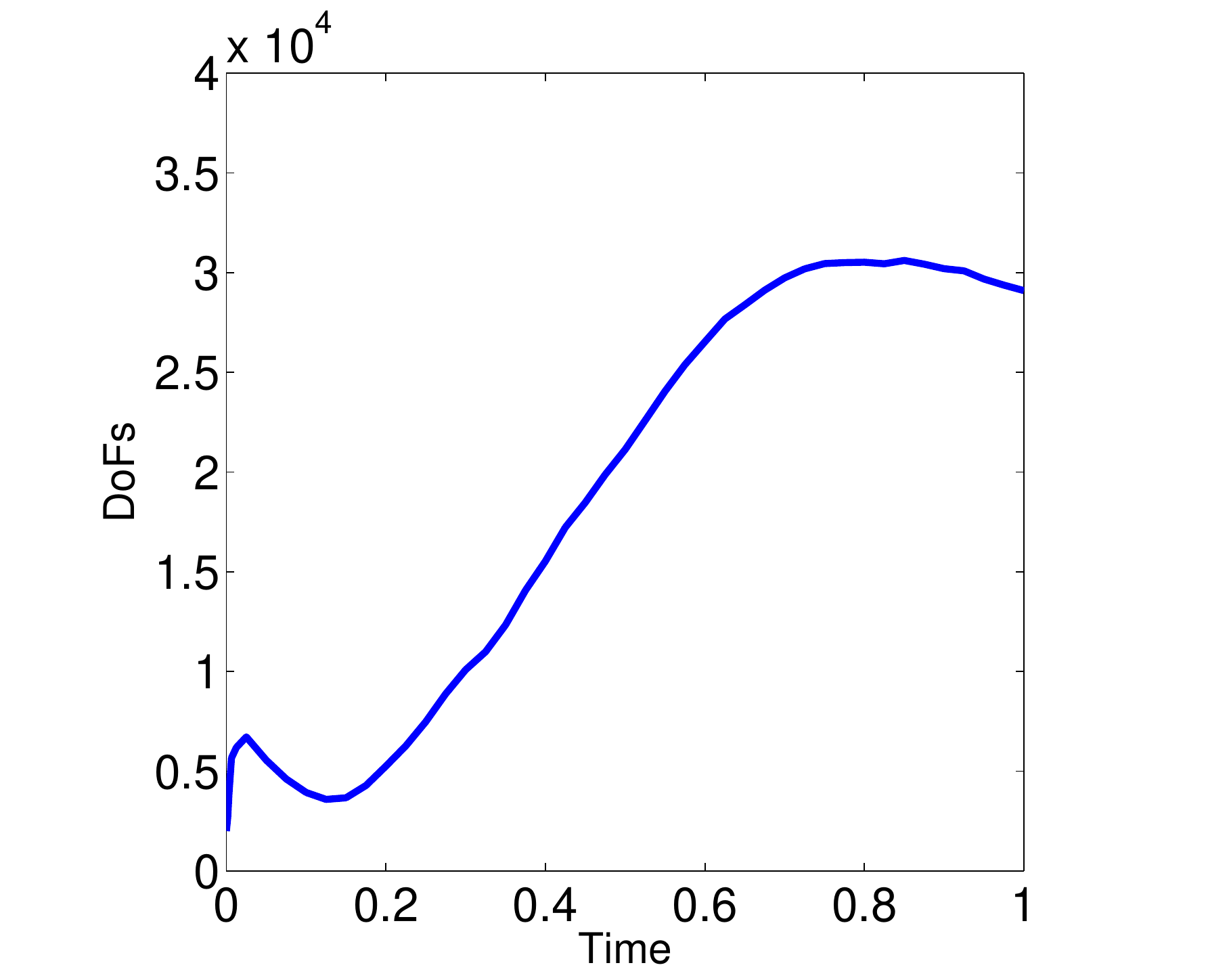}}
\subfloat{\includegraphics[width=0.5\textwidth]{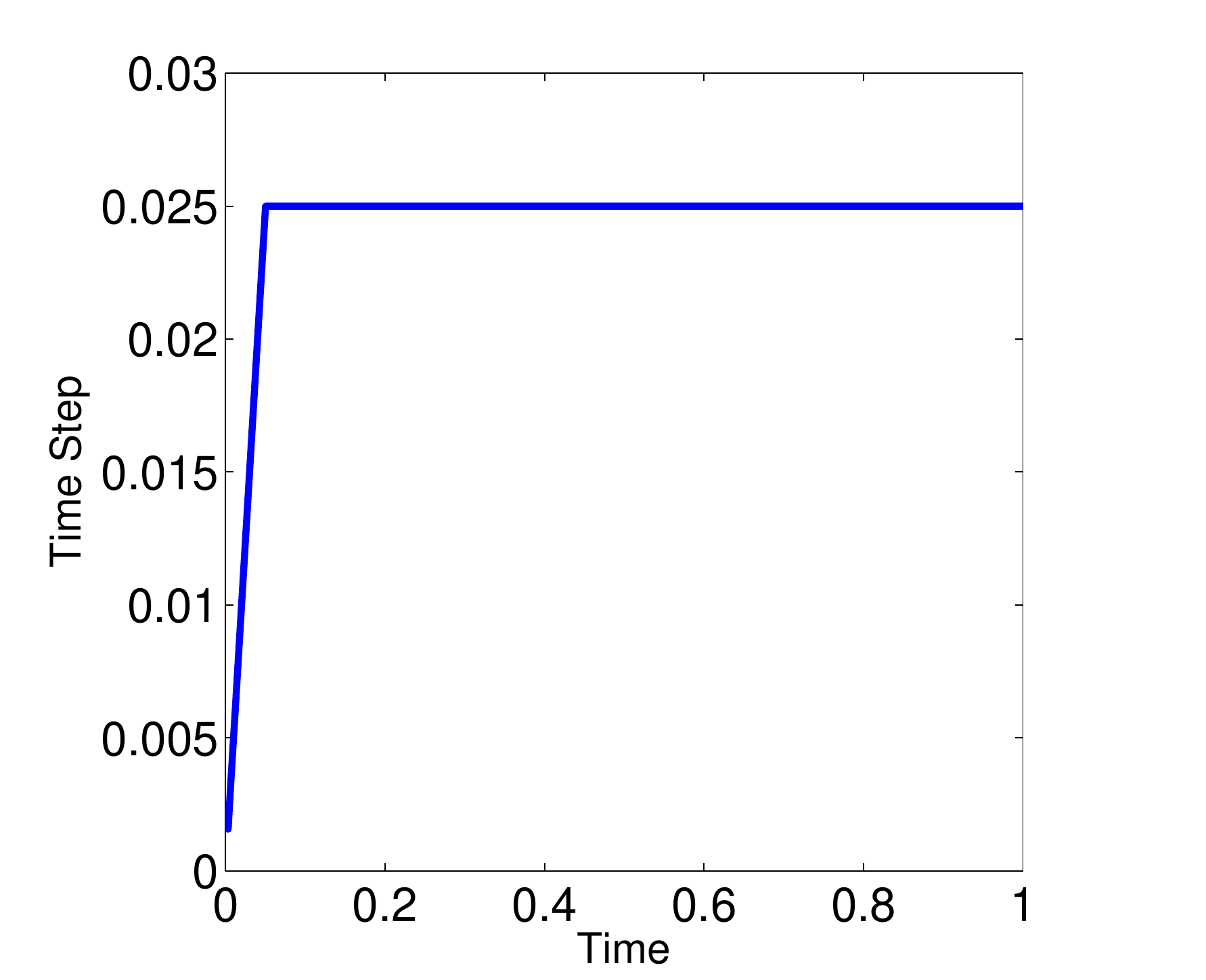}}
\caption{Example \ref{ex4}: Evolution of DoFs (left) and time-steps $\tau$ (right)\label{ex4:dtdof}}
\end{figure}

Time vs DoFs and time vs time step-size plots are given in Fig.~\ref{ex4:dtdof}. We see that initially small time steps are used and then it reaches a steady time step, Fig.~\ref{ex4:dtdof} on the right. The number of DoFs increases (refinement dominates coarsening) monotonically after the meet of first high-permeability rock until the meet of third high-permeability rock and then the increase stops, Fig.~\ref{ex4:dtdof} on the left. This is meaningful since there is no sharp flow canalization after the third high-permeability rock.

\section{Conclusion}

We implemented a time-space adaptive algorithm for non-linear ADR equations  based on utilizing the elliptic reconstruction technique to be able to use the elliptic a posteriori error estimator for the convection dominated parabolic problems with non-linear reaction mechanisms.  We derived a posteriori error estimator in the $L^{\infty}(L^2) + L^2(H^1)$-type norm using backward Euler in time and SIPG in space. We demonstrated the performance of the algorithm by numerical studies.


\begin{thebibliography}{10}
\providecommand{\url}[1]{{#1}}
\providecommand{\urlprefix}{URL }
\expandafter\ifx\csname urlstyle\endcsname\relax
  \providecommand{\doi}[1]{DOI~\discretionary{}{}{}#1}\else
  \providecommand{\doi}{DOI~\discretionary{}{}{}\begingroup
  \urlstyle{rm}\Url}\fi

\bibitem{arnold02uad}
Arnold, D., Brezzi, F., Cockborn, B., Marini, L.: Unified analysis of
  discontinuous {G}alerkin methods for elliptic problems.
\newblock SIAM J. Numer. Anal. \textbf{39}, 1749--1779 (2002)

\bibitem{ayuso09dgm}
Ayuso, B., Marini, L.D.: Discontinuous {G}alerkin methods for
  advection-diffusion-reaction problems.
\newblock SIAM J. Numer. Anal. \textbf{47}, 1391--1420 (2009)

\bibitem{bastian11udg}
Bastian, P., Engwer, C., Fahlke, J., Ippisch, O.: An unfitted discontinuous
  {G}alerkin method for pore-scale simulations of solute transport.
\newblock Mathematics and Computers in Simulation \textbf{81}(10), 2051--2061
  (2011)

\bibitem{bause12ash}
Bause, M., Schwegler, K.: Analysis of stabilized higher-order finite element
  approximation of nonstationary and non-linear convection-diffusion-reaction
  equations.
\newblock Comput. Methods Appl. Mech. Engrg. \textbf{209--212}, 184--196 (2012)

\bibitem{bause13hof}
Bause, M., Schwegler, K.: Higher order finite element approximation of systems
  of convection-diffusion-reaction equations with small diffusion.
\newblock Journal of Computational and Applied Mathematics \textbf{246}, 52--64
  (2013)

\bibitem{burger09pkd}
B\"{u}rger, R., Sep{\`u}lveda, M., Voitovich, T.: On the
  {P}roriol-{K}oornwinder-dubiner hierarchical orthogonal polynomial basis for
  the {DG-FEM} (2009)

\bibitem{cangiani13adg}
Cangiani, A., Georgoulis, E.H., Metcalfe, S.: Adaptive discontinuous {G}alerkin
  methods for nonstationary convection-diffusion problems.
\newblock IMA Journal of Numerical Analysis pp. 1--20 (2013)

\bibitem{castillo12pdg}
Castillo, P.: Performance of discontinuous {G}alerkin methods for elliptic
  {PDE}s.
\newblock SIAM J. Sci. Comput. \textbf{24}, 524--547 (2012)

\bibitem{pietro14rra}
Di~Pietro, D.A., Vohralik, M.: A review of recent advances in discretization
  methods, a posteriori error analysis, and adaptive algorithms for numerical
  modeling in geosciences.
\newblock Geosciences Numerical Methods \textbf{69}, 701--729 (2014)

\bibitem{dobrev08psi}
Dobrev, V.A., Lazarov, R.D., Zikatanov, L.T.: Preconditioning of symmetric
  interior penalty discontinuous {G}alerkin {FEM} for elliptic problems.
\newblock In: Domain Decomposition Methods in Science and Engineering XVII,
  \emph{Lecture Notes in Computer Science and Engineering}, vol.~60, pp.
  33--44. Springer (2008)

\bibitem{epshteyn07epp}
Epshteyn, Y., Rivi\`{e}re, B.: Estimation of penalty parameters for symmetric
  interior penalty {G}alerkin methods.
\newblock J. Comput. Appl. Math. \textbf{206}, 843--872 (2007)

\bibitem{houston02dhp}
Houston, P., Schwab, C., S\"{u}li, E.: Discontinuous hp-finite element methods
  for advection-diffusion-reaction problems.
\newblock SIAM J. Numer. Anal. \textbf{39}, 2133--2163 (2002)

\bibitem{klieber06ast}
Klieber, W., Rivi\`{e}re, B.: Adaptive simulations of two-phase flow by
  discontinuous {G}alerkin methods.
\newblock Computer Methods in Applied Mechanics and Engineering
  \textbf{196}(1--3), 404--419 (2006)

\bibitem{makridakis03era}
Makridakis, C., Nochetto, R.H.: Elliptic reconstruction and a posteriori error
  estimates for parabolic problems.
\newblock SIAM J. Numer. Anal. \textbf{41}(4), 1585--1594 (2003)

\bibitem{proft09dgm}
Proft, J., Riviere, B.: Discontinuous {G}alerkin methods for
  convection-diffusion equations for varying and vanishing diffusivity.
\newblock Int. J. Numer. Anal. Model \textbf{6}(4), 533--561 (2009)

\bibitem{riviere08dgm}
Rivi\`{e}re, B.: Discontinuous {G}alerkin methods for solving elliptic and
  parabolic equations, Theory and implementation.
\newblock SIAM (2008)

\bibitem{schotzau09rae}
Sch\"{o}tzau, D., Zhu, L.: A robust a-posteriori error estimator for
  discontinuous {G}alerkin methods for convection-diffusion equations.
\newblock Applied Numerical Mathematics \textbf{59}, 2236--2255 (2009)

\bibitem{sun05lna}
Sun, S., Wheeler, M.F.: L2({H}1) norm a posteriori error estimation for
  discontinuous {G}alerkin approximations of reactive transport problems.
\newblock Journal of Scientific Computing \textbf{22--23}(1--3), 501--530
  (2005)

\bibitem{tambue10eia}
Tambue, A., Lord, G.J., Geiger, S.: An exponential integrator for
  advection-dominated reactive transport in heterogeneous porous media.
\newblock Journal of Computational Physics \textbf{229}, 3957--3969 (2010)

\bibitem{uzunca14adg}
Uzunca, M., Karas\"{o}zen, B., Manguo\u{g}lu, M.: Adaptive discontinuous
  {G}alerkin methods for non-linear diffusion-convection-reaction equations.
\newblock Computers and Chemical Engineering \textbf{68}, 24--37 (2014)

\bibitem{slingerland14fls}
Van~Slingerland, P., Vuik, C.: Fast linear solver for diffusion problems with
  applications to pressure computation in layered domains.
\newblock Computational Geosciences pp. 1--14 (2014)

\end{thebibliography}

\end{document}